\magnification=\magstep1
\baselineskip=14pt

\def \Bbb{\bf }
\def \A{{\cal A}}
\def \B{{\cal B}}
\def \F{{\cal F}}
\def \G{{\cal G}}

\def \Z{{\Bbb Z}}
\def \Q{{\Bbb Q}}

\def \M{{\cal M}}

\def \Nat{{\Bbb N}}
\def \N{{\cal N}}
\def \R{{\Bbb R}}
\def \torus{{\Bbb T}}
\def \T{{\cal T}}

\def \C{{\cal C}}
\def \lra{\longrightarrow}
\def \lla{\longleftarrow}
\def \llra{-\!\!-\!\!\!\lra}

\def \sn{\bigskip\noindent}
\def \qed{\hfill $\bullet$}
\def \bla#1{\buildrel{#1}\over\lla }
\def \bra#1{\buildrel{#1}\over\lra }

\def \bbra#1{\buildrel{#1}\over\llra }
\def \im{\hbox{\rm Im}}

\def \norm#1{{| \! | #1 | \! |}}
\def \nin{\not\in}
\def \mod{\hbox{\ \rm mod\ }}
\def \rtimes{{\times}}
\def \H{{\cal H}}
\def \P{{\cal P}}
\def \I{{\cal I}}
\def \J{{\cal J}}
\def \JJ{{\cal P}}

\def \dirlim{\displaystyle\mathop{\rm lim\,}_{\rightarrow}}
\def \invlim{\displaystyle\mathop{\rm lim\,}_{\leftarrow}}
\def \cupij{\displaystyle\mathop{\cup}_{j\leq i}}
\def \del{{\partial}}
\def \W{{\cal W}}
\def \S{Section~}
\def \newpage{\vfill\eject}
\def \ho#1#2{H_{#1}(\Gamma,{#2})}
\def \bew{\sn{\bf Proof}}
\def \erz#1{\langle{#1}\rangle}
\def \frac#1#2{{{#1}\over {#2}}}
\def \rk{{\rm rk\, }}
\def \eos#1#2#3#4{\left\{\eqalign{{#1} &\quad
\hbox{ for }\quad {#2},\cr {#3} &\quad \hbox{ for }\quad  {#4}\cr}\right.}
\def \df{{\rm d}}
\def \d{{d}}
\def \dpe{{\dim V}}
\def \x{\vec{x}}

\def \ap{{\bf [AP]}}
\def \bks{{\bf [BKS]}}
\def \bjks{{\bf [BJKS]}}
\def \bsj{{\bf [BSJ]}}
\def \belone{{\bf [B2]}}
\def \beltwo{{\bf [B1]}}
\def \bcl{{\bf [BCL]}}
\def \br{{\bf [Br]}}
\def \db{{\bf [dB1]}}
\def \dbtwo{{\bf [dB2]}}
\def \con{{\bf [C]}}
\def \DuneauKatz{{\bf [DK]}}
\def \fh{{\bf [FH]}}
\def \fhk{{\bf [FHKI-III]}}
\def \FHKphys{{\bf [FHK]}}
\def \furst{{\bf [F]}}
\def \GaKe{{\bf [GK]}}
\def \gs{{\bf [GS]}}
\def \hm{{\bf [HM]}}
\def \hof{{\bf [H]}}
\def \kd{{\bf [KD]}}
\def \kelzero{{\bf [K1]}}
\def \kelone{{\bf [K2]}}
\def \keltwo{{\bf [K3]}}
\def \kn{{\bf [KrNe]}}
\def \Elser{{\bf [Els]}}
\def \ElserHenley{{\bf [EH]}}
\def \HenleyElser{{\bf [HE]}}
\def \KramerPapadopolos{{\bf [KrPa]}}
\def \lag{{\bf [La1]}}
\def \lagtwo{{\bf [La2]}}
\def \le{{\bf [Le]}}
\def \mack{{\bf [Mac]}}
\def \massey{{\bf [Mas]}}
\def \mrw{{\bf [MRW]}}
\def \okd{{\bf [OKD]}}
\def \pt{{\bf [PT]}}
\def \pat{{\bf [Pa]}}
\def \pen{{\bf [Pe]}}
\def \rad{{\bf [Ra]}}
\def \rei{{\bf [Rie]}}
\def \ren{{\bf [Ren]}}
\def \robone{{\bf [R1]}}
\def \robtwo{{\bf [R2]}}
\def \rudolf{{\bf [Ru]}}
\def \schtwo{{\bf [Sch]}} 
\def \sen{{\bf [Se]}}
\def \st{{\bf [ST]}}
\def \soc{{\bf [Soc]}}
\def \solone{{\bf [S1]}}
\def \sol{{\bf [S2]}}
\def \colone{{\bf [1]}}
\def \coltwo{{\bf [2]}}
\def \du{{\bf [Du]}}
\def \fortwo{{\bf [Fo]}}
\def \pr{{\bf [Pr]}}
\def \pss{{\bf [PSS]}}
\def \cbt{{\bf [T]}}
\def \BBG{{\bf [BBG]}}
\def \bhz{{\bf [BZH]}}
\def \kelput{{\bf [KePu]}}
\def \smith{{\bf [Sm]}}
\def \zobetz{{\bf [Z]}}
\def \gljj{{\bf [GLJJ]}}
\def \her{{\bf [He]}}
\def \gps{{\bf [GPS]}}
\def \wil{{\bf [W]}}
\def \schecht{{\bf [SBGC]}}
\def \tableofcontents#1{\halign{##\hfill&&\qquad\qquad\qquad\qquad
##\hfill\cr#1\crcr}}


\font \smallfont=cmr10 at 8pt
\font \chap=cmbx10 at 17pt
\font \sect=cmbx10 at 12pt

\footline={\ifnum\pageno=1\hss\tenrm\folio\hss\else\hfil\fi}
\headline={\ifnum\pageno=1\hfil\else
{\ifodd\pageno\rightheadline\else\leftheadline\fi}\fi}
\def\rightheadline{\tenrm\hfil{}\hfil\folio}
\def\leftheadline{\tenrm\folio\hfil{\smallfont
FORREST HUNTON KELLENDONK}\hfil} \voffset=2\baselineskip

\centerline{\chap TOPOLOGICAL INVARIANTS FOR}
\bigskip

\centerline{\chap PROJECTION METHOD PATTERNS}
\bigskip
\bigskip
\centerline{\sect Alan Forrest$^1$, John Hunton$^2$, Johannes Kellendonk$^3$}
\bigskip

\centerline{October 2000}\bigskip

\centerline{$^1$Department of Mathematics, }
\centerline{ National University of Ireland, }
\centerline{ Cork, Republic of Ireland}
\centerline{ Email: stmt8027@bureau.ucc.ie}
\bigskip

\centerline{$^2$The Department of Mathematics and Computer Science, }
\centerline{University of Leicester, University Road,  }
\centerline{Leicester, LE1 7RH, England}
\centerline{Email: jrh@mcs.le.ac.uk}
\bigskip

\centerline{$^3$Fachbereich Mathematik, Sekr.\ MA 7-2,}
\centerline{Technische Universit\"at Berlin, }
\centerline{10623 Berlin, Germany}
\centerline{Email: kellen@math.tu-berlin.de}

\newpage

\noindent {\sect TABLE OF CONTENTS}
\medskip

\tableofcontents{General Introduction & 03 \cr I Topological
Spaces and Dynamical Systems & 10\cr \quad 1 Introduction &  10\cr
\quad 2 The projection method and associated geometric
constructions & 11\cr \quad 3 Topological spaces for point
patterns &  15\cr \quad 4 Tilings and point patterns & 18\cr \quad
5 Comparing $\Pi_u$ and $\widetilde\Pi_u$  &  22\cr \quad 6
Calculating $M\widetilde P_u$ and $MP_u$ &  24\cr \quad 7
Comparing $MP_u$ with $M\widetilde P_u$ &  27\cr \quad 8 Examples
and Counter-examples & 30\cr \quad 9 The topology of the
continuous hull &  33\cr \quad 10 A Cantor $\Z^d$ Dynamical System
& 36\cr II Groupoids, $C^*$-algebras, and their invariants & 42\cr
\quad 1 Introduction &  42\cr \quad 2 Equivalence of Projection
Method pattern groupoids &  43\cr \quad 3 Continuous similarity of
Projection Method pattern groupoids &  49\cr \quad 4 Pattern
cohomology and K-theory &  52\cr \quad 5 Homological conditions
for self similarity &  54\cr III Approaches to calculation I:
cohomology for codimension 1 &  57\cr \quad 1 Introduction  &
57\cr \quad 2 Inverse limit acceptance domains & 57\cr \quad 3
Cohomology of the case d=N--1 &  58\cr IV Approaches to
calculation II: infinitely generated cohomology & 62\cr \quad 1
Introduction & 62\cr \quad 2 The canonical projection tiling &
62\cr \quad 3 Constructing $\C$-topes &  67\cr \quad 4 The
indecomposable case &  72\cr \quad 5 The decomposable case & 76\cr
\quad 6 Conditions for infinitely generated cohomology & 80\cr V
Approaches to calculation III: cohomology for small codimension &
84\cr \quad 1 Introduction  & 84\cr \quad 2 Set up and statement
of the results &  85\cr \quad 3 Complexes defined by the singular
spaces &  88\cr \quad 4 Group homology & 90\cr \quad 5 The
spectral sequences &  92\cr \quad 6 Example: Ammann-Kramer tilings
& 100\cr References & 103\cr}

\newpage

\headline={\ifnum\pageno=1\hfil\else
{\ifodd\pageno\rightheadline\else\leftheadline\fi}\fi}
\def\rightheadline{\tenrm\hfil{\smallfont
INTRODUCTION}\hfil\folio}
\def\leftheadline{\tenrm\folio\hfil{\smallfont
FORREST HUNTON KELLENDONK}\hfil} \voffset=2\baselineskip

\sn{\chap General Introduction}
\bigskip

Of the many examples of aperiodic tilings or aperiodic point sets
in Euclidean space found in recent years, two classes stand out as
particularly interesting and \ae sthetically pleasing. These are
the {\it substitution} tilings, tilings which are self-similar in
a rather strong sense described in \gs\ \robone\ \solone\ \ap, and
the tilings and patterns obtained by the method of cut and
projection from higher dimensional periodic sets described in \db\
\kn\ \kd. In this memoir we consider the second class. However,
some of the best studied and most physically useful examples of
aperiodic tilings, for example the Penrose tiling \pen\ and the
octagonal tiling (see \soc ), can be approached as examples of
either class. Therefore we study specially those tilings which are
in the overlap of these two classes, and examine some of their
necessary properties.

Tilings and patterns in Euclidean space can be compared by various
degrees of equivalence, drawn from considerations of geometry and
topology \gs.  Two tilings can be related by simple geometric
tranformations (shears or rotations), topological distortions
(bending edges), or by more radical adaptation (cutting tiles in
half, joining adjacent pairs etc). Moreover, point patterns can be
obtained from tilings in locally defined ways (say, by selecting
the centroids or the vertices of the tiles) and vice versa (say,
by the well-known Voronoi construction). Which definition of
equivalence is chosen is determined by the problem in hand.

In this paper, we adopt definitions of equivalence (pointed
conjugacy and  topological conjugacy, I.4.5) which allow us to
look, without loss of generality, at sets of uniformly isolated
points (point patterns) in Euclidean space. In fact these patterns
will typically have the Meyer property \lag\ (see I.4.5).
Therefore in this introduction, and often throughout the text, we
formulate our ideas and results in terms of point patterns and
keep classical tilings in mind as an implicit example.

The current rapid growth of interest in projection method patterns
started with the discovery  of material quasicrystals in 1984
\schecht, although these
patterns had been studied 
before then. Quasicrystaline material surprised the physical world
by showing sharp Bragg peaks under X-ray scattering, a phenomenon
usually associated only with periodic crystals. Projection method
patterns share this unusual property and in recent studies they
have become the preferred model of material quasicrystals \colone\
\coltwo . This model is not without criticism, see e.g.\ \lagtwo.

Whatever the physical significance of the projection method
construction, it also has great mathematical appeal in itself: it
is elementary and geometric and, once the acceptance domain and
the dimensions of the spaces used in the construction are chosen,
has a finite number of degrees of freedom. The projection method
is also a natural generalization of low dimensional examples such
as Sturmian sequences \hm\ which have strong links with classical
diophantine approximation.

\sn{\bf General approach of this book} In common with many papers
on the topology of tilings, we are motivated by the physical
applications and so are interested in the properties of an
individual quasicrystal or pattern in Euclidean space. The
topological invariants of the title refer not to the topological
arrangement of the particular configuration as a subset of
Euclidean space, but rather to an algebraic object (graded group,
vector space {\it etc}.) associated to the pattern, and which in
some way captures its geometric properties. It is defined in
various equivalent ways as a classical topological invariant
applied to a space constructed out of the pattern. There are two
choices of space to which to apply the invariant, the one
$C^*$-algebraic, the other dynamical, and these reflect the two
main approaches to this subject, one starting with the
construction of an operator algebra and  the other with a
topological space with $\R^d$ action.

The first approach, which has the benefit of being closer to
physics and which thus provides a clear motivation for the
topology, can be summarized as follows. Suppose that the point set
$\T$ represents the positions of atoms in a material, like a
quasicrystal. It then provides a discrete model for the
configuration space of particles moving in the material, like
electrons or phonons. Observables for these  particle systems,
like energy, are, in the absence of external forces like a
magnetic field, functions of partial translations. Here a partial
translation is an operator on the Hilbert space of square summable
functions on $\T$ which is a translation operator from one point
of $\T$ to another combined with a range projection which depends
only on the neighbouring configuration of that point. The
appearance of that range projection is directly related to the
locality of interaction. The norm closure $\A_\T$ of the algebra
generated by partial translations can be regarded as the
$C^*$-algebra of observables. The topology we are interested in is
the non-commutative topology of the $C^*$-algebra of observables
and the topological invariants of $\T$ are the invariants of
$\A_\T$. In particular, we shall be interested in the $K$-theory
of $\A_\T$. Its $K_0$-group has direct relevance to physics
through the gap-labelling. In fact, Jean Bellissard's
$K$-theoretic formulation of the gap-labelling stands at the
beginning of this approach \beltwo.

The second way of looking at the topological invariants of a
discrete point set $\T$ begins with the construction of the
continuous hull of $\T$. There are various ways of defining a
metric on a set of patterns through comparison of their local
configurations. Broadly speaking, two patterns are deemed close if
they coincide on a large window around the origin $0\in \R^d$ up
to a small discrepancy. It is the way the allowed discrepancy is
quantified which leads to different metric topologies and we
choose here one which has the strongest compactness properties,
though we have no intrinsic motivation for this. The continuous
hull of $\T$ is the closure, $M\T$, of the set of translates of
$\T$ with respect to this metric. We use the notation $M\T$
because it is essentially a mapping torus construction for a
generalized discrete dynamical system: $\R^d$ acts on $M\T$ by
translation and the set $\Omega_\T$ of all elements of $ M\T$
which (as point sets) contain $0$ forms an abstract transversal
called the discrete hull. If $d$ were $1$ then $\Omega_\T$ would
give rise to a Poincar\'e section, the intersection points of the
flow line of the action of $\R$ with $\Omega_\T$ defining an orbit
of a $\Z$ action in $\Omega_\T$, and $\M\T$ would be the mapping
torus of that discrete dynamical system $(\Omega_\T,\Z)$. For
larger $d$ one cannot expect to get a $\Z^d$ action on $\Omega_\T$
in a similar way but finds instead a generalized discrete
dynamical system which can be summarized in a groupoid $\G\T$
whose unit space is the discrete hull $\Omega_\T$.
Topological invariants for $\T$ are therefore the topological
invariants of $M\T$ and of $\G\T$ and we shall be interested in
particular in their cohomologies. We define the cohomology of
$\T$  to be that of $\G\T$. Under a finite type condition, namely
that for any given $r$ there are only a finite number of
translational congruence classes of subsets which fit inside a
window of diameter $r$, the algebra $\A_\T$ sketched above is
isomorphic to the groupoid $C^*$-algebra of $\G\T$. This links the
two approaches.

Having outlined the general philosophy we hasten to remark that we
will not explain all its aspects in the main text. In particular,
we have nothing to say there about the physical aspects of the
theory and the description of the algebra of observables,
referring here the reader to \bhz\ \kelput, or to the more
original literature \beltwo\ \belone\ \kelzero. Instead, our aim
in this memoir is to discuss and compare the different commutative
and non-commutative invariants, to demonstrate their applicability
as providing obstructions to a tiling arising as a substitution,
and finally to provide a practical method for computing them; this
we illustrate with a number of examples, including that of the
Ammann-Kramer (\lq 3 dimensional Penrose\rq) tiling. Broadly
speaking, one of the perspectives of this memoir is that
non-commutative invariants for projection point patterns can be
successfully computed by working with suitable commutative
analogues.

\sn{\bf The subject of this book} We work principally with the
special class of point sets (possibly with some decoration)
obtained by cut and projection from the integer lattice $\Z^N$
which is generated by an orthonormal base of $\R^N$. Reserving
detail and elaborations for later, we call a {\it projection
method pattern\/} $\T$ on $E=\R^d$ a pattern of points (or a
finite decoration of it) given by the orthogonal projection onto
$E$ of points in a strip $(K\times E)\cap\Z^N\subset\R^N$, where
$E$ is a subspace of $\R^N$ and $K\times E$ is the so-called {\it
acceptance strip}, a fattening of $E$ in $\R^N$ defined by some
suitably chosen region $K$ in the orthogonal complement
$E^{\perp}$ of $E$ in $\R^N$. The pattern $\T$ thus depends on the
dimension $N$, the positioning of $E$ in $\R^N$ and the shape of
the {\it acceptance domain\/} $K$. When this construction was
first made \db\ \kd\ the domain $K$ was taken to be the projected
image onto $E^{\perp}$ of the unit cube in $\R^N$ and this choice
gives rise to the so-called {\it canonical projection method
patterns}, but for the first three chapters we allow $K$ to be any
compact subset of $E^{\perp}$ which is the closure of its interior
(so, with possibly even fractal boundary, a case of current
physical interest \bks\ \smith\ \zobetz\ \gljj ).

It is then very natural to consider not only $\T$ but also all
point patterns which are obtained in the same way but with $\Z^N$
repositioned by some vector $u\in\R^N$, {\it i.e.}, $\Z^N$
replaced by $\Z^N+u$. Completing certain subsets of positioning
vectors $u$ with respect to an appropriate pseudo-metric gives us
the continuous hull $M\T$. This analysis shows in particular that
$M\T$ contains another transversal $X_\T$ which gives rise to $d$
independant commuting $\Z$ actions and hence to a genuine discrete
dynamical system $(X_\T,\Z^d)$ whose mapping torus is also $M\T$.
This is a key point in relating the $K$-theory of $\A_\T$ with the
cohomology of $\T$; in the process, the latter is also identified
with the \v Cech cohomology of $M\T$ and with the group cohomology
of $\Z^d$ with coefficients in the continuous integer valued
functions over $X_\T$.

The space $X_{\T}$ arises in another way. Let $V$ be a connected
component of the euclidian closure of $\pi^\perp(\Z^N)$, where
$\pi^\perp$ denotes the orthoprojection onto $E^\perp$. We first
disconnect $V$ along the boundaries of all
$\pi^\perp(\Z^N)$-translates of $K$ (we speak loosely here, but
make the idea precise in Chapter I). Then $X_\T$ can be understood
as a compact quotient of the disconnected $V$ with respect to a
proper isometric free abelian group action. In the case of
canonical projection method patterns, on which we concentrate in
the last two chapters, the boundaries which disconnect $V$ are
affine subspaces and so define a directed system of locally finite
CW decompositions of Euclidean space. With respect to this
CW-complex, the integer valued functions over $X_\T$ appear as
continuous chains in the limit. This makes the group cohomology of
the dynamical action of $\Z^d$ on $X_{\T}$ accessible through the
standard machinery (exact sequences and spectral sequences) of
algebraic topology.

As mentioned, the interest in physics in the non-commutative
topology of tilings and point sets is based on the observation
that $\A_\T$ is the $C^*$-algebra of observables for particles
moving in $\T$. In particular, any Hamilton operator which
describes this motion has the property that its spectral
projections on energy intervals whose boundaries lie in gaps of
the spectrum belong to $\A_\T$ as well and thus define elements of
$K_0(\A_\T)$. Therefore, the ordered $K_0$-group (or its image on
a tracial state) may serve to \lq count\rq\ (or label) the
possible gaps in the spectrum the Hamilton operator \belone\ \BBG\
\kelzero. One of the main results of this memoir is the
determination in Chapter V of closed formulae for the ranks of the
$K$-groups corresponding to canonical projection method patterns
with small codimension (as one calls the dimension of $V$). These
formul\ae\ apply to all common tilings including the Penrose
tilings, the octagonal tilings and three dimensional icosahedral
tilings. Unfortunately our method does not as yet give full
information on the order of $K_0$ or the image on a tracial state.

Further important results of this memoir concern the structure of
a $K$-group of a canonical projection method pattern. We find that
its $K_0$-group is generically infinitely generated. But when the
rank of its rationalization is finite then it has to be free
abelian. We observe in Chapters III and IV that both properties
are obstruction to some kinds of self-similarity. More precisely,
infinitely generated rationalized cohomology rules out that the
tiling is a substitution tiling. On the other hand, if we know
already that the tiling is substitutional then its $K$-group must
be free abelian for it to be a canonical projection method tiling.

No projection method pattern is known to us which has both
infinitely generated cohomology and also allows for local matching
rules in the sense of \le. Furthermore, all projection method
patterns which are used to model quasicrystals seem to have a
finitely generated $K_0$-group. We cannot offer yet an
interpretation of the fact that some patterns produce only
finitely many generators for their cohomology whereas others do
not, but, if understood, we hope it could lead to a criterion to
single out a subset of patterns relevant for quasicrystal physics
from the vast set of patterns which may be obtained from the
projection method.

We have mentioned above the motivation from physics to study the
topological theory of point sets or tilings. The theory is also
also of great interest for the theory of topological dynamical
systems, since in $d=1$ dimensions the dynamical systems mentioned
above have attracted a lot of attention. In \gps\ the meaning of
the non-commutative invariants for the one dimensional case has
been analysed in full detail. Furthermore, substitution tilings
give rise to hyperbolic $\Z$-actions with expanding attractors
(hyperbolic attractors whose topological dimensions are that of
their expanding direction) \ap\ a subject of great interest
followed up recently by Williams \wil\ who conjectures that
continuous hulls of substitution tilings (called tiling spaces in
\wil) are fiber bundles over tori with the Cantor set as fiber. We
have not put emphasis on this question but it may be easily
concluded from our analysis of Section I.10 that the continuous
hulls of projection method tilings are always Cantor set fiber
bundles over tori (although these tilings are rarely
substitutional and therefore carry in general no obvious
hyperbolic $\Z$-action). Anderson and Putnam \ap\ and one of the
authors \kelone\ have employed the substitution of the tiling to
calculate topological invariants of it.
\bigskip

\sn{\bf Organization of the book} The order of material in this
memoir is as follows. In Chapter I we define and describe the
various dynamical systems mentioned above, and examine their
topological relationships. These are compared with the pattern
groupoid and its associated $C^*$ algebra in Chapter II where
these latter objects are introduced. Also in Chapter II we set up
and prove the equivalence of our various topological invariants;
we end the chapter by demonstrating how these invariants provide
an obstruction to a pattern being self-similar.

The remaining chapters offer three illustrations of the
computability of these invariants. In Chapter III we give a
complete calculation for all \lq codimension one\rq\ projection
patterns -- patterns arising from the projection of slices of
$\Z^{d+1}$ to $\R^d$ for more or less arbitrary acceptance
domains. In Chapter IV we give descriptions of the invariants for
generic projection patterns arising from arbitrary projections
$\Z^N$ to $\R^d$ but with canonical acceptance domain. Here,
applying the result at the end of Chapter II, we prove the result
mentioned above that almost all canonical projection method
patterns have infinitely generated cohomology and so fail to be
substitution tilings. In Chapter V we examine the case of
canonical projection method patterns with finitely generated
cohomology, such as would arise from a substitution system. We
develop a systematic approach to the calculation of these
invariants and use this to produce closed formul\ae\ for the
cohomology and $K$-theory of projection patterns of codimension 1,
2 and 3: in principle the procedure can be iterated to higher
codimensions indefinitely, though in practice the formul\ae\ would
soon become tiresome. Some parameters of these formulae allow for
a simple description in arbitrary codimension, as e.g.\ the Euler
characterisitc (V.2.8). We end with a short description of the
results for the Ammann-Kramer tiling.

There is a separate introduction to each chapter where relevant
classical work is recalled and where the individual sections are
described roughly. We adopt the following system for
crossreferences. The definitions, theorems etc.\ of the same
chapter are cited e.g.\ as Def.~2.1 or simply 2.1. The
definitions, theorems etc.\ of the another chapter are cited e.g.\
as Def.~II.2.1 or simply II.2.1.

\sn{\bf A note on the writing of this book}  Originally this
memoir was conceived by the three authors as a series of papers
leading to the results now in Chapters IV and V, aiming to found a
calculus for projection method tiling cohomology. These papers are
currently available as a preprint-series {\it Projection
Quasicrystals I-III} \fhk\ covering most of the results in this
memoir. The authors' collaboration on this project started in 1997
and, given the importance of the subject and time it has taken to
bring the material to its current state, it is inevitable that
some results written here have appeared elsewhere in the
literature during the course of our research. We wish to
acknowledge these independent developments here, although we will
refer to them again as usual in the body of the text.

The general result of Chapter I, that the tiling mapping torus is
also a discrete dynamical mapping torus, and that the relevant
dynamics is an almost 1-1 extension of a rotation on a torus, has
been known with varying degrees of precision and generality for
some time and we mention the historical developments in the
introduction to Chapter I. Our approach constructs a large
topological space from which the pattern dynamical system is
formed by a quotient and so we follow most closely the idea
pioneered by Le \le\ for the case of canonical projection tilings.
The ``Cantorization'' of Euclidean space by corners or cuts, as
described by Le and others (see \le\ \hof\ etal.), is produced in
our general topological context in sections I.3, I.4 and I.9. In
this, we share the ground with Schlottmann \schtwo\  and Herrmann
\her\ who have recently established the results of Chapter I in
such (and even greater) generality, Schlottmann in order to
generalize results of Hof and describe the unique ergodicity of
the underlying dynamical systems and Herrmann to draw a connection
between codimension $1$ projection patterns and Denjoy
homeomorphisms of the circle. We mention this relation at the end
of chapter~III.

Bellissard, Contensou and Legrand \bcl\ compare the $C^*$-algebra
of a dynamical groupoid with a $C^*$-algebra of operators defined
on a class of tilings obtained by projection, the general theme of
Chapter II. Using a Rosenberg Shochet spectral sequence, they also
establish, for $2$-dimensional canonical projection tilings, an
equation of dynamical cohomology and $C^*$ $K$-theory in that
case. It is the first algebraic topological approach to projection
method tiling $K$-theory found in the literature. We note,
however, that the groupoid they consider is not always the same as
the tiling groupoid we consider, nor do the dynamical systems
always agree; the Penrose tiling is a case in point, where we find
that $K_0$ of the spaces considered in \bcl\ is $\Z^{\infty}$. The
difference may be found in the fact that we consider a given
projection method tiling or pattern and its translates, while they
consider a larger set of tilings, two elements of which may
sometimes be unrelated by approximation and translation parallel
to the projection plane.

\sn{\bf Acknowledgements.}
We thank F.~G\"ahler for verifying and assisting our results (V.6)
on the computer and for various useful comments
and M.~Baake for reading parts of the manuscript.
The third author would like to thank
J.~Bellissard for numerous discussions and constant support.

The collaboration of the first two authors was
initiated by the William Gordon Seggie Brown Fellowship at The
University of Edinburgh, Scotland, and received continuing
support from a Collaborative
Travel Grant from the British Council and the Research Council of Norway
with the generous assistance of The University of Leicester, England, and
the EU Network ``Non-commutative Geometry'' at NTNU Trondheim, Norway. The
collaboration of the first and third authors was supported by the
Sonderforschungsbereich 288, ``Differentialgeometrie und Quantenphysik" at
TU Berlin, Germany, and by
the EU Network and NTNU Trondheim. The third
author is supported by the Sfb288 at TU Berlin. All
three authors are most grateful for the financial help received from these
various sources.

\newpage

\headline={\ifnum\pageno=1\hfil\else
{\ifodd\pageno\rightheadline\else\leftheadline\fi}\fi}
\def\rightheadline{\tenrm\hfil{\smallfont
I TOPOLOGICAL SPACES AND DYNAMICAL SYSTEMS}\hfil\folio}
\def\leftheadline{\tenrm\folio\hfil{\smallfont
FORREST HUNTON KELLENDONK}\hfil} \voffset=2\baselineskip

\sn{\chap I Topological Spaces and Dynamical Systems}
\bigskip

\sn{\sect 1 Introduction}

\sn
In this chapter our broad goal is to study the topology and associated
dynamics of projection method patterns, while imposing only few restrictions
on the freedom of the construction. From a specific set of projection data we
define and examine a number of spaces and dynamical systems and their
relationships; from these constructions, in later chapters, we set up our
invariants, defined via various topological and dynamical cohomology
theories. In Chapter II we shall also compare the commutative spaces of this
chapter with the non-commutative spaces considered by other authors in {\it
e.g.\/} \belone\ \bcl\ \ap\ \kelone.

Given a subspace, $E$, acceptance domain, $K$, and a positioning parameter
$u$, we  distinguish two particular $\R^d$ dynamical systems constructed by
the projection method, $(MP_u, \R^d)$ and $(M\widetilde P_u, \R^d)$, the first
automatically a factor of the second. This allows us to define a {\it
projection method pattern} (with data $(E,K,u)$) as a pattern, $\T$, whose
dynamical system, $(M\T, \R^d)$, is intermediate to these two extreme
systems. Sections 3 to 9 of this chapter provide a complete description of the
spaces and the extension $M\widetilde P_u \lra MP_u$, showing, under further
weak assumptions on the acceptance domain, that it is a finite isometric
extension.  In section 7 we conclude that this restricts $(M\T, \R^d)$ to one
of a finite number of possibilities, and that any projection method pattern is
a finite decoration of its corresponding point pattern $P_u$ (2.1). The
essential definitions are to be found in
Sections~2 and 4.

In section 10 we describe yet another dynamical system connected
with a projection method pattern, this time a $\Z^d$ action on a
Cantor set $X$, whose mapping torus is the space of the pattern
dynamical system. It will be  this dynamical system that, in
chapters 3, 4 and 5, will allow the easiest computation and
discussion of the behaviour of our invariants. For the canonical
case with $E^{\perp} \cap \Z^N = 0$ this is the same system as
that constructed in \bcl .

All the dynamical systems produced in this memoir are almost 1-1
extensions of an action of $\R^d$ or $\Z^d$ by rotations (Def.\
2.15) on a torus (or torus  extended by a finite abelian group).
In each case the dimension of the torus and the generators of the
action can be computed explicitly. This gives a clear picture of
the orbits of non-singular points in the pattern dynamical system.
A precursor to our description of the pattern dynamical system can
be found in the work of Robinson \robtwo , who examined the
dynamical system of the Penrose tiling and showed that it is an
almost 1-1 extension of a minimal $\R^2$ action by rotation on a
$4$-torus. Although Robinson used quite special properties of the
tiling, Hof \hof\ has noted that these techniques are
generalizable without being specific about the extent of the
generalization.

Our approach is quite different from that of Robinson and, by
constructing a larger topological space from which the pattern
dynamical system is formed by a quotient, we follow most closely
the approach pioneered by Le \le\ as noted in the General
Introduction. The care taken here in the topological foundations
seems necessary for further progress and to allow general
acceptance domains. Even in the canonical case, Corollary 7.2 and
Proposition 8.4 of this chapter, for example, require this
precision despite being direct generalizations of Theorem 3.8 in
\le . Also, as mentioned in the General Introduction, many of the
results of this Chapter are to be found independently in \schtwo .

\bigskip

\sn{\sect 2 The projection method and associated geometric constructions}

\sn We use the construction of point patterns and tilings given in
Chapters 2 and 5 of Senechal's monograph \sen\ throughout this
paper, adding some assumptions on the acceptance domain in the
following definitions.

\sn{\bf Definitions 2.1} Consider the lattice $\Z^N$ sitting in
standard position inside $\R^N$ (i.e.\ it is generated by an
orthonormal basis of $\R^N$). Suppose that $E$ is a $d$
dimensional subspace of $\R^N$ and $E^{\perp}$ its
orthocomplement. For the time being we shall make no assumptions
about the position of either of these planes.

Let $\pi$ be the projection onto $E$ and $\pi^{\perp}$ the
projection onto $E^{\perp}$.

Let $Q= \overline{ E+\Z^N}$ (Euclidean closure). This is a closed
subgroup of $\R^N$.

Let $K$ be a compact subset of $E^{\perp}$ which is the closure if
its interior (which we write $Int K$) in $E^{\perp}$. Thus the
boundary of $K$  in $E^{\perp}$ is compact and nowhere dense. Let
$\Sigma = K + E$, a subset of $\R^N$ sometimes refered to as {\it
the strip with acceptance domain $K$}.

A point $v \in \R^N$ is said to be {\it non-singular} if the boundary,
$\partial\Sigma$, of $\Sigma$ does not intersect $\Z^N+v$. We write $NS$
for the set of non-singular points in $\R^N$. These points are also called {\it
regular} in the literature.

Let $\widetilde P_v = \Sigma \cap(\Z^N +v)$, the {\it strip point pattern}.

Define $P_v = \pi(\widetilde P_v)$, a subset of $E$ called the {\it
projection point pattern}.

\sn  In what follows we assume $E$ and $K$ are fixed and
suppress mention of them as a subscript or argument.

\sn{\bf Lemma 2.2} {\it With the notation above,

i/ $NS$ is a dense $G_{\delta}$ subset of $\R^{N}$
invariant under translation by $E$.

ii/ If $u \in NS$, then $NS \cap (Q+u)$ is dense in $Q+u$.

iii/ If $u \in NS$ and $F$ is a vector subspace of $\R^N$ complementary to
$E$, then $NS \cap (Q+u) \cap F$ is dense in $(Q+u) \cap F$.}

\sn{\bf Proof} i/ Note that $\R^N \setminus NS$ is a translate of
the set $\cup_{v \in \Z^N} (\partial K + E + v)$ (where the
boundary is taken in  $E^{\perp}$) and our conditions on $K$
complete the proof.

ii/ $NS \cap (Q+u) \supset E + \Z^N + u$.

iii/ $\overline{(E + \Z^N + u) \cap F} = (Q+u) \cap F$. \qed

\sn{\bf Remark 2.3} The condition on the acceptance domain $K$ is
a  topological version of the condition of \hof . We note that our
conditions include the examples of acceptance domains with fractal
boundaries which have recently interested quasicrystalographers
\bks\ \smith\ \zobetz\ \gljj.

In the original construction \db\ \kd\  $K=\pi^{\perp}([0,1]^N)$.
We call this the {\it canonical acceptance domain} and we reserve
the name {\it canonical projection method pattern} for the
patterns $P_u$ produced from this acceptance domain. Sometimes
this is shortened to {\it canonical pattern} for convenience.

This is closely related to the {\it canonical projection tiling},
defined by \okd\ formed by a canonical acceptance domain, $u \in
NS$ and projecting onto $E$ those $d$-dimensional faces of the
lattice $\Z^N +u$ which are contained entirely in $\Sigma$. We
write this tiling $\T_u$.

\sn The following notation and technical lemma makes easier some
calculations in future sections.

\sn{\bf Definition 2.4} If $X$ is a subspace of $Y$, both
topological spaces, and $A \subset X$, then we write $Int_XA$ to
mean the interior of $A$ in the subspace topology of $X$.

Likewise we write $\partial_XA$ for the boundary of $A$ taken in
the subspace topology of $X$.

\sn{\bf Lemma 2.5} {\it a/ If $u \in NS$, then $ (Q+u)\cap Int K =
Int_{(Q+u)  \cap E^{\perp}}((Q+u) \cap K)$ and $ (Q+u)\cap
\partial_{E^{\perp}}K = \partial_{(Q+u) \cap E^{\perp}}((Q+u) \cap
K)$.

b/ If $u \in NS$, then $((Q+u) \cap E^{\perp}) \setminus NS =
\partial_{(Q+u) \cap E^{\perp}}((Q+u) \cap K) +
\pi^{\perp}(\Z^N)$.}

\sn{\bf Proof} a/ To show both facts, it is enough to show that
$(\partial_{E^{\perp}} K) \cap (Q+u)$ has no interior as a
subspace of $(Q+u) \cap E^{\perp}$.

Suppose otherwise and that $U$ is an open subset of  $\partial K
\cap (Q+u)$ in $(Q+u) \cap E^{\perp}$. By the density of
$\pi^{\perp}(\Z^N)$ in $Q \cap E^{\perp}$, we find $v \in \Z^N$
such that $u \in U+\pi^{\perp}(v)$. But this implies that $u \in
\partial K + \pi^{\perp}(v)$ and so $u \nin NS$ - a contradiction.

b/ By definition the left-hand side of the equation to be proved
is equal to $(\partial_{E^{\perp}} K + \pi^{\perp}(\Z^N)) \cap
(Q+u)$ which equals $(\partial_{E^{\perp}} K \cap (Q+u)) +
\pi^{\perp}(\Z^N)$ since $\pi^{\perp}(\Z^N)$ is dense in $Q \cap
E^{\perp}$. By part a/ therefore we obtain the right-hand side of
the equation. \qed

\sn{\bf Condition 2.6} We exclude immediately the case  $(Q+u)
\cap Int K = \emptyset$ since, when $u \in NS$, this is equivalent
to $P_u = \emptyset$.

\sn{\bf Examples 2.7} We note the parameters of two well-studied examples,
both with canonical acceptance domain (2.3).

The {\it octagonal tiling} \soc\ has $N=4$ and $d=2$, where $E$ is a
vector subspace of $\R^4$ invariant under the action of the linear map
which maps orthonormal basis vectors $e_1 \mapsto e_2$, $e_2 \mapsto
e_3$, $e_3 \mapsto e_4$, $e_4 \mapsto -e_1$. Its orthocomplement,
$E^{\perp}$, is the other invariant subspace. Here $Q = \R^4$ and so many
of the distinctions made in subsequent sections are irrelevant to this
example.

The {\it Penrose tiling} \pen\ \db\ 
has $N=5$ and $d=2$ (although we note that there is an elegant
construction using the root lattice of $A_4$ in $\R^4$ \bjks). The
linear map which maps $e_i \mapsto e_{i+1}$ (indexed modulo $5$)
has two $2$ dimensional and one $1$ dimensional invariant
subspaces. Of the first two subspaces, one is chosen as $E$ and
the other we name $V$. Then in fact $Q = E \oplus V \oplus
\widetilde \Delta$, where $\widetilde \Delta = {1 \over 5}(e_1 +
e_2 + e_3 + e_4 + e_5)\Z$, and $Q$ is therefore a proper subset of
$\R^5$, a fact which allows the construction of {\it generalized
Penrose tilings} using a parameter $u \in NS \setminus Q$.

Note that we speak of tilings and yet only consider point
patterns.  In both examples, the projection tiling \okd\ is
conjugate to both the corresponding strip point pattern and
projection point pattern, a fact proved in greater generality in
section 8.

\sn We develop these geometric ideas in the following lemmas. The next
is Theorem 2.3 from \sen .

\sn{\bf Theorem 2.8}  {\it Suppose that $\Z^N$ is in standard
position in $\R^N$ and suppose that $\phi\colon \R^N \lra \R^n$ is
a surjective linear map. Then there is a direct sum decomposition
$\R^n = V \oplus W$ into real vector subspaces such that
$\phi(\Z^N) \cap V$ is dense in $V$, $\phi(\Z^N) \cap W$ is
discrete and $ \phi(\Z^N) = (V \cap \phi(\Z^N)) +  (W \cap
\phi(\Z^N))$. \qed}

\sn We proceed with the following refinement of Proposition 2.15 of \sen .

\sn{\bf Lemma 2.9} {\it Suppose that $\Z^N$ is in standard position in
$\R^N$ and suppose that $\phi\colon  \R^N \lra F$ is an orthogonal
projection onto $F$ a subspace of $\R^N$. With the decomposition of $F$
implied by Theorem 2.8, $(F \cap \Z^N) + (V \cap \phi(\Z^N)) \subset
\phi(\Z^N)$ as a finite index subgroup.

Also, the lattice dimension of $F \cap \Z^N$ equals $\dim F - \dim
V$ and the real vector subspace generated by  $F\cap\Z^N$ is
orthogonal to $V$.}

\sn{\bf Proof} Suppose that $U$ is the real linear span of $\Delta=F \cap
\Z^N$. Since $\Delta$ is discrete, the lattice dimension of
$\Delta$ equals the real space dimension of $U$.

 The argument of the proof of Proposition 2.15 in \sen\ shows that each
element of $F\cap \Z^N$ is orthogonal to $V$. Therefore we have
$\dim_{\R}(U) \leq \dim_{\R}(F) - \dim_{\R}(V)$ immediately.

Consider the rational vector space $\Q^N$, contained in $\R^N$ and
containing $\Z^N$, both in canonical position. Let $U'$ be the rational
span of $\Delta$ and note that $U' = U \cap \Q^N$ and that $\dim_{\Q}(U') =
\dim_{\R}(U)$. Let $U'^{\perp}$ be the orthocomplement of $U'$ with respect
to the standard inner product in $\Q^N$ so that, by simple rational vector
space arguments, $\Q^N = U' \oplus U'^{\perp}$. Thus $(U' \cap \Z^N) +
(U'^{\perp} \cap \Z^N)$ forms a discrete lattice of dimension $N$.

Extending to the real span, we deduce that $(U \cap \Z^N)+ (U^{\perp} \cap
\Z^N)$ is a discrete sublattice of $\Z^N$ of dimension $N$, hence a
subgroup of finite index. Also the lattice dimension of $U \cap \Z^N$ and
$U^{\perp} \cap \Z^N$ are equal to $\dim_{\R}(U)$ and $\dim_{\R}
(U^{\perp})$ respectively.

Let $L=(U^{\perp} \cap \Z^N)$ be considered as a sublattice of $U^{\perp}$.
It is integral (with respect to the restriction of the inner product on
$\R^N$) and of full dimension. The projection $\phi$ restricts to an
orthogonal projection $U^{\perp} \lra U^{\perp} \cap F$ and, by
construction, $U^{\perp} \cap F \cap L = 0$. Therefore Proposition 2.15 of
\sen\ applies to show that $\phi(L)$ is dense in $U^{\perp} \cap F$ and
that $\phi$ is 1-1 on $L$.

However $\phi(L) \subset \phi(\Z^N)$ and so, by the characterisation of
Theorem 2.2, we deduce that $U^{\perp} \cap F \subset V$. However, since
$U^{\perp} \supset V$, we have  $U^{\perp} \cap F = V$.

We have $U \cap \Z^N = F \cap \Z^N$ and $\phi(U^{\perp}\cap \Z^N) =
\phi(\Z^N) \cap V$ automatically. Therefore $(\phi(\Z^N) \cap V) + (F \cap
\Z^N) = \phi((U^{\perp} \cap \Z^N) + (U \cap \Z^N))$. As proved above, this
latter set is the image of a finite index subgroup of the domain, $\Z^N$,
and therefore it is a finite index subgroup of the image $\phi(\Z^N)$ as
required.

The remaining properties follow quickly from the details above. \qed

\sn{\bf Definition 2.10} Let $\Delta = E^{\perp} \cap \Z^N$ and
$\widetilde \Delta = U \cap \overline{\pi^{\perp}(\Z^N)}$ where
$U$ is the  real vector space generated by $\Delta$.

\sn Note that the discrete group $\Delta$ defined here is not the real
vector space $\Delta (E)$ defined in \le , but it is a cocompact sublattice
and so the dimensions are equal.

\sn{\bf Corollary 2.11} {\it With the notation of Thm.~2.8 and
taking $\phi=\pi^\perp$, $\overline{\pi^{\perp}(\Z^N)} = V \oplus
\widetilde\Delta$ and $Q = E \oplus V \oplus \tilde\Delta$ are
orthogonal direct sums. Moreover, $\Delta$  is a subgroup of
$\widetilde\Delta$ with finite index. \qed}

\sn{\bf Example 2.12} For example the octagonal tiling has $\Delta = 0$ and
the Penrose tiling has $\Delta = (e_1+ e_2 + e_3 + e_4 + e_5)\Z$, a
subgroup of index $5$ in $\widetilde \Delta$.

\sn And finally a general result about isometric extensions of dynamical
systems.

\sn{\bf Definition 2.13} Suppose that $\rho\colon  (X,G) \lra (Y,G)$ is a
factor map of topological dynamical systems with group, $G$, action. If
every fibre $\rho^{-1}(y)$ has the same finite cardinality, $n$, then we
say that
$(X,G)$ is an $n$-to-1 extension.

\sn The structure of such extensions, a special case of isometric
extensions, is well-known \furst.

\sn{\bf Lemma 2.14} {\it Suppose that $\rho\colon  (X,G) \lra (Y,G)$ is an
$n$-to-1 extension and that $(X,G)$ is minimal. Suppose further that there
is an abelian group $H$ which acts continuously on $X$, commutes with
the $G$ action, preserves $\rho$ fibres and acts transitively on each
fibre. If $(X,G) \bra {\rho'} (Z,G) \bra {\rho''} (Y,G)$ is an intermediate
factor, then $(Z,G)$ is an $m$-to-1 extension where $m$ divides $n$, and we
can find a subgroup, $H'$ of $H$, so that

i/ $H/H'$ acts continuously on $Z$, commutes with the $G$ action, preserves
$\rho''$ fibres and acts transitively on each fibre and

ii/ $H'$ acts on $X$ as a subaction of $H$, preserving $\rho'$ fibres and
acting transitively on each fibre.}

\sn{\bf Proof} Given $h \in H$, consider $X_h = \{x\ |\  \rho'(x) =
\rho'(hx)\}$ which is a closed $G$-invariant subset of $X$. Therefore, by
minimality, $X_h = \emptyset$ or $X$. Let $H' = \{h \in H \ |\  X_h = X\}$
which can be checked is a subgroup of $H$. The properties claimed follow
quickly. \qed

\sn{\bf Definitions 2.15} We will call an extension which obeys
the conditions  of Lemma 2.14 a {\it finite isometric extension}.

An {\it almost 1-1 extension} of topological dynamical systems
$\rho \colon  (X,G)\lra (Y,G)$ is one in which the set
$\rho^{-1}(y)$ is a singleton for a dense $G_{\delta}$ of $y \in
Y$. In the case of minimal actions, it is sufficient to find just
one point $y \in Y$ for which $\rho^{-1}(y)$ is a singleton.

We say that an abelian topological group, $G$, acting on a compact
abelian topological group, $Z$ say, acts {\it by rotation} if
there is a continuous  group homomorphism, $\psi : G \lra Z$ such
that $gz = z + \psi(g)$ for all $z \in Z$ and $g \in G$.
\bigskip

\sn{\sect 3 Topological spaces for point patterns}

\sn When $v$ is non-singular, $P_v$ forms an almost periodic
pattern of points in the sense that each spherical window, whose
position is shifted over the infinite pattern, reveals the same
configuration at a syndetic (relatively dense) set of positions
\sen . We may formulate this fact precisely in terms of  minimal
dynamics in a well-known process. Here we note the relevant
constructions and lemmas.

\sn{\bf Definition 3.1} Let $B(r)$ be the closed ball in $E$, centre $0$
and of radius $r$ with boundary $\partial B(r)$. Given a closed subset,
$A$, of $\R^N$, define $A[r] = (A \cap B(r)) \cup \partial B(r)$, a closed
subset of $B(r)$.
Consider the Hausdorff metric $d_r$ defined among closed subsets of $B(r)$
and define a metric (after \robone , \sol ) on closed subsets of the plane
by $$D(A,A') = \inf\{ 1/(r+1)\ |\  d_r(A[r], A'[r]) < 1/r \}.$$

\sn We are grateful to Johansen for pointing out that the topology induced by
$D$ on subsets of $E$ is precisely the topology induced when $E$ is embedded
canonically in its one-point compactification, the sphere of dimension $\dim
E$, and the Hausdorff metric is used to compare subsets of this sphere. Such
an observation proves quickly the following Proposition which appears first in
\rudolf\ (see also \robone\ and \rad ).

\sn{\bf Proposition 3.2} {\it If $u \in NS$, then the sets $\{P_v \ |\  v \in
NS \}$ and $\{P_v \ |\  v \in u + E \}$ are precompact with respect to $D$.
\qed}

\sn{\bf Definition 3.3} Define
$$MP=\overline{\{P_v \ |\  v \in NS \}}$$ and
$$MP_u=\overline{\{P_v \ |\  v \in u + E \}},$$
the completions of the above sets  with respect to $D$.
The symbol $M$ is used throughout this paper to indicate a construction
such as this: a ``Mapping Torus" or {\it continuous hull}.

\sn{\bf Remark 3.4} The term ``continuous hull''
(of $P_u$ with respect to $D$) simply refers to the fact that $MP_u$ is the
$D$-closure of the orbit of $P_u$ under the continuous group $E$.
A similar construction starts not
with a  point set or a tiling but with an operator on a Hilbert space
\belone, this is where the name came from. See \bhz\ for details and
a comparison.

Note that $\Delta =0$ if and only if $MP = MP_u$ for
all $u \in NS$, which happens if and only if $MP = MP_u$ for some $u \in NS$.

Also $P_v$ forms a Delone set (see \sol ), so we deduce that, for
$w \in E$ and  $\norm w$ small enough, $D(P_v,P_{v+w}) = \norm w /
( 1 + \norm w )$.

\sn{\bf Proposition 3.5} {\it Suppose that $w \in E$, then the map
$ P_v \mapsto P_{v+w}$, defined for $v \in NS$, may be extended to
a homeomorphism of $MP$, and the family of homeomorphisms defined
by taking all choices of  $w \in E$ defines a group action of
$\R^d \cong E$ on $NS$.

 Also for each $u \in NS$, $MP_u$ is invariant under this action of $E$ and
$E$ acts minimally on $MP_u$. \qed}

\sn The dynamical system $MP_u$ with the action by $E \cong \R^d$
is the dynamical system associated with the point pattern $P_u$,
analogous to that constructed by Rudolf \rudolf\ for tilings.  We
modify this to an action by $E$ on a non-compact cover of $MP_u$
as follows.

\sn{\bf Definition 3.6} For $v,v' \in \R^N$, write $\overline D(v,v') =
D(P_v,P_{v'}) + \norm {v-v'} $; this is clearly a metric. Let $\Pi$ be the
completion of $NS$ with respect to this metric.

\sn The following lemma starts the basic topological description of these
spaces.

\sn{\bf Lemma 3.7} {\it a/ The canonical injection $NS \lra \R^N$ extends
to a continuous surjection $\mu\colon  \Pi \lra \R^N$. Moreover, if $v \in
NS$, then $\mu^{-1}(v)$ is a single point.

b/ The map $v \mapsto P_v$, $v \in NS$, extends to a continuous
$E$-equivariant surjection, $\eta\colon  \Pi \lra MP$, which is an open map.

c/ The action by translation by elements of $E$ on $NS$ extends to a
continuous action of $\R^d \cong E$ on $\Pi$.

d/ Similarly the translation by elements of $\Z^N$ is $\overline
D$-isometric and extends to a continuous action of $\Z^N$ on $\Pi$. This
action commutes with the action of $E$ found in part c/.

e/ If $a \in MP$ and $b \in \R^N$, then $|\eta^{-1}(a) \cap \mu^{-1}(b)|
\leq 1$. }

\sn{\bf Proof} a/ The only non-elementary step of this part is the latter sentence.

We must show that if $v \in NS$ then for all $\epsilon > 0$ there is a
$\delta >0$ such that $\norm {w-v} < \delta$ and $ w \in NS$ implies that
$D(P_w,P_v) < \epsilon$. However, we know that if $B$ is a ball in $\R^N$
of radius much bigger than $1/(2 \epsilon)$, then $(\Z^N + v) \cap B$ is of
strictly positive distance, say at least $2 \delta$ with $\delta >0$ chosen
$< \epsilon /2$, from $\partial \Sigma$. Therefore, whenever $\pi(v-w)=0$
and $\norm {v-w} < \delta$, we have $P_v \cap B = P_w \cap B$ and hence
$D(P_v,P_w) < \delta$ . On the other hand, if $\pi(v-w) \neq 0$ but $\norm
{v-w} < \delta$ then we may replace $w$ by $w' = w + \pi(v-w)$, a
displacement by less than $\delta$. By the remark (3.4), we deduce that
$D(P_w, P_{w'}) < \delta$ and so we have $D(P_w,P_v) < 2 \delta < \epsilon$
in general, as required.

b/ The extension to $\eta$, and the equivariance and surjectivity,
are immediate.  The open map condition is quickly confirmed using
remark (3.4).

 c/ follows from  the uniform action of $E$ noted in Remark
3.4. d/ follows similarly where uniform continuity is immediate from the
isometry.

Note that e/ is a direct consequence of the definition of the metric
$\overline D$. \qed

\sn{\bf Definition 3.8} For $u \in NS$, let $\Pi_u$ be the completion of
$E+\Z^N+u$ with respect to the $\overline D$ metric.

\sn {\bf Lemma 3.9} {\it For $u \in NS$, $\Pi_u$ is a closed
$E+\Z^N$-invariant subspace of $\Pi$. If $x \in \Pi_u$, then $(E + \Z^N)x$,
the orbit of $x$ under the $E$ and $\Z^N$ actions, is dense in $\Pi_u$.
Consequently

a/ The injection, $E+\Z^N+u \lra \R^N$  extends to a continuous map, equal
to the restriction of $\mu$ to $\Pi_u$, $\mu_u\colon  \Pi_u \lra \R^N$,
whose image is $Q+u$.

b/ By extending the action by translation by elements of $E+\Z^N$ on
$E+\Z^N+u$, $E + \Z^N$ acts continuously and minimally on $\Pi_u$. This is
the restriction of the action of Lemma 3.7 c/ and d/.

c/ The map $v \mapsto P_v$, $v \in E+\Z^N+u$, extends to an open continuous
$E$-equivariant surjection, $\eta_u\colon  \Pi_u \lra MP_u$, which is the
restriction of $\eta$.

d/ If $x \in \Pi_u$ and $v \in E+\Z^N$ acts on $\Pi_u$ fixing $x$, then in
fact $v=0$.}

\sn{\bf Proof} The first sentence is immediate since, by definition,
$\Pi_u$ is the closure of an $E+\Z^N$ orbit in $\Pi$.

Suppose that $x \in \Pi_u$ and that $y \in E+\Z^N+u$ which we consider as a
subset of $\Pi_u$. Then there are $x_n \in E+\Z^N+u$ such that $x_n
\rightarrow x$ in the $\overline D$ metric. Write $\beta_n\colon  \Pi_u \lra
\Pi_u$ for the translation action by $-x_n$ and write $\alpha$ for the translation
action by $y$. Then we have $\mu(\beta_n(x)) \rightarrow 0$ and so
$\mu(\alpha\beta_n(x)) = y + \mu(\beta_n(x)) \rightarrow y$.

But, since $\mu$ is 1-1 at $y \in NS$ by Lemma 3.7 a/, we deduce that
$\overline D(\alpha\beta_n(x),y) \rightarrow 0$ and so $y$ is in the
closure of the $E+\Z^N$ orbit of $x$. However the orbit of $y$ is dense and
so we have the density of the $x$ orbit as well.

The lettered parts follow quickly from this. \qed

\sn By the results of parts b/ and c/ of Lemma 3.9, we may drop the suffix
$u$ from the maps $\mu_u$ and $\eta_u$ without confusion, and this is what
we do unless it is important to note the domain explicitly.

The aim of the next few sections is to fill in the fourth corner of the
commuting square $$\matrix{ \Pi_u & \bbra {\eta} & MP_u \cr && \cr
\biggl\downarrow {}^{\mu} &&\biggr\downarrow {}^{?}  \cr && \cr (Q+u) &
\bbra {?} & ?
\cr}$$ in a way which illuminates the underlying structure.

\bigskip

\sn{\sect 4 Tilings and Point Patterns}

\sn
We now connect the original
construction of projection tilings due to Katz and Duneau \kd\ with the point
patterns that we have been considering until now. We refer to \okd\ and
\sen\ for precise descriptions of the construction; we extract the points
essential for our argument below.

We note two developments of the $D$ metric (3.1) which will be used ahead.
The first development is also E.A.Robinson's original application of $D$ \robone .

\sn{\bf Definition 4.1}  We suppose that we have a finite set of
pointed  compact  subsets of $E$ which we call the {\it units},
and we suppose that we have a uniformly locally finite subset of
$E$, a {\it point pattern}. A {\it pattern}, $\T$, in $E$ (with
these units and underlying point pattern) is an arrangement of
translated copies of the units in $E$, the distinguished point of
each copy placed over a point of the point pattern, no point of
the point pattern uncovered and no point of the point pattern
covered twice. Sometimes, an underlying point pattern is not
mentioned explicitly.

For example we could take a
tiling of $E$ and let the pattern consist of the boundaries of the tiles with
superimposed decorations, {\it i.e.\ }  small compact sets, in their
interior giving further asymmetries or other distinguishing features. Or we
could take a point pattern, perhaps replacing each point with one of a finite
number of decorations.
See \gs\ for a thorough discussion of this process in general.

By taking the union of all the units of the pattern, we obtain a
locally compact  subset $P(\T)$ of $E$ which can be shifted by
elements of $E$, $P(\T) \mapsto P(\T)+v$, and these various
subsets of $E$ can be compared using $D$ literally as defined
(3.1) (the addition of further decorations can also solve the
problem of confusing overlap of adjacent units of the pattern, a
complication which we ignore therefore without loss of
generality). Under natural conditions (see \rudolf\ \sol ), which
are always satisfied in our examples, the space $\{P(\T)+v \ |\  v
\in E\}$ is precompact with respect to the $D$ metric and its
closure, the {\it continuous hull} of $\T$ written $M\T$ here,
supports a natural continuous $E$ action. The {\it pattern
dynamical system} of $\T$ is this dynamical system $(M\T,E)$.

\sn{\bf Definition 4.2} The second development adapts $D$ to
compare subsets of $\Sigma$. Recall the notation $B(r)$ for the
closed Euclidean $r$-ball  in $E$ (3.1). Let $C(r) =
\pi^{-1}(B(r)) \cap \Sigma$ and let $dC(r) = \pi^{-1}(\partial
B(r)) \cap \Sigma$.

Given a subset, $A$, of $\Sigma$ define $A[r] = (A \cap C(r) )\cup dC(r)$.
Let $d'_r$ be the Hausdorff metric defined among closed subsets of $C(r)$
and define a metric on subsets of $\Sigma$ by $$D'(A,A') = \inf \{
1/(r+1) \ |\  d'_r(A[r], A'[r]) < 1/r \}$$ Let $\overline D' (v,w) =
D'(\widetilde P_v,\widetilde P_w) + \norm{v-w}$, where we recall that
$\widetilde P_v = \Sigma \cap (\Z^N + v)$.

Let $M\widetilde P_u$ be the $D'$-closure of the space $\{ \widetilde P_v \ |\
v \in E+u \}$, and let $\widetilde \Pi_u$ be the $\overline D'$ completion
of $NS \cap (Q+u)$. Let $M\widetilde P$ be the $D'$ closure of the space
$\{ \widetilde P_v \ |\  v \in NS \}$.

The analogues of Proposition 3.5 and Lemma 3.9 with respect to $\widetilde
P$, $M\widetilde P$, $\widetilde \Pi_u$, $M\widetilde P_u$ and $Q+u$,
continue to hold and so we define maps $\widetilde \mu \colon \widetilde \Pi_u
\lra Q+u$ and $\widetilde \eta \colon \widetilde \Pi_u \lra M\widetilde P_u$.

\sn We use the projection $\pi$ to compare the strip pattern with
the projection pattern. It will turn out that $\pi$ is a
homeomorphism between $\Pi_u$  and $\widetilde \Pi_u$, but that
the definition of $\widetilde \Pi_u$ will be more convenient than
that of $\Pi_u$. Using $\pi$ we may work with either space.

\sn{\bf Theorem 4.3} {\it There are $E$-equivariant maps $\pi_*$ induced by
the projection $\pi$ which complete the commuting square
$$\matrix{ \widetilde \Pi_u & \bbra {\pi_*} & \Pi_u \cr && \cr \biggl\downarrow
{}^{\widetilde \eta} && \biggr\downarrow {}^{\eta}  \cr && \cr M\widetilde P_u &
\bbra {\pi_*} & MP_u \cr}$$ Furthermore we have the following commuting square
$$\matrix{ \widetilde \Pi_u & \bbra {\pi_*} & \Pi_u \cr && \cr \biggl\downarrow
{}^{\widetilde \mu} && \biggr\downarrow {}^{\mu}  \cr && \cr Q+u&=\! =\! =\! = &
Q+u \cr}$$  in which all the labelled maps are 1-1 on $NS$. \qed}

\sn Consider the example of the canonical projection tiling,
$\T_u$ (2.3). If we  know $\widetilde P_u$ then we have all the
information needed to reconstruct $\T_u$ by its definition.
Conversely, the usual assumption that the projected faces are
non-degenerate (see \le\ (3.1)) allows us to distinguish the
orientation of the lattice face (in $\Z^N$) from which a given
tile came. Piecing together all the faces defined this way obtains
$\widetilde P_u$. So the canonical projection tiling is conjugate
(in the sense defined ahead in 4.5) to $\widetilde P_u$.

On the other hand, the well-known Voronoi or Dirichlet tiling \gs\
obtained from a  point pattern in $E$ is a tiling conjugate to the
original point pattern provided we decorate each tile with the
point which generates it.

With these two examples of tiling in mind, we consider the pattern
$\widetilde P_u$  to represent the most elaborate tiling or
pattern that can be produced by the projection method, without
imposing further decorations not directly connected with the
geometry of the construction, and at the other extreme, the point
pattern, $P_u$, represents the least decorated tiling or pattern
which can be produced by the projection method.

\sn{\bf Definition 4.4}  For a given $E$ and $K$ as in (2.1),  we
include in the {\it class of projection method patterns} all those
patterns, $\T$, of $\R^d$ such that there is a $u \in NS$ and two
$E$-equivariant  surjections $$M\widetilde P_u \lra M\T \lra
MP_u$$ whose composition is $\pi_*$.

We call $(E,K,u)$ the {\it data} of the projection method and by
presenting these  data we require tacitly that $K$ has the
properties of Definition 2.1, that $u \in NS$ and that $(Q+u) \cap
Int K \neq \emptyset$ (2.6).

\sn Thus the tilings of \okd\ and the Voronoi tilings discussed
above are examples from this class when $K =
\pi^{\perp}([0,1]^N)$. In order to compare these  two
constructions, or to consider projection method patterns in the
general sense of (4.4), we aim to describe $\pi_*\colon
M\widetilde P_u \lra MP_u$.

First we adopt the following definitions which possibly duplicate
notions already  existing in the literature.

\sn{\bf Definitions 4.5} Adapting a definition of Le \le , we say
that two patterns, $\T,\T'$, in $E$ are  {\it topologically
conjugate} if there is an $E$-equivariant homeomorphism, $M\T
\leftrightarrow M\T'$.

The two patterns, $\T,\T'$, are {\it pointed conjugate} if there
is an  $E$-equivariant homeomorphism, $M\T \leftrightarrow M\T'$,
which maps $\T$ to $\T'$.

A pattern $\T'$ is a {\it finite decoration} of a pattern $\T$ if
there are real  numbers $r$ and $s$ so that the following happens:
i/ $\T$ may be constructed from $\T'$ by a transformation which
alters the unit of $\T$ at a point $v \in \R^d$ according only to
how $\T$ appears in the ball $v + B(r)$: and ii/, conversely, if,
for any choice of $w \in \R^d$ we know what $\T'$ looks like in
the ball $w+B(s)$, then we can construct the remainder of $\T'$
from $\T$ by a transformation (depending perhaps on the appearance
of $\T'$ in the ball $w+B(r)$) which alters the unit of $\T'$ at a
point $v \in \R^d$ according only to how $\T'$ appears in the ball
$v + B(r)$.

Finally, we say that a pattern is a {\it Meyer pattern} if the
underlying  point pattern is a Meyer set, that is a set $M$ for
which we can find $R$ and $r$ so that $M-M$ intersects every
$R$-ball in at least one point and intersects every $r$ ball in at
most one point.

\sn{\bf Remark 4.6} Note that all the patterns we consider in this
paper are (pointed conjugate to) Meyer patterns. This is the
starting point for Schlottmann's analysis of the projection method
\schtwo .

To tie these definitions in to the existing literature, we note
that topological conjugacy is strictly weaker than local
isomorphism (as in \le\ for example) and strictly stronger than
equal quasicrystal type \robone . Pointed conjugacy is strictly
stronger than mutual local derivability \bsj\ and topological
equivalence \keltwo , but has no necessary relation with local
isomorphism and quasicrystal type. Finite decoration is strictly
weaker than local derivability \bsj .

However, we have the following, an immediate application of the definitions
to the fact that an $n$-to-$1$ factor map (see 2.13) is an open map \furst .

\sn{\bf Lemma 4.7} {\it Suppose we have two Meyer patterns,
$\T,\T'$, in  $E$ and a continuous $E$-equivariant surjection
$M\T' \lra M\T$ which is n-to-1, sending $\T'$ to $\T$. Then $\T'$
is a finite decoration of $\T$. \qed}

\newpage

\sn{\sect 5  Comparing $\Pi_u$ and $\widetilde \Pi_u$}

\sn We start by examining $\pi_*\colon  \widetilde \Pi_u \lra
\Pi_u$ from (4.3) and seek conditions under which it is a
homeomorphism. As the section proceeds we shall find  that the
conditions can be whittled away to the minimum possible. Recall
the space $V$, one of the orthocomponents of the decomposition of
$Q$ in Corollary 2.11. It is the connected component of
$\overline{\pi^\perp(\Z^N)}$ containing $0$.

\sn{\bf Lemma 5.1} {\it Suppose  that $u \in NS$ and that, for all $v
\in Q+u$ such that $v \in \partial( (V+v) \cap Int K)$ (the boundary taken
in $V+v$), we have $(\Delta + v) \cap K = \{v\}$; then $\pi_*\colon
\widetilde
\Pi_u \lra \Pi_u$ is an $E$-equivariant homeomorphism.}

\sn{\bf Proof} We ask under what circumstances could we find $x \in \Pi_u$
with two preimages under $\pi_*$ in $\widetilde \Pi_u$? We would need two
sequences $v_n,w_n \in (Q+u)\cap NS$ both converging to $x$ in the
$\overline D$ metric such that $\widetilde P_{v_n}$ and $\widetilde
P_{w_n}$ have different $\overline D'$ limits, say $A$ and $B$
respectively. From this we see that $A \Delta B \subset \partial \Sigma$
(symmetric difference) and yet $\pi(A) = \pi(B)$.

Let $p \in \pi(A \Delta B)$ and consider the set $(A \Delta B) \cap
\pi^{-1}(p)$. As noted above, this set is a subset of the boundary of
$\Sigma \cap \pi^{-1}(p) \equiv K$ and each pair of elements is separated
by some element of $\Delta$.

Suppose that $a \in (A \setminus B) \cap \pi^{-1}(p)$. By construction,
there are $a_n \in (Q+u) \cap NS \cap \widetilde P_{v_n}$ converging to $a$
implying that $a \in \partial( (Q+u) \cap Int K)$. But by hypothesis, we
deduce $B \cap \pi^{-1}(p) = \emptyset$ - a contradiction to the fact that
$p \in \pi(A) = \pi(B)$.

A symmetric argument produces a contradiction from $b \in (B \setminus A)
\cap \pi^{-1}(p)$. \qed

\sn Note that if $\Delta = 0$ or, more generally, if $K \cap (K + \delta) =
\emptyset$ whenever $\delta \in \Delta$, $\delta \neq 0$, then the
hypothesis of the Lemma is satisfied trivially.

\sn{\bf Corollary 5.2} {\it If $\Delta = 0$, then $\pi_*\colon
\widetilde
\Pi_u \lra \Pi_u$ is an $E$-equivariant homeomorphism. \qed}

\sn In special cases the
hypothesis is satisfied less trivially. We give a slightly more special
condition here.

\sn{\bf Lemma 5.3} {\it Suppose that $J$ is the closure of a
fundamental domain for $\Delta$ in $E^{\perp}$,  and that $J =
\overline{Int J}$ in $E^{\perp}$. If $K$ is contained in some
translate of $J$, then $\pi_*\colon \widetilde \Pi_u \lra \Pi_u$
is a homeomorphism. In particular, if $K = \pi^{\perp}([0,1]^N)$,
then $\pi_*\colon \widetilde \Pi_u \lra \Pi_u$ is a
homeomorphism.}

\sn{\bf Proof} For the first part, suppose that $a,b \in K$ and $0 \neq a-b
= \delta \in \Delta$, then, by construction, $a$ and $b$ sit one in each of
two hyperplanes orthogonal to $\delta$ between which $K$ lies. Note that
then these hyperplanes are therefore both parallel to $V$ and each
intersects $K$ only in a subset of $\partial K$. Therefore, $a,b \in
\partial K$ and further, since $V+a$ and $V+b$ are contained one in each of
the hyperplanes, we have $a \nin \partial ((V+a) \cap Int K)$
(boundary in $V+a$) and $b \nin \partial ((V+b) \cap Int K)$ ({\it
similis}). Therefore the conditions of Lemma 5.1 are fulfilled
vacuously.

In the second part, suppose that $K = \pi^{\perp}([0,1]^N)$ and that
$\delta = (\delta_1,\delta_2,...,\delta_N) \in \Delta$, $\delta \neq 0$
(the case $\Delta = 0$ is easy). Consider the set $I = \{\langle \delta, t
\rangle \ |\  t \in K\}$, where $\langle .,. \rangle$ is the inner product on
$\R^N$. This is a closed interval. Also, since $\delta$ is fixed by the
orthonormal projection $\pi^{\perp}$, $I = \{\langle \delta, s \rangle \ |\  s
\in [0,1]^N\}$, from which we deduce that the length of $I$ is $\sum_j
|\delta_j|$. But since $|\delta_j| < 1$ implies that $\delta_j = 0$, we
have $\langle \delta,\delta \rangle = \sum |\delta_j|^2 \geq \sum
|\delta_j|$ and so $K$ can be fitted between two hyperplanes orthogonal to
$\delta$ and separated by $\delta$.

Therefore $K$ is contained in a translate of $\cap_{\delta
\in \Delta, \delta \neq 0} \{v \in E^{\perp} \ |\ | \langle v,\delta \rangle | \leq
(1/2)\langle \delta, \delta \rangle \ \}$, which in turn is contained in
the closure of a fundamental domain for $\Delta$. So we have confirmed the
conditions of the first part. \qed

\sn{\bf Remark 5.4} Using Lemma 2.5, the condition of Proposition
5.3 is  equivalent to the following condition: $((IntK) - (IntK))
\cap \Delta = \{0\}$, where we write $A-A = \{a-b \ |\ a,b \in
A\}$ for the arithmetic (self-)difference of $A$, a subset of an
abelian group. Compare with 8.2 ahead.

All of these results say that if $K$ is small enough relative to
$\Delta$ then $\pi_*$ is a homeomorphism. In fact, we can
dispense with all such conditions, and the following construction
gives a procedure to reduce the size of a general acceptance
domain appropriately.

\sn{\bf Theorem 5.5} {\it Suppose that $E$ and $K$ are as in definition 2.1,
and that $u \in NS$. Then $\pi_* : \widetilde \Pi_u \lra \Pi_u$ is an
$E$-equivariant homeomorphism. }

\sn{\bf Proof} Suppose that $E$, $K$ and $u$ are chosen as
required and  that $\Delta \neq 0$. The case $\Delta = 0$ is
covered by Corollary 5.2.

Suppose that $J$ is the closure of a fundamental domain for
$\Delta$ in $E^{\perp}$, such that $J = \overline{Int J}$ in
$E^{\perp}$ and suppose,  as we always can by shifting $J$ if
necessary, that $\partial J \cap (Q+u) = \emptyset$.

Let $K' = (K+\Delta) \cap J$. Then $K'$ is a subset of $E^{\perp}$
which obeys the conditions required in the original definition of
(2.1). Also the placement of $J$ ensures that the points in $Q+u$,
in particular $u$ itself, which are non-singular with respect to
$K$  are also non-singular with respect to $K'$.

Moreover, if we define $\Sigma' = K' + E$, then, by construction,
$\pi(\Sigma' \cap (v + \Z^N)) = \pi(\Sigma \cap (v+\Z^N)) = P_u$
for all $v \in \R^N$. Therefore, working with $\Sigma'$ instead of
$\Sigma$, we can retrieve the projection point pattern, $P_u$. So,
by Lemma 5.3 and the fact that $K' \subset J$, we see that $\pi_*$
is a homeomorphism between the  spaces $\Pi_u$ and $\widetilde
\Pi_u(\Sigma')$, the construction of (4.2) with respect to
$\Sigma'$.

However, for any $v \in E+\Z^N+u$, we have the equalities:
$\Sigma'\cap (v + \Z^N) = \Sigma' \cap ((\Sigma \cap (v + \Z^N))+
\Delta)$ and $\Sigma\cap (v + \Z^N) = \Sigma \cap ((\Sigma' \cap
(v + \Z^N))+ \Delta)$. Therefore the set $\Sigma' \cap (v + \Z^N)$
can be constructed from $\Sigma \cap (v + \Z^N) = \widetilde P_v$
and vice versa. Moreover, this correspondence defines a $\overline
D$ metric isometry, between $\widetilde \Pi_u(\Sigma')$ and
$\widetilde \Pi_u$, which intertwines $\pi_*$. Completing the
correspondence gives an isometry between $\widetilde
\Pi_u(\Sigma')$ and $\widetilde \Pi_u$ which intertwines $\pi_*$,
as required. \qed

\sn{\bf Remark 5.6} We note a second process of reduction without loss of
generality. Until now we have assumed nothing about the rational position
of $E$, but it is convenient to assume and is often required in the
literature that $E \cap \Z^N = 0$: the {\it irrational} case.

If we do not assume this then we can always reduce to the
irrational  case by quotienting out the rational directions. A
simple argument allows us to find in the most general case of
projection method pattern, $P$ say, an underlying irrational
projection method pattern, $P'$, with, $MP = MP' \times \torus^k$
for some value of $k$; and this torus factor splits naturally with
respect to the various constructions and group actions we find
later. We leave the details to the reader.

\bigskip

\sn{\sect 6 Calculating $M\widetilde P_u$ and $MP_u$}

\sn We now describe $M\widetilde P_u$ and $MP_u$ as quotients of
$\widetilde \Pi_u$ and $\Pi_u$ respectively. Although the last
section established an  equivalence between $\widetilde \Pi_u$ and
$\Pi_u$, we find it useful to distinguish the two constructions.

First we examine $M\widetilde P_u$ and prove a generalisation of (3.8) of
\le .

\sn{\bf Proposition 6.1} {\it Suppose that $u \in NS$, then there is an
isometric action of $\Z^N$ on $\widetilde \Pi_u$, which factors by
$\widetilde \mu$ to the translation action by $\Z^N$ on $Q+u$, and
$M\widetilde P_u = \widetilde \Pi_u/\Z^N$. Thus we obtain a commutative
square of $E$ equivariant maps $$\matrix{ \widetilde \Pi_u & \bbra
{\widetilde \eta} & M\widetilde P_u \cr && \cr \biggl\downarrow
{}^{\widetilde\mu} && \biggr\downarrow {}^{\widetilde \mu}  \cr && \cr Q+u
& \llra & (Q+u)/\Z^N.
\cr}$$ The left vertical map is 1-1 precisely at the points in $NS \cap
(Q+u)$. The right vertical map is 1-1 precisely on the same set, modulo
the action of $\Z^N$.}

\sn{\bf Proof} The action of $\Z^N$ on $\widetilde \Pi_u$, as an extension
of the action on $Q+u$ by translation, is easy to define since the maps are
$\overline D'$-isometries.

If $v,w \in NS$ then it is clear that $\widetilde P_v = \widetilde P_w$ if
and only if $v-w \in \Z^N$. Moreover, there is $\delta > 0$ so that
$\norm{v-w} < \delta$ implies that $D'(\widetilde P_v, \widetilde P_w) \geq
\norm{v-w} /2$.

{From} this we see that, if $\widetilde P_v = \widetilde P_w$ and $\widetilde
P_{v'} = \widetilde P_{w'}$ and $\norm{v-v'} < \delta/2$ and $\norm{w-w'} <
\delta/2$,  then $v-w = v'-w'$. The uniformity of $\delta$ irrespective of
the choice of $v,w,v'$ and $w'$ shows that the statement $\widetilde
\eta(v)=\widetilde \eta(w)$ implies $\widetilde \mu(v)-\widetilde \mu(w)
\in \Z^N$, which is true for $v,w \in NS \cap (Q+u)$, is in fact true for
all pairs in $\widetilde \Pi_u$, the $\overline D'$ closure.

To show the 1-1 properties for the map on the left, suppose that $v \in
Q+u$ and that $p \in \partial \Sigma \cap (\Z^N+v)$, i.e.\ $v \nin NS$. Then
since $K$ is the closure of its interior and since $NS$ is dense in $\R^N$
(Lemma 2.2), there are two sequences $v_n, v'_n \in NS $ both converging to
$v$ in Euclidean topology and such that $p+(v_n-v) \in \Sigma$ and $p +
(v'_n-v) \nin \Sigma$. This implies that any $D'$ limit point of
$\widetilde P_{v_n}$ contains $p$ and any $D'$ limit point of $\widetilde
P_{v'_n}$ does not contain $p$. But both such limit points (which exist by
compactness of $M \widetilde P$) are in ${\widetilde \mu}^{-1}(v)$ which is
a set of at least two elements therefore.

The 1-1 property for the map on the right follows directly from this and
the commuting diagram. \qed

\sn The space $(Q+u)/\Z^N$ and its $E$ action, which is being compared with
$M\widetilde P_u$, also has a simple description.

\sn{\bf Lemma 6.2} {\it With the data above, $(Q+u)/\Z^N$ is a coset of the
closure of $E \mod \Z^N$ in $\R^N/\Z^N \equiv \torus^N$. Therefore $(Q+u)/\Z^N$ with
its $E$ action is isometrically conjugate to a minimal action of
$\R^d$ by rotation on a torus of dimension $N-\dim\Delta$.}

\sn{\bf Proof} The space $Q \mod \Z^N$ is equal to the closure of
$E \mod \Z^N$ and its translate by $u \mod \Z^N$ is an isometry
which is $E$ equivariant. The action of $E$ on its closure is
isometric and transitive, hence minimal, and is by translations.
$E$ is a connected subgroup of $\torus^N$ and so also is the
closure of $E$, which is therefore equal to a torus of possibly
smaller dimension. The codimension  of this space agrees with the
codimension of $V+E$ (the continuous component of $Q$) in $\R^N$
which, by Lemma 2.9 and Corollary 2.11, equals $\dim \Delta$ as
required. \qed

\sn Now we turn to a description of $MP_u$ which is similar in form to that
of $M\widetilde P_u$, but  as to be shown in examples 8.7 and 8.8, need not be
equal to $M \widetilde P_u$.

\sn{\bf Lemma 6.3} {\it Suppose that $u \in NS$. If $v,w \in NS\cap (Q+u)$
and $P_v = P_w$  then there are $v^* \in v + \Z^N$ and $w^* \in w + \Z^N$
such that $v^*,w^* \in \Sigma$ and $\pi(v^*)=\pi(w^*)$, and with this
choice $\widetilde P_v+\Delta-\pi^{\perp}(v^*) = \widetilde P_w + \Delta -
\pi^{\perp}(w^*)$.}

\sn{\bf Proof} Fix $p_o \in P_v=P_w$ and let $v^* \in \widetilde P_{v}$ be
chosen so that $\pi(v^*) = p_o$ and similarly, let $w^* \in \widetilde
P_{w}$ be chosen so that $\pi(w^*)=p_o$. Clearly $v^*$ and $w^*$ obey the
conditions required. Also $\widetilde P_w - \pi^{\perp}(w^*)$ and
$\widetilde P_v - \pi^{\perp}(v^*)$ are both contained in $p_o + \Z^N$ and
project under $\pi$ to the same set $P_v$. Thus the difference of two
points, one in $\widetilde P_w - \pi^{\perp}(w^*)$ and the other in
$\widetilde P_v - \pi^{\perp}(v^*)$, and each with the same image under
$\pi$, is an element of $\Delta$ as required.
 \qed

\sn{\bf Proposition 6.4} {\it Suppose that $x,y \in \Pi_u$ and
that $\eta(x) = \eta(y)$, then there is a $v \in Q$ and a
$\overline D$ isometry $\phi\colon  \Pi_u \lra \Pi_u$ so that
$\phi(x) = \phi(y)$ and the following  diagram commutes
$$\matrix{\Pi_u &  \bbra {\phi} & \Pi_u  \cr && \cr
\biggl\downarrow {}^{\mu} && \biggr\downarrow {}^{\mu} \cr && \cr
Q+u & \bbra {w \mapsto w+v} & Q+u \cr}$$ and $\eta_u \phi =
\eta_u$ (here the restriction to $\Pi_u$ is important to note,
c.f. 3.9). In this case we deduce $v+u \in NS$.

 Conversely, if we have such an isometry in such a diagram and if
$P_{u+v}=P_u$, then $v+u \in NS$ and $\eta_u \phi = \eta_u$ automatically.}

\sn{\bf Proof} Suppose $w \in E+\Z^N$ and that $\alpha_w \colon
\Pi_u \lra \Pi_u$ is the map completed from the map $z \mapsto
z+w$ defined first for  $z \in NS \cap (Q+u)$ (see Proposition
3.5). Then, since $\eta(x) = \eta(y)$ and $\eta$ is
$(E+\Z^N)$-equivariant, we have $\eta(\alpha_w(x)) =
\eta(\alpha_w(y))$ for all $w \in E+\Z^N$. So, by definition, the
map $\alpha_w(x) \mapsto \alpha_w(y)$ defined point-by-point for
$w \in E + \Z^N$ is a $\overline D$ isometry from the
$(E+\Z^N)$-orbit of $x$ onto the $(E+\Z^N)$-orbit of $y$ (By Lemma
3.9 d/ the mapping is well-defined). These two orbits being dense
(Lemma 3.9) in $\Pi_u$,  this map extends as an isometry onto,
$\phi\colon  \Pi_u \lra \Pi_u$.

Since $\mu$ is $(E+\Z^N)$-equivariant, we deduce the intertwining with
translation by $v = \mu(y)-\mu(x)$. Also since $\eta \phi = \eta$ on the
$E+\Z^N$ orbit of $x$, the $E$-equivariance of $\eta$ extends this equality
over all of $\Pi_u$.

Conversely suppose we have an isometry which intertwines the translation by
$v$ on $Q+u$. Then for general topological reasons the cardinality of the
$\mu$ preimage of a point in $Q+u$ is preserved by translation by $v$ and
we deduce that $NS \cap (Q+u)$ is invariant under the translation by $v$.
In particular $u+v \in NS$. The equation follows since it applies, by
hypothesis and Lemma 6.3, at $u$ and therefore, by equivariance, at all
points in $E+\Z^N+u$, a dense subset. \qed

\sn{\bf Definition 6.5} For $u \in NS$, let $R_u = \{v \in Q \ |\  v+u \in NS,
P_{u+v} = P_u \}$.

\sn{\bf Corollary 6.6} {\it Suppose that $u \in NS$ and $w \in NS \cap
(Q+u)$, then $R_w = R_u$. Therefore, if $v \in R_u$, then $v+w \in NS \cap
(Q+u)$ for all $w \in NS \cap (Q+u)$.}

\sn{\bf Proof} By Lemma 3.9, we know that $\Pi_w \bra {\eta_w} Q+w$ equals
$\Pi_u \bra {\eta_u} Q+u$ and so any isometry of $\Pi_u$ which factors by
$\eta$ through to a translation by $v$ also does the same for $\Pi_w$.
Proposition 6.4 completes the equivalence.

The second sentence follows directly from the definition of $R_w$. \qed

\sn{\bf Remarks 6.7} It would be natural to hope that the
condition $u+v \in NS$ could be removed from the definition of
$R_u$. We have been unable to do this in general. But since $NS$
is a dense $G_{\delta}$ set (2.2) and, anticipating Theorem 7.1,
$R_u$ is countable, we see that for a dense $G_{\delta}$ set of
$u\in NS$ ({\it generically}) we can indeed equate $R_u = \{v \in
Q : P_{u+v} = P_u\}$.

This is bourne out in Corollary 6.6 where we see that $R_u$ is
defined independently of the choice of $u$ generically, and $R_u$
can be thought of as an invariant of $\Pi_u$. This result also
shows that $R_u$ is a subset of the translations of $\R^N$ which
leave $NS \cap (Q+u)$ invariant.

Note that, since $\mu$ is 1-1 only on $NS$, $R_u$ could as well
have been defined as $\{v \in Q \ |\ \eta\mu^{-1}(u+v)= \{P_u\} \
\}$.

It is clear that $\Z^N \subset R_u$.

\sn{\bf Theorem 6.8} {\it If $u \in NS$, then $R_u$ is a closed
subgroup of $Q$. Also $R_u$ acts by $\phi$ isometrically on $\Pi_u$ and
defines a homeomorphism $\Pi_u/R_u \equiv MP_u$. Moreover the $R_u$ action
commutes with the $E$-action, so the homeomorphism is $E$-equivariant.}

\sn{\bf Proof} The main point to observe is that $R_u$ consists precisely
of those elements $v$ such that there is an isometry $\phi_v$ as in
Proposition 6.4 with $\eta_u \phi_v = \eta_u$. Since the inverse of such an
isometry is another such, and the composition of two such isometries
produces a third, we deduce the group property for $R_u$ immediately. The
isometric action is given to us and Proposition 6.4 shows directly that
$\Pi_u/R_u \equiv MP_u$.

Closure of $R_u$ is more involved. Suppose that $v_n \in R_u$ and that $v_n
\rightarrow v$ in the Euclidean topology. Then $\phi_{v_n}$ is uniformly
Cauchy and so converges uniformly to a bijective isometry, $\psi$, of
$\Pi_u$ which intertwines the translation by $v$ on $Q+u$. Therefore, if
$\mu^{-1}(u+v)$ has at least two elements, then so also does
$\psi^{-1}\mu^{-1}(u+v)$, but this set is contained in $\mu^{-1}(u)$, a
contradiction since $\mu^{-1}(u)$ is a singleton.  Therefore $u+v \in NS$
and $\mu^{-1}(u+v) = \psi \mu^{-1}(u) = \lim \phi_{v_n}\mu^{-1}(u) = \lim
\mu^{-1}(u+v_n) = \mu^{-1}(u)$.  Thus $v \in R_u$ and so $R_u$ is closed.

The commutativity with the $E$ action on $\Pi_u$ is immediate from the
corresponding commutativity on $Q+u$. \qed

\sn{\sect 7 Comparing  $MP_u$ with $M\widetilde P_u$}

\sn
In Section~4 we defined projection method patterns as those whose
dynamical system sits intermediate to $M\widetilde P_u$ and $MP_u$. We
discover in this section how closely these two spaces lie and circumstances
under which they are equal.

To compare $MP_u$ with $M\widetilde P_u$ we start with the fact
that  $\widetilde \Pi_u = \Pi_u$ (5.5). By Proposition 6.1 and
Theorem 6.8, therefore, the problem becomes the comparison of
$R_u$ with $\Z^N$. Perhaps surprisingly, under general conditions
we find that $R_u$ is not much larger than $\Z^N$ and under normal
conditions the two groups are equal.

\sn{\bf Theorem 7.1} {\it For all $u \in NS$, $\Z^N \subset R_u$ as a
finite index subgroup. In fact, with the notation of (2.10), $R_u \subset
\Z^N+\widetilde\Delta$.}

\sn{\bf Proof} Suppose that $v \in R_u$. Then in particular, by (6.7),
$P_{v+u} = P_u$. Therefore there is an $a \in \Z^N$ such that $\pi(v+u+a) =
\pi(u)$ and so by translating if necessary, we may assume without loss of
generality that $v \in E^{\perp}$; and this defines $v$ uniquely $\mod
\Delta$.

With this assumption we deduce from Lemma 6.3 that $\widetilde
P_{u+v}+\Delta= \widetilde P_u + \Delta +v$ In particular,
$\pi^{\perp}(\widetilde P_{u+v})+\Delta= \pi^{\perp}(\widetilde P_u) +
\Delta +v$.

Now each of $\pi^{\perp}(\widetilde P_{u+v})$ and
$\pi^{\perp}(\widetilde P_u)$ is contained in $K$ a compact set.
Suppose that $\alpha \in \Delta^{\perp}$ and that $\langle
v,\alpha \rangle \neq 0$, then there is $t \in \Z$ such that
$|\langle tv,\alpha \rangle | > 2 \norm {\alpha} \ diam K$.
However, since $tv \in R_u$ by Theorem 6.3, we have
$\pi^{\perp}(\widetilde P_{tv+u})+\Delta= \pi^{\perp}(\widetilde
P_u) + \Delta +tv$. Applying the function $\langle .,\alpha
\rangle$ to both sets produces a contradiction by construction.
Thus we have $v \in U$, the space  generated by $\Delta$ (see
2.9). But if $v \in R_u$ then $v \in Q$ by definition, and so we
have $v \in \widetilde\Delta$, a group which, by Corollary 2.11,
contains $\Delta$ with finite index. \qed

\sn We note, for use in section 10, that therefore $R_u$ is free
abelian on $N$ generators.

\sn{\bf Corollary 7.2} {\it If $\Delta = 0$, then $R_u = \Z^N$ and
$\pi_* \colon M\widetilde P_u \leftrightarrow MP_u$. \qed}

\sn The following combines Propositions 6.1 and 6.8 and fits the present
circumstances to the conditions of Lemma 2.14. Recall the definitions of
$n$-to-1 extension (2.13) and finite isometric extension (2.15).

\sn{\bf Proposition 7.3} {\it Suppose that $u \in NS$. The map
$\pi_*\colon  M\widetilde P_u \lra MP_u$ is $p$-to-1 where $p$ is
the index of $\Z^N$ in $R_u$. The group $R_u$ acts isometrically
on $M\widetilde P_u$, commuting with the $E$ action, preserving
$\pi_*$ fibres and acting transitively on each fibre. This action
is mapped almost 1-1 by $\widetilde \mu$ to an action by $R_u$ on
$(Q+u)/\Z^N$ by rotation, and so we complete a commuting square
$$\matrix{M\widetilde P_u &  \bbra {\pi_*} & MP_u  \cr && \cr
\biggl\downarrow {}^{\widetilde \mu} && \biggr\downarrow {}^{\mu}
\cr && \cr (Q+u)/\Z^N & \bbra {mod\ R_u/\Z^N} & (Q+u)/R_u.
 \cr}$$
\qed}

\sn From this and the construction of Lemma 2.14 applied to the case $G=E$
and $H = R_u $, we deduce the main theorem of the section.

\sn{\bf Theorem 7.4} {\it Suppose that, $E$, $K$ and $u \in NS$
are  data for $\T$, a projection method pattern. Then there is a
group $S_{\T}$, intermediate to $\Z^N < R_u$, which fits into a
commutative diagram of $E$ equivariant maps $$\matrix{M\widetilde
P_u & \llra & M\T & \llra & MP_u  \cr &&&& \cr \biggl\downarrow
{}^{\widetilde \mu} && \biggr\downarrow && \biggr\downarrow
{}^{\mu} \cr &&&& \cr (Q+u)/\Z^N & \bbra {mod\ S_{\T}/\Z^N} &
(Q+u)/S_{\T} & \bbra {mod\ R_u/S_{\T}} & (Q+u)/R_u \cr}$$ where
the top row maps are finite isometric extensions and the bottom
row maps are group quotients.

Conversely, every choice of group $S'$ intermediate to $\Z^N <
R_u$  admits a projection method pattern, $\T$, fitting into the
diagram above with $S'=S_{\T}$. \qed}

\sn With the considerations of section 4 (in particular using 4.7)
we can count the projection method patterns up to topological
conjugacy or pointed conjugacy in the following corollary.

\sn{\bf Corollary 7.5} {\it With fixed projection data and the
conditions of Theorem 7.4, the set of topological conjugacy
classes of projection method patterns is in bijection with the
lattice of subgroups of $R_u/\Z^N$. Moreover, each projection
method pattern, $\T$, is pointed conjugate to a finite decoration
of $P_u$, and $\widetilde P_u$ is pointed conjugate to a finite
decoration of $\T$. \qed}

\sn We deduce the following result also found by Schlottmann
\schtwo ,  which is in turn a generalisation of the result of
Robinson \robtwo\ and the topological version of the result of Hof
\hof.

\sn{\bf Corollary 7.6} {\it With the conditions assumed in Theorem 7.4, the
pattern dynamical system $M\T$ is an almost 1-1 extension (2.15) of a minimal
$\R^d$ action by rotation on a $(N-\dim\Delta)$-torus.}

\sn{\bf Proof} It suffices to show that the central vertical arrow in the
diagram of Theorem 7.4 is 1-1 at some point. But this is immediate since
each of the end arrows is 1-1 at $u$ say. \qed

\newpage

\sn{\sect 8 Examples and Counter-examples}

\sn
In this section we give sufficient conditions,
similar to and stronger than 5.3, under which $P_u$ or $\T$ is
pointed conjugate to $\widetilde P_u$, and show why these conjugacies
are not true in general.

\sn{\bf Definition 8.1} For data $(E,K,u)$, define  $B_u =
\overline{(Q+u) \cap IntK}$ (Euclidean closure in $E^{\perp}$).

\sn{\bf Proposition 8.2} {\it Suppose that $E$, $K$ and $u \in NS$
are chosen so that $E \cap \Z^N = 0$ and $\Delta \cap
[(B_u-B_u)-(B_u-B_u)] = \{0\}$, then $R_u = \Z^N$. In this case,
therefore, $P_u$ is pointed conjugate to $\widetilde P_u$.}

\sn{\bf Proof} This follows from the fact, deduced directly from
the condition given, that if $v \in NS \cap (Q+u)$ and $a,b \in
P_v$, then we can determine $w-w'$ whenever $w,w' \in \widetilde
P_v$ are such that $a=\pi(w)$ and $b=\pi(w')$. Knowing the
differences of elements of $\widetilde P_v$ forces the position of
$\widetilde P_v$ in $\Sigma$ by the density of
$\pi^{\perp}(\widetilde P_v)$ in $B_u$. So we can reconstruct
$\widetilde P_v$ uniquely from $P_v$ and we have $R_u = \Z^N$.
\qed

\sn{\bf Corollary 8.3} {\it In the case that the acceptance domain
is canonical, the condition that the points $\pi(w) \ |\ w \in
\{-1,0,1\}^N$ are all distinct is sufficient to show that $R_u =
\Z^N$ for all $u \in NS$. In this case, therefore, $P_u$ is
pointed conjugate with $\widetilde P_u$. }

\sn{\bf Proof} The condition implies that $\Delta \cap [(K-K) -
(K-K)] = \{0\}$  and this gives the condition in the proposition
since $B_u \subset K$. \qed

\sn If we are interested  merely in the equation between $\T$  and
$\widetilde P_u$, then the case that the acceptance domain is
canonical also allows simple sufficient conditions weaker than
8.3, as we show below.

We observe first that the construction of \okd\ can be extended
to admit non-generic parameters, provided that we are comfortable
with ``tiles'' which, although they are convex polytopes, have no
interior in $E$ and are unions of faces of the tiles with
interior. We retain these degenerate tiles as components of our
``tiling'', {\it i.e.\ }  as units of a pattern, giving essential
information about the pattern dynamics. We call such patterns {\it
degenerate} canonical projection tilings.

We write $e_j$ with $1 \leq j \leq N$ for the canonical unit basis of $\Z^N$.

\sn{\bf Proposition 8.4} {\it In the case of a canonical acceptance domain,
the condition that no two points from $\{\pi(e_j) \ |\ 1 \leq j \leq N\}$ are
collinear is sufficient to show that $\T_u$, the canonical (but possibly
degenerate) tiling, is pointed conjugate to $\widetilde P_u$ for all $u \in
NS$.}

\sn{\bf Proof} We show that the conditions given imply that the
shape of a tile (even in degenerate cases) determines from which
face of the lattice cube it is projected. In fact we shall show
that if $I \subset \{1,2,...,N\}$ then knowing $\pi(\gamma^I)$ and
the cardinality of $I$ determines $I$ (we write $\gamma^I = \{
\sum_{i \in I} \lambda_ie_i \ |\ 0 \leq \lambda_i \leq 1, \
\forall\ i \in I\}$).

Suppose that $\pi(\gamma^I) = \pi(\gamma^J)$ and $I,J \subset
\{1,2,..,N\}$ are of the same cardinality. It is possible always
to distinguish an edge on the polyhedron $\pi(\gamma^I)$ which is
parallel to a vector $\pi(e_i)$ for some $i \in I$; and $i$ is
determined from this edge by hypothesis. The same is true of this
same edge with respect to $J$ and so $i \in J$ also.

Writing $I' = I \setminus \{i\}$ and $J' = J \setminus \{i\}$  we
deduce that $\pi(\gamma^{I'}) = \pi(\gamma^I) \cap ( \pi(\gamma^I)
- \pi(e_i)) = \pi(\gamma^J) \cap ( \pi(\gamma^J) - \pi(e_i)) =
\pi(\gamma^{J'})$. Now we can apply induction on the cardinality
of $I$, and deduce that $I' = J'$ and so $I=J$. Induction starts
at cardinality $1$ by hypothesis.

Now, given this correspondence between shape of tile and its
preimage under $\pi$, we reconstruct $\widetilde P_u$ from $\T_u$
much as we did in Proposition 8.2 above. To complete the argument
we must check that no other element of $M \widetilde P_u$ maps
onto $\T_u$ in $M\T_u$. But if there were such an element, then
the argument above shows that it cannot be of the form $\widetilde
P_{u'}$ for $u' \in NS$. Also, by Theorem 7.4, we deduce that the
map $M\widetilde P_u \lra M\T_u$ is $p$-to-$1$ with $p \geq 2$,
and so, using Lemma 2.2ii/, we find $\T_v$ with two preimages of
the form $\widetilde P_v$ and $\widetilde P_{v'}$. But this
contradicts the principle of the previous sentence. \qed

\sn The conditions of this proposition include all the
non-degenerate cases of the canonical projection tiling usually
treated in the literature (including the Penrose tiling), so from
the equation $M\T_u = \widetilde \Pi_u/\Z^N$, deduced from
Proposition 8.4 as a consequence, we retrieve many of the results
stated in section 3 of \le .

\sn Now we turn to conditions under which $R_u$ differs  from
$\Z^N$. We can extend the argument of 7.1 to give a geometric
condition for elements of $R_u$, of considerable use in computing
examples.

\sn{\bf Lemma 8.5} {\it Suppose that $E \cap \Z^N = 0$,  $u \in
NS$ and $v \in E^{\perp}$. Then $v \in R_u \cap E^{\perp}$ if and
only if $v+u \in NS$ and $v + (B_u + \Delta) = B_u + \Delta$.}

\sn{\bf Proof} We suppose that $u,v+u \in NS$.  Then $\overline{
\pi^{\perp}(\widetilde P_u)} = \overline{(Q+u) \cap IntK}$ and
$\overline{ \pi^{\perp}(\widetilde P_{u+v})} = \overline{(Q+u+v)
\cap IntK}$.

If $v \in R_u \cap E^{\perp}$, then, by Lemma 6.3,  we have
$\widetilde P_u + \Delta +v = \widetilde P_{u+v} + \Delta$. Also,
since by 7.1, $v$ is in $\widetilde \Delta$, we have $Q+u = Q+u+v$
and $\pi^{\perp}(\widetilde P_u) + \Delta + v =
\pi^{\perp}(\widetilde P_{u+v}) + \Delta$. Putting all these
together gives the required equality $v + (B_u + \Delta) = B_u +
\Delta$.

Conversely, if $v + (B_u + \Delta) = B_u + \Delta$,  then, by the
argument of 7.1, $v \in \widetilde \Delta$ and so, as above, $Q+u
= Q+u+v$ and $\pi^{\perp}(\widetilde P_u) + \Delta + v =
\pi^{\perp}(\widetilde P_{u+v}) + \Delta$. So, if $a \in
\widetilde P_u$, then there is a $b \in \widetilde P_{u+v}$ such
that $\pi^{\perp}(a) - \pi^{\perp}(b) \in \Delta -v$. However,
since $\pi^{\perp}$ is 1-1 on $\Z^N$, we can retrieve the set
$\widetilde P_u$ as the inverse $\pi^{\perp}$ image of
$\pi^{\perp}(\widetilde P_u) \cap IntK$ and similarly for
$\widetilde P_{u+v}$. This forces $a-b \in \Delta -v$ therefore,
and so $\pi(a) = \pi(b)$. Thus we see that $P_u = P_{u+v}$, as is
required to show that $v \in R_u$ \qed

\sn{\bf Example 8.6} By Corollary 7.2 and the fact  that $\Delta =
0$, the octagonal tiling is pointed conjugate to both its
projection and strip point patterns.

Also it is from Lemma 8.5 (and not 8.3) that we can  deduce that
$R_u = \Z^5$ for the (generalised) Penrose tiling for all choices
of $u \in NS$, and so the generalised Penrose tiling is pointed
conjugate to both its projection and strip point patterns.

\sn Lemma 8.5 makes it quite easy to construct  counter-examples
to the possibility that $R_u = \Z^N$ always.

\sn{\bf Example 8.7} We start with a choice of $E$ for which
$\Delta \neq 0$ and $\widetilde \Delta$ contains $\Delta$
properly. The $E$ used to construct the Penrose tiling is such an
example. As in the proof of Lemma 2.6, let $U$ be the real span of
$\Delta$ and $V$ the orthocomplement of $U$ in $E^{\perp}$. Choose
a closed unit disc, $I$, in $V$ and let $J$ be the closure of a
rectangular fundamental domain for $\Delta$ in $U$. Let $K = I +J
\subset E^{\perp}$ and note that $K$ has all the properties
required of an acceptance domain in this paper and that, by Lemma
5.1, we have $\Pi_u = \widetilde \Pi_u$.

With $u \nin V + \widetilde \Delta$ (equivalently $u \in NS$), the
rectilinear construction of $K$ ensures that $((Q+u) \cap Int K) +
\Delta$ is invariant under the translation by any element, $v$, of
$\widetilde \Delta$. Also $v+u \in NS$ since the boundary of
$((Q+u) \cap Int K) + \Delta$ in $Q+u$ is invariant under
translations by $\widetilde \Delta$. So, with the
characterisation of Lemma 8.5, this shows that $R_u =
\widetilde\Delta + \Z^N$ which is strictly larger than $\Z^N$.

By varying the shape of $J$ in this example, we can get $R_u/\Z^N$ equal to
any subgroup of $(\widetilde\Delta + \Z^N)/\Z^N$, and we can make it a
non-constant function of $u \in NS$ as well.

\sn{\bf Example 8.8} Take $N=3$ with unit vectors $e_1$, $e_2$ and $e_3$.
Let $L$ be the plane orthonormal to $e_1-e_2$ and let $E$ be a line in $L$
placed so that $E \cap \Z^3 = \{0\}$. Then $E^{\perp}$ is a plane which
contains $e_1-e_2$ and we have $\Delta = \{n(e_1-e_2) \ |\  n \in \Z\}$.

Write $e^{\perp}_1$, $e^{\perp}_2$, and $e^{\perp}_3$ for the
image under $\pi^{\perp}$ of $e_1$, $e_2$ and $e_3$ respectively.
Then $e^{\perp}_1+e^{\perp}_2$ and $e^{\perp}_3$ are collinear in
$E^{\perp}$ and they are both contained in $V$ (the continuous
subspace of $V + \widetilde \Delta =
\overline{\pi^{\perp}(\Z^3)}$). $\Delta$ is orthogonal to $V$ and
$e^{\perp}_1-e^{\perp}_2 = e_1-e_2$. However $\widetilde \Delta =
\{n(e_1-e_2)/2 \ |\  n \in \Z\}$ which contains $\Delta$ as an
index $2$ subgroup.

The set $K = \pi^{\perp}([0,1]^3)$ is a hexagon in $E^{\perp}$ with a
centre of symmetry. It is contained in the closed strip defined by two
lines, $V+a$ and $V+b$, where $b-a = e_1-e_2$, and it is in fact
reflectively symmetric around an intermediate line, $V+c$ where $c-a =
(e_1-e_2)/2$. The boundary of the hexagon on each of $V+a$ and $V+b$ is an
interval congruent to $e^{\perp}_3$ ({\it i.e.\ } a translate of
$\{te^{\perp}_3 \ |\  0\leq t \leq 1 \}$). The four other sides are intervals
congruent to $e^{\perp}_1$ or $e^{\perp}_2$, two of each. The vertices of
the hexagon are on $V+a$, $V+b$ or $V+c$, two on each.

The point of all this is that there is a choice of non-singular
$u$ (in $E^{\perp}$ without loss of generality) such that  $B_u =
\overline{(Q+u) \cap Int K}$ consists of the two intervals $K\cap(
V+a')$ and $K\cap( V+b')$, where $2a' = a+c$ and $2b'=b+c$ (we can
choose $u \in (V+a')\cap NS$ for example), and these intervals are
a translate by $\pm(e_1-e_2)/2$ of each other. Thus we deduce that
$B_u + \Delta = B_u + \Delta +v$ for all $v \in \widetilde
\Delta$.

Upon confirming that $v+u \in NS$ for all $v \in \widetilde
\Delta$  as well, we use 8.5 to show that $R_u = \Z^3 + \widetilde
\Delta$, which contains $\Z^3$ with index $2$.

\sn{\bf Remark 8.9}  We note that Example 8.8
is degenerate and Proposition 8.4 shows why this must be the case. However,
under any circumstances, there exist projection method tilings, in the sense
of 4.4, pointed conjugate to $P_u$ or to $\widetilde P_u$. The point here is
that these tilings will not necessarily be constructed by the special method
of Katz and Duneau.

Also, leaving the details to the reader, we mention that Example
8.8 and its analogues in higher dimensions are the only
counter-examples to the assertion $R_u = \Z^N$ when the acceptance
domain is canonical and $\Delta$ is singly generated (and here we
find always that $R_u/\Z^N$ is at most a cyclic group of 2
elements). When $\Delta$ is higher dimensional we have no concise
description of the exceptions allowed.

\bigskip

\sn{\sect 9 The topology of the continuous hull}

\sn Sections 6 and 7 tell us that the continuous hulls $MP_u$ and
$M\widetilde P_u$ are quotients of the same space $\Pi_u \equiv
\widetilde \Pi_u$. One advantage of this equality is that the
topology of $\widetilde \Pi_u$ is more easily described than
$\Pi_u$ {\it a priori}.

\sn{\bf Definition 9.1} Let $F$ be a plane complementary, but not
necessarily orthonormal, to $E$ and let $\pi'$ be the skew projection
(idempotent map) onto $F$ parallel to $E$.

Let $K' = \pi'(K)$.

Set $F_u^o = NS \cap (Q+u) \cap F$ and let $F_u$ be the $\overline
D'$-closure (Def.~4.2) of $F_u^o$ in $\widetilde \Pi_u$.

Note that, since $R_u \subset Q$, $F_u^o$ is invariant under
translation by $\pi'(r),\ r \in R_u$ and by extension $F_u$
supports a continuous $R_u$ action. $R_u$ acts freely on $F_u$
when $E \cap \Z^N = 0$ ({i.e.\ }  with $x \in F_u$, the equation
$gx = x$ implies $g=0$).

Similarly, $R_u$ acts on $E$ by translation by $\pi(r),\ r \in R_u$.

\sn{\bf Lemma 9.2} {\it With the data above, $F_u= {\widetilde
\mu}^{-1}(F \cap (Q+u))$ (c.f. (7.3)) and there is a natural
equivalence $\widetilde \Pi_u \equiv F_u \times E$ and a
surjection $\nu \colon F_u \lra ((Q+u) \cap F)$  which fits into
the following commutative square $$\matrix{ \widetilde \Pi_u &
\longleftarrow \!\! \longrightarrow & F_u \times E \cr
&&\cr\biggl\downarrow{}^{\widetilde \mu}&&\biggr\downarrow {}^{\nu
\times id}  \cr && \cr Q+u & \longleftarrow \!\! \longrightarrow
& ((Q+u)\cap F) \times E. \cr}$$ Moreover these maps are
$E$-equivariant where we require that $E$ acts trivially on $F_u$.
The set $\nu^{-1}(v)$ is a singleton whenever $v \in NS \cap F
\cap (Q+u)$.

The canonical action of $R_u$ on $\widetilde \Pi_u$ is represented
in this equivalence as the direct sum ({\it i.e.\ } diagonal) of
the action of $R_u$ on $F_u$ and $E$ described in (9.1).}

\sn{\bf Proof} This follows quickly from the observation that there is a
natural, $\overline D'$-uniformly continuous and $E$ equivariant
equivalence $NS \cap (Q+u) = F_u^o + E \equiv F_u^o \times E$, which
can be completed. \qed

\sn{\bf Definitions 9.3} Let $\A_u$ be the algebra of subsets
({\it i.e.\ } closed under finite union, finite intersection and
symmetric difference) of $F_u^o$ generated by the sets $(NS \cap
(Q+u) \cap K') + \pi'(v)$ as $v$  runs over $\Z^N$. It is clear
that this algebra is countable and invariant under the action of
$R_u$.

Write $C^*(\A_u)$ for the smallest $C^*$ algebra which contains the
indicator functions of the elements of $\A_u$.

Let $\Z\A_u$ be the ring (pointwise addition and multiplication) generated
by this same collection of indicator functions.

Let $CF_u$ be the group of continuous integer valued functions
compactly supported on $F_u$.

These three algebraic objects support a canonical $R_u$ action
induced by the action of $R_u$ on $F_u$ described in (8.1) and so
we define three $\Z[R_u]$ modules. As $\Z^N$ sits inside $R_u$,
this action can be restricted to a canonical subaction by $\Z^N$
thereby defining three $\Z[\Z^N]$ modules.

Let $\B_u = \{ \overline A \ |\  A \in \A_u\}$ where the bar refers to
$\overline D'$ closure in $F_u$.

\sn And finally, we give the main theorem of this section which is of
utmost importance for the remaining chapters.

\sn{\bf Theorem 9.4} {\it With data $(E,K,u)$ and $\widetilde \Pi_u = \Pi_u$,

a/ The collection $\B_u$ is a base of clopen neighbourhoods which generates
the topology of $F_u$.

b/ We have the $*$-isomorphisms of $C^*$-algebras $C_o(F_u) \cong
C^*(\A_u)$ and $C_o(\widetilde \Pi_u) \cong C^*(\A_u) \otimes C_o(E)$
which respect the maps defined in Lemma 8.2.

c/ $CF_u \cong \Z\A_u$ as a $\Z[R_u]$ module (and by pull-back as
a  $\Z[\Z^N]$ module).

d/ $F_u$ is locally a Cantor Set.}

\sn First we have a lemma also of independent interest in the next section.

\sn{\bf Definition 9.5} Write $\overline K$ for the $\overline D'$-closure
of the set $K' \cap NS \cap (Q+u)$ (recall $K'$ from 9.1).

\sn{\bf Lemma 9.6} {\it $\overline K$ is a compact open subset of $F_u$. }

\sn{\bf Proof} Closure is by definition so compactness follows immediately
on observing that $K' \cap (Q+u) \cap NS$ is embedded $\overline
D'$-isometrically in the space $M\widetilde P_u \times K'$ with metric $D'
+ \norm .$ as the closed subset $\{(\widetilde P_v, v) \ |\  v \in K' \cap
(Q+u) \cap NS\}$. But $M\widetilde P_u \times K'$ is compact.

For openness, we appeal to an argument similar to that of (6.1). Suppose,
for a contradiction, that $v_n$ is a $\overline D'$-convergent sequence in
$(F \cap NS ) \cap K'$ and that $v'_n$ is a $\overline D'$ convergent
sequence in $(F \cap NS ) \setminus K'$ and that both sequences have the
same limit $x \in \widetilde \Pi_u$. Therefore $v = \widetilde \mu(x)$ is
the Euclidean limit of the $v_n$ and $v'_n$ and so $v \in \partial K$. But
by construction $\widetilde P_{v_n}$ and $\widetilde P_{v'_n}$ have a
different $D'$ limit -- a contradiction since the limit in each case must be
$\widetilde \eta(x)$. \qed

\sn{\bf Proof of Theorem 9.4} a/
The sets in $\B_u$ are clopen by Lemma 9.6 above. Therefore the metric
topology on $F_u$ we are considering, let us call it $\tau$, is finer than the
topology $\tau'$ generated by $\B_u$. Both topologies are invariant under the
action of $\Z^N$ so that it suffices to show their equivalence on some
closed $r$-ball $X$ of $F_u$. Suppose that $\tau'$ were a Hausdorff topology.
Then this equivalence would follow
from the continuity of the identity map
$(X,\tau)\to (X,\tau')$, because compactness of $(X,\tau)$ implies that of
$(X,\tau')$ and so the map would automatically be a homeomorphism.
On the other hand, we will proof below that $A \mapsto \overline A$
yields an isomorphism of Boolean algebras between $\A_u$ and $\B_u$.
In particular, $\B_u$ is closed under symmetric differences so
that $\tau'$ being $T_1$ already implies that it is Hausdorff.

It remains therefore to check that $\tau'$ is $T_1$, i.e.\ that the
collection $\B_u$ contains a decreasing
set of neighbourhoods around any point in $F_u$.

Certainly, if $a\neq b$ with $a,b \in F \cap (Q+u)$, then the
assumption that $Int(K) \cap (Q+u) \neq \emptyset$ (2.6) (hence
$Int(K') \cap (Q+u) \neq \emptyset$, interior taken in $F$) and
the facts that $K'$ is  bounded and that $\pi'(\Z^N)$ is dense in
$Q \cap F$, imply that there is some $v \in \Z^N$ such that $a \in
(Int(K') \cap (Q+u)) + \pi'(v)$ and $b \nin \overline {(K' \cap
(Q+u))} + \pi'(v)$ (Euclidean closure in $F \cap (Q+u)$). {\it
i.e.\ } $a$ and $b$ are separated by the topology induced by
$\widetilde \mu(\B_u)$.

In particular, if $x,y \in F_u$ and $y \in \cap\{B \in \B_u \ |\ x
\in B\}$ then $\widetilde \mu(x) = \widetilde \mu(y)$.

But, if $x \neq y$ and $\widetilde \mu(x) = \widetilde \mu(y)=v$,
then, by (3.7)e/, $\widetilde \eta(x) \neq \widetilde \eta(y)$ and
we may suppose that there is a point $p \in \widetilde \eta(x)
\setminus \widetilde \eta(y)$. We use the argument of Lemma 8.6 to
show that then there are two sequences $v_n, v'_n \in \pi'(\Z^N)$
both converging to $v$ in Euclidean topology and such that
$p+(v_n-v) \in \Sigma$ and $p + (v'_n-v) \nin \Sigma$ for all $n$
large enough. But then, for such a choice of $n$,  $y \nin
\overline A$ (closure of $A = (NS \cap (Q+u) \cap K') +
\pi'(p+v_n)$ in $\overline D'$ metric) and $x \in \overline A$, a
contradiction to the construction of $y$.

Therefore $x=y$ and so, by the local compactness (Lemma 9.6) of $F_u$ we
have the required basic property of the collection $\B_u$.

b/ This will follow from a/ and the equivalences in Lemma 9.2 if we can
show that $\A_u$ is isomorphic to $\B_u$ as a Boolean algebra. To show
this, it is enough to show that $A \mapsto \overline A$ (closure in
$\overline D'$ metric) is 1-1 on $\A_u$; and for this it suffices to prove
that if $A \in \A_u$ is non-empty, then its Euclidean closure has interior
(in $(Q+u) \cap F$).

Note that $NS \cap K' = NS \cap Int(K')$, so that if $A \in \A_u$ then $A$
is formed of the union and intersection of sets of the form $(NS \cap (Q+u)
\cap Int(K')) + v$ ($v \in \pi'(\Z^N)$), and the subtraction of unions and
intersections of sets of the form $(NS \cap (Q+u) \cap K') + v$. With this
description and Lemma 2.5, $A$ is equal to $NS \cap Int(\overline A)$ (Euclidean
closure, and interior in $(Q+u) \cap F$), and from this our conclusion follows.

c/ Elements of $CF_u$ are finite sums of integer multiples of
indicator functions of compact open sets. Such sets are finite unions of
basic clopen sets from the collection in part a/. The isomorphism in part
b/ completes the equation.

d/ Given the results of a/ and Lemma 9.6 it is sufficient to show that
$F_u$ has no isolated points. However, by the argument of part b/ and Lemma 2.5 we see
that every clopen subset of $F_u$ has $\widetilde \mu$ image with
Euclidean interior (in $(Q+u) \cap F$)
and so cannot be a single point. \qed

\bigskip

\sn{\sect 10 A Cantor $\Z^d$ Dynamical System}

\sn
In this section we describe a $\Z^d$ dynamical system whose mapping torus is
equal to $M\T$.
First, assuming projection data, $(E,K,u)$ and
$E \cap \Z^N=0$, we find a
suitable $F$ to which to apply the construction of the previous section.

\sn{\bf Definition 10.1} Suppose that $G$ is a group intermediate
to $\Z^N$ and $R_u$. The examples in our applications ahead are
$G=\Z^N$, $G=R_u$ and $G = S_{\T}$, found in Theorem 7.4.

Fix a free generating set, $r_1,r_2,...,r_N$, for $G$ and suppose that the
first $\dim \Delta$ of these generate the subgroup $G \cap E^{\perp}$ (this
can be required by Lemma 2.9).

Let $F$ be the real vector space spanned by $r_1,r_2,...,r_n$, where $n = N
-d$.

Note that, since $E \cap G = 0$ (by Lemma 2.9 and the assumption $E \cap
\Z^N = 0$), $F$ is complementary to $E$ and, since $n \geq \dim \Delta$,
$F$ contains $\Delta$.

\sn Note that by Theorem 7.1, any two groups intermediate to $R_u$
and  $\Z^N$, $G$ and $G'$ say, differ only by some elements in
$R_u \cap E^{\perp}$, a complemented subgroup of $R_u$. Thus we
may fix their generating sets to differ only among those elements
which generate $G \cap E^{\perp}$ or $G' \cap E^{\perp}$
respectively. With this convention therefore, the construction of
$F$ is independent of the choice of group intermediate to $R_u$
and $\Z^N$ and $F$ depends only on the data $(E,K,u)$ chosen at
the start.

Also, more directly, the elements $r_{\dim \Delta + 1},...,r_N$
depend  only on the data $(E,K,u)$ chosen at the start.

\sn{\bf Definition 10.2} Suppose that $G_0$ is the subgroup of $G$
generated by $r_1,r_2,...,r_n$, and that $G_1$ is the
complementary subgroup generated by the other $d$ generators.
Since $n \geq \dim \Delta$, $G_1$  depends only on the data
$(E,K,u)$ and not on the choice of $G$ intermediate to $\Z^N$ and
$R_u$.

Both groups, $G_1$ and $G_0$ act on $F_u$ and $E$ as
subactions of $R_u$ (Definition 9.1).

Let $X_G = F_u/G_0$, a space on which $G_1$ acts continuously.

\sn{\bf Theorem 10.3} {\it Suppose that we have projection data
$(E,K,u)$  such that $E \cap \Z^N = 0$ and $G$, a group
intermediate to $\Z^N$ and $R_u$. Then $X_G$ is a Cantor set on
which $G_1$ acts minimally and there is a commutative square of
$G_1$ equivariant maps $$\matrix{ F_u & \bbra {q} & X_G \cr && \cr
\biggl\downarrow {}^{\nu} &&\biggr\downarrow {}^{\nu'}  \cr && \cr
F \cap (Q+u) & \llra & (F \cap (Q+u))/G_0. \cr}$$ The set
$\nu'^{-1}(v)$ is a singleton whenever $v \in (NS \cap F \cap
(Q+u))/G_0$.

The space $(F \cap (Q+u))/G_0$ is homeomorphic to a finite union of tori
each of dimension $(N-d-\dim \Delta)$. Indeed, this space can be considered
as a topological group, in which case it is the product of a subgroup of
$\widetilde \Delta / \Delta$ with the $(N-d-\dim \Delta)$-torus. The action
of $G_1$ on this space is by rotation and is minimal.}

\sn{\bf Proof} Assuming we have proved the fact that $X_G$ is compact then
the commuting square and its properties follow quickly. Therefore we look
at $X_G$.

Since $G_0$ acts isometrically on $F_u$ with uniformly discrete
orbits (Theorem 6.8), $q$ is open and locally a homeomorphism and
so $X_G$ inherits the metrizability of $F_u$, a base of clopen
sets and the lack of isolated  points (see Theorem 9.4 d/).

Now, let $Y_o = \{ \sum_{1 \leq j \leq n} \lambda_jr'_j \ |\  0
\leq \lambda_j < 1 \} \cap (Q+u)$, a subset of $F \cap (Q+u)$.
Choose $J \subset \Z^N$ finite but large enough that $Y_1 =
\cup_{v \in J}((K' \cap (Q+u)) + \pi'(v))$ contains $\overline
Y_o$ (Euclidean closure). In particular  $q( \cup_{v \in J}
(\overline K + \pi'(v))) = X_G$, the image of a compact set (Lemma
9.6) under a continuous map. So $X_G$ is also compact.

Therefore, we have checked all conditions
that show $X_G$ is a Cantor set.

Minimality follows from the minimality of the $G$ action on $F_u$ which in
turn follows from the minimality of the $\Z^N+E$ action on $\widetilde
\Pi_u$, proved analogously to (3.9).

The structure of the rotational factor system follows quickly from
the  first part of this lemma, the structure of $F \cap (Q+u)$,
and Lemma 9.2. \qed

\sn Now we describe $X_G$ as a fundamental domain of the action of $G_0$.

\sn{\bf Definition 10.4} From the details of 10.3 we constuct a clopen
fundamental domain for the action of $G_0$ on $F_u$. We let $G_0^+ = \{
\sum_{1\leq j \leq n} \alpha_jr_j \ |\  \alpha_j \in \Nat \}$ and set $Y^+ =
\cup_{r \in G_0^+} (Y_1 + r)$ and define $Y = Y^+ \setminus \cup_{r \in
G_0^+, r \neq 0} (Y^+ +r)$.

Define $Y_G = \overline {\nu^{-1}(Y \cap NS)}$ (closure in the
$\overline D'$ metric), a subset of $F_u$.

\sn The following is immediate from this construction, using Lemma 9.6 and
the equivariance of $\nu$ and $\nu'$ in Lemma 9.3 with respect to the $R_u$
action.

\sn{\bf Lemma 10.5} {\it With data $(E,K,u)$, $E \cap \Z^N = 0$,
and the  definitions above, $Y$ is a fundamental domain for the
translation action by $G_0$ on $F \cap (Q+u)$. Moreover, $Y_G$ is
a compact open subset of $F_u$, and a fundamental domain for the
action by $G_0$ on $F_u$.

There is a natural homeomorphism
$X_G \leftrightarrow Y_G$ which is $G_1$ equivariant. \qed}

\sn{\bf Definition 10.6} Define $CX_{G}$ to be the ring
of continuous integer valued functions defined on $X_G$.
Also define $C(F_u;\Z)$ to be the ring
of continuous integer valued functions defined on $F_u$
without restriction on support, uniformity or magnitude (c.f.\ the definition
of $CF_u$ from 9.3).
As $\Z[G_0]$
modules, the first is trivial and the second is defined as usual using the
subaction of the $R_u$ action on $F_u$. Both are non-trivial $\Z[G_1]$
modules.

\sn The following combines Lemmas 10.3, 10.5 and Proposition 9.4 and will be
of much importance in Chapter II.

\sn{\bf Corollary 10.7} {\it With the data of Lemma 10.5, $$CF_u \cong
CX_{G} \otimes \Z[G_0]$$ and $$C(F_u;\Z) \cong
Hom_{\Z}(\Z[G_0],CX_{G})$$ as $\Z[G_0]$ modules. \qed}

\sn{\bf Definition 10.8} Let $E'$ be the real span of
$r_{n+1},...,r_N$ selected in Definition 10.1. This space contains
$G_1$ (the integer span  $r_{n+1},...,r_N$) as a cocompact
subgroup.

Recall the definition of dynamical mapping torus \pt\ which for the $G_1$
action on $X_{G}$ may be equated with $$MT(X_{G},G_1) = (X_{G} \times E') /
\langle (gx,v) - (x,  v-g) \ |\  g \in G_1 \rangle .$$

\sn{\bf Proposition 10.9} {\it Suppose that we have data $(E,K,u)$ such that
$E \cap \Z^N = 0$, and $G$, a group intermediate to $\Z^N$ and $R_u$. With the
definitions above, $E'$ is a $d$-dimensional
subspace of $\R^N$ complementary
to both $F$ and $E^{\perp}$. Also $MT(X_{G},G_1)
\equiv \widetilde \Pi_u/G$.}

\sn{\bf Proof}  The transformation $E' \lra E$ defined by $r_j \mapsto
\pi(r_j)$, $n+1 \leq j \leq N$ is bijective since $G_1$ is complementary to
$G_0$ and hence to the subset $G \cap E^{\perp}$ of $G_0$. From this we
deduce the complementarity immediately.

{From} Lemma 10.3 we see that $G_0$ acts naturally on $MT(F_u,G_1)$ and that
$MT(X_{G},G_1)$ $ \equiv MT(F_u,G_1)/G_0$.

To form $MT(F_u,G_1)$ we take $F_u \times E'$ and quotient by the relation
$(ga,v) \sim (a,  v-g),\ g \in G_1 , a \in F_u, v \in \R^d$. Applying the
inverse of the map of the first paragraph, we can re-express the mapping
torus as $F_u \times E$ quotiented by the relation $(ga,w+\pi(g)) \sim (a,
w ),\ g \in G_1  , a \in F_u, w \in E$.

However, the action of $G_1$ on the $F_u$ is that induced by translation on
$F_u^o$ by elements $\pi'(g) \ |\  g \in G_1$. So, working first on the space
$F_u^o$, we have the equations $$MT(F_u^o,G_1) = (F_u^o \times E)/\langle
(a + \pi'(g),w+\pi(g)) - (a, w ) \ |\  g \in G_1 \rangle$$ and, since $F_u^o
\times E = NS$, we may write $ (a,w) =v \in NS$ and so the quotient equals
 $$NS/\langle v+g - v \ |\  g \in G_1 \rangle = NS/G_1$$
(recall that $NS$ is $G_1$ invariant
by Corollary 6.6). Then, by completing, we deduce the equation $MT(F_u,G_1)
= \widetilde \Pi_u / G_1$ directly. A further quotient by $G_0$ completes the
construction. \qed

\sn{\bf Corollary 10.10} {\it Suppose that $\T$ is a projection
method  pattern with data $(E,K,u)$ such that $E \cap \Z^N = 0$.
Then $(X_{\T},G_1)$ is a minimal Cantor $\Z^d$ dynamical system,
whose mapping torus $MT(X_{\T},G_1)$ is homeomorphic to $M\T$. The
pattern dynamical system, $(M\T, E)$ is equal to the canonical
$\R^d$ action on the mapping torus $(MT(X_{\T},G_1), \R^d)$ up to
a constant time change.}

\sn{\bf Proof} Choose $G = S_{\T}$ from Theorem 7.4 which gives
$M\T \equiv \widetilde \Pi_u/G$. From this, all but the time
change information  follows quickly from 10.9 and 10.3, noting
that $G_1 \cong \Z^d$.

To compare the two $\R^d$ actions, we apply the constant time change which
takes the canonical $\R^d (\cong E')$ action on $MT(X_{\T},\Z^d)$ to the
canonical $\R^d (\cong E)$ action on $M\T$ by the isomorphism
$\pi|_{E'}\colon  E' \lra E$, mapping generators of the $G_1$ action $r_j
\mapsto \pi(r_j)$ for $n+1 \leq j \leq N$. \qed

\sn{\bf Examples 10.11} The dynamical system of 10.10 for the octagonal
tiling is a $\Z^2$ action on a Cantor set, an almost 1-1 extension of a
$\Z^2$ action by rotation on $\torus^2$ (see \bcl ). For the Penrose
tiling, it is also a Cantor almost 1-1 extension of a $\Z^2$ action by
rotation on $\torus^2$ (see \robone ), where we must check carefully that the torus
factor has only one component (the only alternative of $5$ components is
excluded {\it ad hoc}).

\sn The correspondence in 10.9 and 10.10 respects the structures
found in  Theorem 7.4 and so we deduce an analogue.

\sn{\bf Definition 10.12} With data $(E,K,u)$ and the selection $G
= \Z^N$,  perform the constructions of 10.1 and 10.2, writing the
$G_0$ and $X_G$ obtained there as $G_u$ and $X_u$ respectively.

Similarly, with $G= S_{\T}$, write $G_0$ as $G_{\T}$ and $X_G$ as $X_{\T}$.

Also, with $G = R_u$, write $G_0$ as $\widetilde G_u$ and $X_G$ as
$\widetilde X_u$.

Note again that $F$ and $G_1$ are the same for all three choices of $G$.

\sn{\bf Corollary 10.13} {\it Suppose that we have data $(E,K,u)$
such  that $E \cap \Z^N = 0$. Then we can construct two Cantor
dynamical systems, $(X_u,G_1)$ and $(\widetilde X_u,G_1)$, the
latter a finite isometric extension of the former, together with a
compact abelian group, $M$, which is a finite union of $(N-d-\dim
\Delta)$-dimensional tori (independent of $u$) on which $G_1$ acts
minimally by rotation, and a finite subgroup, $Z_u$, of $M$.

These have the property that, if $\T$ is a projection method
pattern with data $(E,K,u)$, then there is a subgroup $Z_{\T}$ of
$Z_u$ and a commuting diagram of $G_1$-equivariant surjections
$$\matrix{\widetilde X_u & \llra & X_{\T} & \llra & X_u  \cr &&&&
\cr \biggl\downarrow {}^{ } && \biggr\downarrow &&
\biggr\downarrow {}^{ } \cr &&&& \cr M & \bbra {} & M/Z_{\T} &
\bbra { } & M/ Z_u \cr}$$ In this diagram, the top row consists of
finite isometric extensions, the bottom row of group quotients and
the vertical maps are almost 1-1.

Taking the $G_1$-mapping torus of this diagram produces the
diagram  of Theorem 7.4.}

\sn{\bf Proof} With the groups $G_u,G_{\T}$ and $\widetilde G_u$
defined in (10.12). By the remark after 10.1, we know that $G_u <
G_{\T} < \widetilde G_u$. Using the notation of sections 7, 9 and
10, set $ Z_u = \widetilde G_u/G_u \cong R_u/\Z^N$ and $Z_{\T} =
G_{\T}/G_u \cong S_{\T}/\Z^N$. Also set $M = ((Q+u)\cap F)/G_u$,
and note that, since $\Z^N = G_1 + G_u$, a direct sum, $G_u$ is in
fact independent of the choice of non-singular $u$. Thus $M$ is
independent also of $u$ (the uniform translation by $u$ is
irrelevant), so we attach no subscript.

With this notation and the equivalences above, together with the
description of the systems, $X_u$ and $\widetilde X_u$, using
Theorem 10.3, we complete the result. \qed

\newpage

\headline={\ifnum\pageno=1\hfil\else
{\ifodd\pageno\rightheadline\else\leftheadline\fi}\fi}
\def\rightheadline{\tenrm\hfil{\smallfont
II GROUPOIDS, $C^*$-ALGEBRAS, AND THEIR INVARIANTS}\hfil\folio}
\def\leftheadline{\tenrm\folio\hfil{\smallfont
FORREST HUNTON KELLENDONK}\hfil} \voffset=2\baselineskip
\sn{\chap II Groupoids, $C^*$-algebras, and their invariants}
\bigskip

\sn{\sect 1 Introduction}

\sn In this chapter we develop the connections between the pattern dynamical
systems described in Chapter I and the pattern groupoid. As with the
continuous hull, a pattern groupoid, which we write $\G\T^*$, can be defined
abstractly for any pattern, $\T$, of Euclidean space and we refer to \kelzero\
\keltwo\ for the most general definitions. Here we give a special form for
projection method patterns.

The (reduced) $C^*$-algebra, $C^*(\G\T^*)$, of this groupoid is a
non-commutative version of the mapping torus; this should be regarded as a more
precise detector of physical properties of the quasicrystal and the discrete
Schr\"odinger operators for the quasicrystal are naturally members of this
algebra.

The initial purpose of this chapter is to compare the
non-commutative structure, {\it i.e.\ }, the groupoid, of a
pattern with the dynamical systems constructed before. In this
regard, we cover similar ground to the work of Bellissard {\it et
al.\/} \bcl\ but, as noted in the introduction, applied to a
groupoid sometimes different. Our Theorem 4.2 is the complete
generalization of  their connection between $C^*$-algebraic
$K$-theory and dynamical cohomology, found for projection method
tilings in $2$ dimensions.

The dynamical system and $C^*$-algebras we associate to a pattern
give rise to  a range of ways of attaching an invariant to the
pattern $\T$. These include $C^*$ and topological $K$-theories,
continuous groupoid cohomology, \v Cech cohomology and the
dynamical or group (co)homology. The second main result of the
chapter is to set up and demonstrate that all these invariants are
isomorphic as groups, though the non-commutative invariants
contain the richer structure of an ordered group. This structure
appears likely to contain information relevant to subsequent
investigations, only some of which is recoverable from the other
invariants. On the other hand, as we shall see in later chapters,
the group (co)homology invariants admit greater ease of
computation and is often sufficient for other applications.

The final result of the chapter shows that all these common
invariants provide an  obstruction to the property of self
similarity of a pattern. We shall use this obstruction in Chapter
IV to show that almost all canonical projection method patterns
fail to be substitution systems.

The organisation of this chapter is as follows. In \S2 we define
the  various groupoids considered and their $C^*$ algebras, and
discuss their equivalences. In \S3 we consider the notion of
continuous similarity of topological groupoids. This is an
important equivalence relation for us as continuously similar
groupoids have the same groupoid cohomology. In \S4 we define our
invariants and prove them to be additively equivalent and in \S5
we establish the role our invariants play in discussing self
similarity.

\newpage

\sn{\sect 2 Equivalence of Projection Method pattern groupoids}

\sn First we develop some general results about topological
groupoids,  appealing to the definitions in \ren . These will lead
us to the notion of Equivalence of Groupoids which compares most
naturally the groupoid $C^*$-algebras.

Also in this section, we define several groupoids which can be
associated to a projection method pattern. We will show that many
of these groupoids are related by equivalence.

\sn{\bf Definition 2.1} We write the unit space of a groupoid $\G$
as $\G^o$, and write the range and source maps, $r,s \colon \G
\lra \G^o$ respectively. Both these maps are continuous and, due
to the existence of a  Haar System in all our examples, we note
that they are open maps as well.

Recall the {\it reduction} of a groupoid. Given a groupoid $\G$ with unit
space, $\G^o$, and a subset, $L$, of $\G^o$, define the reduction of $\G$
to $L$ as the subgroupoid ${}_L\G_L = \{ g \in \G \ |\  r(g), s(g) \in L \}$ of
$\G$, with unit space, $L$.

If $L$
is closed then ${}_L\G_L$ is a closed subgroupoid of $\G$.

We also define $\G_L = \{ g \in \G \ |\  s(g) \in L \}$ and note
the maps  $\rho \colon \G_L \lra \G^o$ and $\sigma \colon \G_L
\lra L$ defined by $r$ and $s$ respectively.

We say $L \subset \G^o$ is {\it range-open} if, for all open $U
\subset \G$,  we have $r(\{x \in U : s(x) \in L\})$ open in
$\G^o$.

Suppose a topological abelian group, $H$, acts by homeomorphisms
on a topological  space $X$, then we define a groupoid called the
{\it transformation groupoid}, $\G(X,H)$, as the topological
direct product, $X \times H$, with multiplication $(x,g)(y,h) =
(x,g+h)$ whenever $y = gx$, and undefined otherwise. The unit
space is $X \times\{0\}$.

\sn This last construction is sometimes called the transformation
group  \ren\ or even the transformation group groupoid, but we
prefer the usage to be found in \pat .

We note that if $H$ is locally compact, then naturally
$C^*(\G(X,H)) =  C_o(X) \rtimes H$, the crossed product \ren .

\sn{\bf Lemma 2.2} {\it Suppose that $H$ is an abelian metric
topological  group acting homeomorphically on $X$. Let $\G =
\G(X,H)$ be the transformation groupoid and suppose that $L$ is a
closed subset of $X \equiv \G^o$.

a/ If $H$ is discrete and countable, then $\G^o$ is a clopen
subset of  $\G$, and $L$ is range-open if and only if it is clopen
in $X$.

b/ If there is an $\epsilon > 0$ such that for all neighbourhoods,
$B \subset B(0,\epsilon)$, of $0$ in $H$ and all $A$ open in $L$,
we have $BA$ open in $X$, then $L$ is range-open.}

\sn{\bf Proof} Only part b/ presents complications.  Suppose that
$U$ is open in $X \times H$. We want to show that $r((L \times H)
\cap U)$ is open. Pick $x = r(y,h) \in r(L \times H \cap U)$ and
let $C \times (B+h)$ be a neighbourhood of $(y,h)$ inside $U$,
with $B$ sufficiently small. Then $A = s(C) \cap L$ is open in $L$
and $x \in (B+h)A = h(BA) $ an open subset of $X$ by hypothesis.
However, $(B+h)A \subset r((L \times H) \cap U)$ by construction,
and so we have found an open neighbourhood of $x$ in $r((L \times
H) \cap U)$ as required. \qed

\sn We use the constructions from \mrw\ \rei\ without comment.  In
particular, we do not repeat the definition of (strong Morita)
{\it equivalence} of groupoids or of $C^*$-algebras, which is
quite complicated. For separable $C^*$-algebras strong Morita
equivalence implies stable equivalence and equates the ordered
$K$-theory (without attention to the scale). All our examples are
separable.

\sn{\bf Lemma 2.3} {\it Suppose that $\G$ is a locally compact
groupoid and that $L \subset \G^o$ is a closed, range-open subset
which intersects every orbit of $\G$. Then ${}_L\G_L$ is
equivalent to $\G$ (in the sense of \mrw ) and the two $C^*$
algebras, $C^*({}_L\G_L)$ and $C^*(\G)$ are strong Morita
equivalent.}

\sn{\bf Proof} It is sufficient to show that $\G_L \bra {\rho}
\G^o$  is a left ($\G \bra {r,s} \G^o$)-module whose $\G \bra
{r,s} \G^o$ action is free and proper, and that $\G_L \bra
{\sigma} L$ is a right (${}_L\G_L \bra {r,s} L$)-module whose
${}_L\G_L \bra {r,s} L$ action is free and proper. In short,
${}_L\G_L$ is an abstract transversal of $\G$ and $\G_L$ a $(\G,
{}_L\G_L )$-equivalence bimodule from which we can construct the
$(C^*(\G),C^*({}_L\G_L))$-bimodule which shows strong Morita
equivalence of the two algebras directly, {\it c.f.} \mrw\ Thm
2.8.

The definition of these actions is canonical and the freedom and
properness  of the actions is automatic from the fact that $L$
intersects every orbit and from the properness and openness of the
maps $r,s$. Indeed all the conditions follow quickly from these
considerations except for the fact that $\G_L \bra {\rho} \G^o$ is
a left ($\G \bra {r,s} \G^o$)-module; and the only trouble here is
in showing that $\rho$ is an open map. However, this is precisely
the problem that range-openness is defined to solve. \qed

\sn Together with Lemma 2.2 above, this result gives a convenient
corollary  which unifies the r-discrete and non-r-discrete cases
treated separately in \ap .

\sn{\bf Corollary 2.4} {\it Suppose that $(X,H)$ and $L \subset X$
obey either  of the conditions of Lemma 2.2, then, writing $\G =
\G(X,H)$, $C^*({}_L\G_L)$ and $C^*(\G)$ are strong Morita
equivalent. \qed}

\sn Before passing to more special examples, we remark that there
is no  obstruction to the generalisation of results 2.2, 2.3 and
2.4 to the case of non-abelian locally compact group actions,
noting only that, for notational consistency with the definition
of transformation groupoid, the group action on a space should
then be written on the right.

Now we define a selection of groupoids associated with projection
method  patterns, all of them transformation groupoids.

\sn{\bf Definition 2.5} Given a projection method pattern, $\T$, with data
$(E,K,u)$,  fix $G =S_{\T}$ as the group obtained from $\T$ (Theorem I.7.4) so
that $M\T = \widetilde \Pi_u/G$. Recall the definitions of $X_{\T}= X_{S_\T}$,
$Y_{\T}= Y_{S_\T}$ and $G_1$
from I.10.2 and I.10.4.

We define in turn: $\G X_{\T} = \G(X_{\T},G_1)$, from the $G_1$
action on $X_{\T}$, and $\G F_{\T} = \G(F_u,G)$, using the action
of $G$ on $F_u$, both defined in (I.8.1).

Also define $\G \widetilde \Pi_{\T} = \G(\widetilde \Pi_u, E+G)$.

All but the last of these groupoids are r-discrete (see \ren ).

\sn{\bf Lemma 2.6} {\it Suppose that $\T$ is a projection method
pattern with data $(E,K,u)$ such that $E \cap \Z^N = 0$. The
groupoids $\G F_{\T}$ and $\G X_{\T}$ are each a reduction of $\G
\widetilde \Pi_{\T}$ to the  closed range-open sets, $F_u$ and
$Y_{\T}$.}

\sn{\bf Proof} It is clear that $\G F_{\T}$ is a reduction of  $\G
\widetilde \Pi_{\T}$ to $F_u$. To prove that this a range-open set
using Lemma 2.2, we take an open subset of $F_u$ and examine the
action of small elements of $E+G$ on it. Only the $E$ action
enters our consideration and then it is clear from (I.9.2) that if
$B$ is an open subset of $E$ and $A$ is an open subset of $F_u$,
then as topological spaces, $BA \equiv A \times B$, which is
clearly open in $\widetilde \Pi_u$.

Recall the notation of (I.10.12), and the homeomorphism  $X_{\T}
\leftrightarrow Y_{\T}$, found by Lemma I.10.5, putting
$G=S_{\T}$. In effect, this homeomorphism equates $X_{\T}$ with a
fundamental domain of the $G_{\T}$ action on $F_u$. This
homeomorphism is $G_1$-equivariant if we equate the $G_1$ action
on $Y_{\T} $ with the induced action of $G/G_0$ on $F_u/G_0 =
X_{\T} \equiv Y_{\T} $. But this is precisely the correspondence
needed to equate $\G(X_{\T}, G_1)$ with the reduction of $\G
F_{\T}$ to $Y_{\T} $ considered as a subset of the unit space of
$\G F_{\T}$. Thus $\G X_{\T}$ is the reduction of $\G \widetilde
\Pi_{\T}$ to $Y_{\T}$, and since $Y_{\T}$ is clopen in $F_u$ the
same argument as above shows that $Y_{\T}$ is closed and
range-open in $\widetilde \Pi_u$. \qed

\sn Now we define a groupoid connected more directly with the pattern, $\T$.

\sn{\bf Definition 2.7}  From Definition I.4.4, recall the
two maps $M\widetilde P_u \lra M\T \lra^* MP_u$
whose composition is $\pi_*$. Without confusion we name the second
(starred) map $\pi_*$ as well.

We also define a map $\eta_{\T}$ which is the composite
$\widetilde \Pi_u \bra {\widetilde \eta} M\widetilde P_u \lra
M\T$.

Note that $\eta(x) = \pi_*(\eta_{\T}(x))$ for all  $x \in
\widetilde \Pi_u$ and that, being a composition of open maps
(I.3.9), $\eta_{\T}$ is an open map.

Define the {\it discrete hull} of $\T$ as
$\Omega_{\T} = \{ S \in M\T \ |\  0 \in
\pi_*(S)\}$.

The {\it pattern groupoid}, $\G\T$, is the space $\{ (S , v) \in
\Omega_{\T} \times E \ |\    v \in \pi_*(S)\}$ inheriting the
subspace topology of $\Omega_{\T}\times E$. The restricted
multiplication operation is $(S', v')(S, v) = (S', v+v')$, if $S =
v'S'$,  undefined otherwise. The unit space is $\G\T^o = \{(S,0) \
|\ S \in \Omega_{\T}\}$, homeomorphic to $\Omega_{\T}$.

Also define $E^{\perp}_u = \widetilde\mu^{-1}(E^{\perp})$, a
space which is naturally homeomorphic to $F_u$; a correspondence
made by extending the application of $\pi^{\perp}$, inverted by
the extension of $\pi'$.

\sn{\bf Lemma 2.8} {\it Suppose that $\T$ is a projection method
pattern with data $(E,K,u)$ such that $E \cap \Z^N = 0$. The
groupoid $\G\T$ is isomorphic to a reduction of $\G \widetilde
\Pi_{\T}$ to a closed range-open set.}

\sn{\bf Proof} Let $L$ be a compact open subset of $E_u^{\perp}$
so that $\eta_{\T}(L) = \Omega_{\T}$ and $\eta_{\T}$ is 1-1 on
$L$. This can be constructed as follows. Define $L_o = \overline{
NS \cap K \cap (Q+u)}$ where the closure is taken with respect to
the $\overline D'$ metric - a clopen subset of $E^{\perp}_u$ by
(I.9.6). Let $L = L_o \setminus \cup\{ gL_o \ |\ g \in G \cap
E^{\perp}, g \neq 0\}$ (using the $G$ action on $\widetilde
\Pi_u$).

We claim that the reduction of $\G \widetilde \Pi_{\T} $ to $L$ is isomorphic
to the pattern groupoid defined above.

Suppose that $(x;g,v) \in \G \widetilde \Pi_{\T}$ and $x \in L$
and $(g+v)x \in L$, then $ 0 \in \eta(x) = \pi_*(\eta_{\T}(x))$
and $0 \in \eta((g+v)x)$.  But note that the action by $v \in E$
on $x \in NS$ is $vx = x-v$ and so $\eta((g+v)x) = \eta(gx) - v  =
\pi_*(\eta_{\T}(gx)) - v = \pi_*(\eta_{\T}(x)) - v$. Thus $0,v \in
\pi_*(\eta_{\T}(x))$ and the map $\psi\colon  (x;g,v) \mapsto
(\eta_{\T}(x),v)$ is well defined ${}_L \G \widetilde \Pi_{\T}{}_L
\lra \G\T$. The $E$ and $G$ equivariance of the maps used to
define $\psi$ show that the groupoid structure is preserved.

Conversely, if $0,-v \in \pi_*(\eta_{\T}(x))$, then there are, by
construction  of $L$, $g,g' \in G$ such that $gx, (g'+v)x \in L$.
Thus $(gx;g'-g,v) \in {}_L \G \widetilde \Pi_{\T}{}_L$ showing
that $\psi$ is onto. Also, the $g,g'$ are unique by the
construction of $L$ above, and so $\psi$ is 1-1. The continuity of
$\psi$ and its inverse is immediate, so we have a topological
groupoid isomorphism, as required.

Thus we have shown that $\G\T$ is isomorphic to a reduction of $\G
\widetilde  \Pi_{\T}$ to the set $L$ which is clearly closed.

Also, $L$ is a subset of $E_u^{\perp}$, transverse to $E$, so that the same argument as 2.6 shows that $L$ is range open.

It remains to show that $L$ hits every orbit of $\G \widetilde \Pi_{\T}$ and
for this it is sufficient to show that for any $x \in \widetilde \Pi_u$, $Gx \cap  (L \times E) \neq 0$ (where we exploit the equivalence:
$\widetilde \Pi_u \equiv E^{\perp}_u \times E$ (I.9.2). But this is
immediate from the fact that $L \times E$ is a clopen subset of $\widetilde
\Pi_u$ (I.9.2), and by minimality of the $G+E$ action on $\widetilde \Pi_u$
(as in I.3.9). \qed

\sn Combining the Lemmas above, we obtain the following.

\sn{\bf Theorem 2.9} {\it Suppose that $\T$ is a projection method
pattern with data $(E,K,u)$ such that $E \cap \Z^N = 0$. The
$C^*$-algebras  $C^*(\G\T)$, $C_o(F_u) \rtimes G$ and $C(X_{\T})
\rtimes G_1$ are strong Morita equivalent and thus their ordered
$K$-theory (without attention to scale) is identical. \qed}

\sn{\bf Remark 2.10}  We can compare the construction above with
the ``rope''  dynamical system constructed by the third author
\kelone\ exploiting the generalised grid method introduced by de
Bruijn \db . The rope construction actually shows that, in a wide
class of tilings including the canonical projection tilings, there
is a Cantor minimal system $(X,\Z^d)$ such that $\G(X,\Z^d)$ is a
reduction of $\G\T$. By comparing the details of the proof above
with \kelone\ it is possible to show directly that, in the case of
non-degenerate canonical projection tilings, the rope dynamical
system is conjugate to $(X_{\T},G_1)$.

\sn We note that the construction of Lemma 2.8 depends only on the
data  $(E,K,u)$ and on $G$ and from this we deduce the following.

\sn{\bf Corollary 2.11} {\it Suppose we specify projection data
$(E,K,u)$,  such that $E \cap \Z^N = 0$. Then, among projection
method patterns, $\T$, with these data, the dynamical system
$(M\T,E)$ determines $\G\T$ up to groupoid isomorphism.

Thus among projection method patterns with fixed data, the
dynamical  invariants are at least as strong as the
non-commutative invariants.\qed}

\sn Finally we reconnect the work of this section with the
original  construction of the tiling groupoid \kelzero .

\sn{\bf Definition 2.12} Recall the notation $A[r] = (A \cap B(r)) \cup
\partial B(r)$ {\it etc.} defined in I.3.1 and I.4.2, for $r \geq 0$ and
$A\subset \R^N$ or $E$. Given two closed sets, $A,A'$ define the distance
$D_o(A,A') = \inf\{1/(r+1) \ |\ r> 0, A[r] = A'[r] \}$.

As a metric this can be used to compare point patterns in $E$ or
$\R^N$  (as in I.3.1 and I.4.2), or decorated tilings in $E$ as
described in I.4.1.

Given a tiling, $\T$, of $E$, bounded subsets which are the closure of their
interior,
the construction of the discrete hull in \kelzero\ starts
by placing a single puncture generically in the interior of each tile
according to local information (usually just the shape, decoration and
orientation of the tile itself--the position of the puncture would then depend
only on the translational congruence class of the tile).
So we form the collection of punctures,
$\tau(\T)$ of $\T$, a discrete subset of points in $E$.

We consider the set $\Omega_{\T}^o =  \{ \T+x \ |\ 0 \in \tau(\T+x) =
\tau(\T)+x\}$, and define a modified hull, which we write $\Omega_{\T}^*$ in
this section, as the $D_o$ completion of this selected set of shifts of $\T$.

As in \kelzero\ we consider only tilings $\T$ which are
of {\it finite type} (or of finite pattern type as in \ap\ or of finite local
complexity) which means that for each $r$ the set
$\{\T'[r]|\T'\in \Omega_{\T}^o\}$ is finite. It is well known that canonical
projection tilings are finite type.

{}From this hull, we define the groupoid, $\G\T^*$ exactly as for
$\G\T$: $\G\T^* = \{ (S , v) \in \Omega_{\T}^* \times E \ |\    v
\in \tau(S)\}$ with the analogous rule for partial multiplication.

The assumption of local information dictates more precisely that
up to isomorphism $\G\T^*$ is independent of the punctures chosen
and the map $\tau$ is continuous, $E$-equivariant, and 1-1 from
$\Omega_{\T}^*$ with $D_o$ metric to the space of Delone subsets
of $E$ also with $D_o$ metric.

\sn{\bf Remark 2.13} Although phrased in terms of tilings, this
definition can in fact be applied to patterns as well, where the
idea of puncture becomes now the association of a point with each
unit of the pattern (I.4.1). In this case the finite type
condition is equated with the condition that $\tau(\T)$ is Meyer
(see \lag ), and this is sufficient to prove the analogues of all
the Lemmas below. However, we continue to use the language of
tilings and, since every projection method pattern is pointed
conjugate to a decorated finite type tiling (decorating the
Voronoi tiles for $P_u$ for example (I.7.5)), we lose no
generality in doing so.

\sn We note that when a projection method pattern $\T$ is in fact
a tiling, the two definitions of hull (2.7  and above) given here
seldom coincide nor do we obtain the same groupoids (but we note
the important exception of the canonical projection tiling in 2.16
ahead). The remainder of this section shows that, never-the-less,
the two groupoids, $\G\T$ and $\G\T^*$, are equivalent. We start
by comparing $D$ and $D_o$.

\sn{\bf Lemma 2.14} {\it Suppose that $\T$ is a finite type tiling
as above, then $\Omega_{\T}^o$ is precompact with respect to
$D_o$. Further $D$ and $D_o$ generate the same topology on
$\Omega_{\T}^o$.}

\sn{\bf Proof} The precompactness of $\Omega_{\T}^o$  is proved in
\kelzero .

For any two tilings, we have $D(\T,\T') \leq D_o(\T,\T')$  by
definition, and so the topology of $D_o$ is always finer than that
of $D$.

Conversely, as a consequence of the finite type condition of $\T$
there is a number $\delta_o < 1$ such that if $0 < \epsilon <
\delta_o$, then $\T+x, \T+x' \in \Omega_{\T}^o$ and $D(\T+x,
\T+x') < \epsilon$ together imply that $\T+x$ and $\T+x'$ actually
agree up to a large radius ($1/\epsilon -1$ will do) and we
conclude $D_o(\T+x, \T+x') < 2\epsilon$ as required.  \qed

\sn Consequently, $\Omega_{\T}^*$ is canonically a  subspace of
$M\T$ and we can consider its properties as such.

\sn{\bf Lemma 2.15} {\it With respect to the $E$  action on $M\T$,
both $\Omega_{\T}$ and $\Omega_{\T}^*$ are range-open.}

\sn{\bf Proof} With the notation of the proof of  Lemma 2.8,
$\Omega_{\T} = \eta_{\T}(L)$, where $L$ is a compact open subset
of $E_u^{\perp}$. As in 2.8, $L$ is range-open in $\widetilde
\Pi_u$  and, since $\eta_{\T}$ is an open $E$-equivariant map, we
deduce the same of $\eta_{\T}(L)$.

For $\Omega_{\T}^*$, we note that (as in \S1, 3.4 before) by the
finite type condition of $\T$, there is a number $\delta_o$ so
that, if $x,x' \in E$, $\T' \in M\T$ and $0 < \norm{x-x'} <
\delta_o$, then $D(\T'-x,\T'-x') \geq \norm{x-x'} /2$ and
$D_o(\T'-x,\T'-x') =1$. In particular, since $\Omega_{\T}^*$ is
$D_o$-compact and hence a finite union of radius $1/2$
$D_o$-balls, we deduce that the map $\Omega_{\T}^* \times
B(\delta_o/2) \lra M\T$, defined as $(\T',x) \mapsto \T'+x$, is
locally injective and hence open. From here the range-openness of
$\Omega_{\T}^*$ is immediate. \qed

\sn{\bf Theorem 2.16} {\it  If $\T$ is at once a tiling and a
projection method pattern with data $(E,K,u)$, such that $E \cap
\Z^N = 0$, then the tiling groupoid $\G\T^*$ as defined in
Def.~2.12 is equivalent to $\G\T$ (Def.~2.7). Thus the respective
$C^*$-algebras are strong Morita equivalent also.

If $\T$ is a non-degenerate canonical projection tiling there is
a puncturing procedure inducing an isomorphism between the two groupoids.}

\sn{\bf Proof} Recall the definition of transformation groupoid
and the  action of $E$ on $M\T$ and consider $\G(M\T,E)$. Using
Lemmas 2.3 and 2.15, it suffices to show that the groupoids $\G\T$
and $\G\T^*$ are each a reduction of $\G(M\T,E)$ to the sets
$\Omega_{\T}$ and $\Omega_{\T}^*$ respectively. But this is
immediate from their definition.

To treat the case of canonical projection tilings, we note that
its tiles are parallelepipeds, and that the point pattern consists
of their vertices. Hence if we fix a small generic vector then
adding this vector to each vertex gives exactly one puncture for
each tile. Translating an element of $\Omega_{\T}$ by that vector
therefore produces an element of $\Omega_{\T}^*$ if we use these
punctures to define the latter. The map so defined is clearly a
continuous bijection which is $E$-equivariant (with respect to the
restricted $E$-action). Therefore it induces an isomorphism
between $\G\T$ and $\G\T^*$. \qed

\bigskip

\sn{\sect 3 Continuous similarity of Projection Method pattern groupoids}

\sn The aim of this section is to compare our groupoids in a
second way: by {\it continuous similarity}. This gives most
naturally an equivalence between groupoid cohomology.

We will show that many of the groupoids we
associate with a projection pattern are related in this way also.
Further background facts
about groupoids and their cohomology and the idea of similarity
may be found in \ren.

\sn{\bf Definition 3.1} Two homomorphisms,
$\phi,\psi\colon\G\lra\H$ between topological groupoids are {\it
continuously similar} if there is a function,
$\Theta\colon\G^o \lra \H$ such that $$\Theta(r(x))\phi(x) =
\psi(x)\Theta(s(x)).$$ Two topological groupoids are {\it
continuously similar} if there exist homomorphisms
$\phi\colon\G\lra\H$, $\psi\colon\H\lra \G$ such that
$\Phi_{\G}=\psi\phi$ is continuously similar to $id_{\G}$ and
$\Phi_{\H}=\phi\psi$ is continuously similar to $id_{\H}$.

\sn In all our examples, produced by Lemma 3.3 ahead, the function
$\Theta$  is continuous but note that this is not required by
Definition 3.1. Our interest in this relation lies in the
following fact which we exploit in \S4; see \ren\ for the
definition of continuous cohomology $H^*(\G;\Z)$ of a topological
groupoid $\G$.

\sn{\bf Proposition 3.2} (\ren, with necessary alterations for the
continuous category) {\it If $\G$ and $\H$ are continuously
similar then $H^*(\G;\Z) = H^*(\H;\Z)$. \qed}

\sn It turns
out that the construction of continuous similarities follows
closely the reduction arguments in the examples that interest us.

\sn{\bf Lemma 3.3} {\it Suppose that $(X,H)$ is a free topological
dynamical system (i.e., $hx=x$ implies that $h$ is the identity),
with transformation groupoid $\G = \G(X,H)$, and that $L,L'$ are
two closed subsets of $\G^o$. Suppose there are continuous
functions, $\gamma\colon L\lra H$, $\delta\colon L'\lra H$ which
define continuous maps $\alpha\colon L\lra L'$ and $\beta\colon L'
\lra L$ by $\alpha x = \gamma(x)x$ and $\beta x =\delta(x)x$. Then
${}_L\G_L$ and ${}_{L'}\G_{L'}$ are continuously similar.}

\sn{\bf Proof} We construct the two homomorphisms, $\phi\colon
{}_L\G_L\lra {}_{L'}\G_{L'}$ and $\psi\colon{}_{L'}\G_{L'} \lra
{}_L\G_L$ by putting $\phi(x,g) = (\alpha x, \gamma(gx) + g
-\gamma(x))$ and $\psi(y,h) = (\beta y,\delta(hy)+g-\delta(y))$.

A quick check confirms that these are homomorphisms, and they are
both clearly continuous. Moreover, $\phi\psi(y,h) = (\alpha \beta
y,\gamma((\delta(hy)+h-\delta(y))\beta(y))+\delta(hy)+h-\delta(y)-\gamma(\beta
y)$, a rather complicated expression which can be simplified if we
note that $\gamma((\delta(hy)+h-\delta(y))\beta(y))
=\gamma((\delta(hy)+h)y) = \gamma(\beta hy)$, and define
$\sigma(y)$ to be the element of $H$ such that $\sigma(y)y =
\alpha \beta y$. Then $\sigma(y)=\delta(y)+\gamma(\beta y)$, by
definition, and so $\sigma(hy)=\gamma(h\beta(y))+\delta(hy)
=\gamma((\delta(hy)+h-\delta(y))\beta(y))+\delta(hy)$. This gives
$\phi \psi(y,h)=(\alpha\beta y,\sigma(hy)+h-\sigma(y))$.

It is now easy to see that $\phi\psi$ is continuously similar to
the identity on ${}_{L'}\G_{L'}$ using the transfer function,
$\Theta: L' \lra {}_{L'}\G_{L'}$ given by $\Theta(y) =
(\alpha\beta y, -\sigma(y))$. This $\Theta$ happens also to be continuous.

Reciprocal expressions give the similarity between $\psi \phi$ and
the identity on ${}_{L}\G_{L}$. \qed

\sn{\bf Remark 3.4}  Note that if $L' = \G^o$, then Lemma 3.3 can
be reexpressed in the following form. If $L$ is a closed subset of
$\G^o$ for which there is a continuous map $\gamma\colon\G^o\lra
H$ such that $\gamma(x)x\in L$ for all $x \in\G^o$, then
${}_{L}\G_{L}$ is continuously similar to $\G$. (The condition on
$L$ implies that $L$ intersects every $H$-orbit of $(\G^o,H)$, but
the converse is not true.)

We apply this lemma and remark in two ways as we examine
continuous similarities between the various groupoids of
section 2.

\sn{\bf Lemma 3.5} {\it Suppose that $\T$ is a projection
method pattern and write $\G$ for $\G \widetilde \Pi_{\T}$. If $L$
is a clopen subset of $F_u$ then ${}_{L}\G_{L}$ is continuously
similar to $\G$.}

\sn{\bf Proof} It suffices to find the function $\gamma$ in the
remark.

Pick an order $\succ$ on $G=S_{\T}$ in which every non-empty set
has a minimal element. The set $EL = \{vy : v \in E, y \in L\}$ is
naturally homeomorphic to $E \times L$  by Lem.~I.9.2, and hence
is clopen in $\widetilde \Pi_u$. By the minimality of the $E+G$
action on $\widetilde \Pi_u$, (which is proved analogously to
Lem.~I.3.9), we have $\cup_{h\in G} hEL=\widetilde \Pi_u$, so that
for each $x \in \widetilde \Pi_u$, there is a $\succ$-minimal
$h\in G$ such that $hx \in EL$. Let $\gamma_0(x)$ be this $g$, and
note that by the freedom and isometric action of $G$ and the
clopenness of $EL$, this function $x \mapsto \gamma_0(x)$ is
continuous and maps $\widetilde\Pi_u$ to $EL$.

Now, given $\gamma_0(x)x\in EL$, there is a unique $\gamma_1(x)\in
E$ such that $\gamma_1(x)\gamma_0(x)x\in L$, and it is clear that
$x \mapsto\gamma_1(x)$ is continuous as a map $\widetilde\Pi_u\lra
E$. The desired map $\gamma$ can now be taken as this
composite.\qed

\sn{\bf Lemma 3.6} {\it Suppose that $\T$ is a finite type
tiling for which we have chosen a puncturing and which is
also a projection method pattern. Then $\G\T$ and $\G\T^*$ are
continuously similar.}

\sn{\bf Proof} We construct the maps $\gamma$ and $\delta$ as
follows.

Recall the metric $D_o(A,B) = \inf \{1/(r+1) : A[r] = B[r]\}$ and
the argument of Lemma 2.14 which shows that $D$ and
$D_o$ are equivalent on each of the sets $\Omega_{\T}$ and
$\Omega_{\T}^*$. (Actually the argument refers only to the second
space, but the fact that $\pi_*\T$ is a Meyer pattern (see \lag\
and Remark 2.13) allows it to be applied directly to
the first space as well.)

We may assume without loss of generality that, for each $S \in
M\T$, each point of $\pi_*(S)$ is in the interior of a tile of $S$
(if not we shift all the tiles in $\T$ by a uniform short generic
displacement and start again equivalently).

Suppose that $S \in \Omega_{\T}$. We know that $0 \in \pi_*(S)$
and that by assumption there is a unique tile in $S$ which
contains the origin in its interior. This tile has a puncture at a
point $v$ say, and so $S-v \in \Omega_{\T}^*$. So we have defined
a map from $\Omega_{\T}$ to $E$, $\gamma\colon S \mapsto -v$ which
is clearly continuous with respect to the $D_o$ metric. Moreover
the map, $\alpha\colon S \mapsto S-v$ has range $\Omega_{\T}^*$.

Conversely, let $r$ be chosen so that every ball in $E$ of radius
$r$ contains at least one point of $\pi_*(\T) = P_u$. Consider the
sets $\pi_*(S) \cap B(r)$, as $S$ runs over $\Omega_{\T}^*$ and
note that there are only finitely many possibilities, i.e.\ the set
$J = \{\pi_*(S) \cap B(r) : S \in \Omega_{\T}^*\}$ is a finite
collection of non-empty finite subsets of $B(r)$. Furthermore, by
the continuity of $\pi_*$ on $\Omega_{\T}^*$ with respect to the
$D_o$ metric, for each $C \in J$, the set $\{S \in \Omega_{\T}^* :
\pi_*(S) \cap B(r) = C\}$ is clopen.

For each $C \in J$ choose an element $v = v(C) \in C$, and define
$\delta(S) = -v( \pi_*(S) \cap B(r))$; this is continuous by
construction. Then $S + \delta(S) \in \Omega_{\T}$. \qed

\sn Now, appealing to 2.6 and 2.8, we can gather the results of
this section into the following corollary.

\sn{\bf Corollary 3.7} {\it  Suppose that $\T$ is a
projection method pattern. Then $\G\T$, $\G \widetilde \Pi_{\T}$,
$\G F_{\T}$ and $\G X_{\T}$ are all continuously similar. If $\T$
is also a finite type tiling, then these are all
continuously similar to the tiling groupoid, $\G\T^*$, of \kelone
.}

\sn{\bf Proof} Lemmas (3.5) and (2.6) and (2.8) show that $\G\T$,  $\G
F_{\T}$ and $\G X_{\T}$ are all continuously similar to $\G
\widetilde \Pi_{\T}$. The second part is a restatement of (3.6).
\qed

\bigskip

\sn{\sect 4 Pattern cohomology and K-theory}

\sn
We are now in a position to define our topological invariants for
projection method patterns and prove their isomorphism as groups.

\sn{\bf Definition 4.1} For a projection method pattern, $\T$, in
$\R^d$, we define for each $m\in\Z$ the following groups.
\item{{\bf (a)}} $H^m(\G\T,\Z)$, the continuous groupoid cohomology
of the pattern groupoid $\G\T$;
\item{{\bf (b)}} $H^m(M\T)$, the \v Cech cohomology of the space $M\T$;
\item{{\bf (c)}} $H_{d-m}(G_{1},CX_{\T})$ and $H^m(G_{1} ; CX_{\T})$,
the group homology and cohomology of $G_1$ with coefficients
$CX_{\T}$ (I.10.6);
\item{{\bf (d)}} $H_{d-m}(S_{\T}; CF_u)$ and $H^m(S_{\T};C(F_u;\Z))$,
the group homology and cohomology of $S_{\T}$ with coefficients
$CF_u$ (I.9.3) or $C(F_u;\Z)$ (I.10.6);
\item{{\bf (e)}} $K_{d-m}(C^*(\G\T))$, the $C^*$ $K$-theory of $C^*(\G\T)$;
\item{{\bf (f)}} $K_{d-m}(C(X_{\T}) \times G_{1})$, the $C^*$ $K$-theory
of the crossed product $C(X_{\T}) \times G_{1}$;

\noindent and, for finite type tilings (2.12),

\item{{\bf (g)}} the continuous groupoid cohomology $H^m(\G\T^*, \Z)$.

\sn{\bf Theorem 4.2} {\it For a projection method pattern
$\T$ and for each value of $m$, the invariants defined in (4.1)(a)
to (d) are all isomorphic as groups. If $\T$ is also a
finite type tiling, then these are also isomorphic to that
defined in (4.1)(g).

The invariants defined in (4.1)(e) and (f) are each isomorphic as
ordered groups. Finally, all these invariants are related via
isomorphisms of groups such as $$K_m(C(X_{\T}) \times G_{1})
\cong \bigoplus_{j=-\infty}^{\infty}H^{m+d+2j}(G_{1};CX_{\T}).$$
These invariants are, in all cases, torsion free, and those in
parts (a) to (d) and (g) are non-zero only for integers $m$ in the
range $0\leq m\leq d$.}

\sn{\bf Proof} It is immediate from the definition \ren\ that if
$W$ is a locally compact space on which a discrete abelian group
$G$ acts freely by homeomorphisms then the continuous groupoid
cohomology $H^*(\G(W,G);\Z)$ is naturally isomorphic to the group
cohomology $H^*(G, C(W;\Z))$ with coefficients the continuous
compactly supported integer-valued functions on $W$, with
$\Z[G]$-module structure dictated naturally by the $G$ action on
$W$. This proves the equality of (a) (and (g) where appropriate)
with the cohomology versions of (c) and (d) from (3.7) and the
fact that $\G F_{\T}$ and $\G X_{\T}$ are transformation
groupoids.

By (I.10.10) $M\T$ is homeomorphic to the mapping torus
$MT(X_{\T},G_1)$ and as noted in \fh\ the \v Cech cohomology
$H^*(MT(X_{\T},G_1))$ is isomorphic to the group cohomology
$H^*(G_{1}, CX_{\T})$ (this is standard and follows, for example,
by induction on the rank of $G_1$ with the induction step passing
from $\Z^r$ to $\Z^{r+1}$ coming from the comparison of the
Mayer-Vietoris decomposition  of $MT(X_{\T},G_1)$ along one
coordinate with the long exact sequence in group cohomology coming
from the extension $\Z^r\rightarrow\Z^{r+1}\rightarrow\Z$). This
proves the isomorphism of (4.1)(b) with the cohomological
invariant of (4.1)(c).

The isomorphism of $H^m(G_{1},CX_{\T})$ with  $H_{d-m}(G_{1} ;
CX_{\T})$ is simply Poincar\'e duality for the group $G_{1} \cong
\Z^d$.

By Cor.~I.10.7, a decomposition of $CF_u$ as a
$\Z[G_{\T}]$ module is given by
$CX_{\T}\otimes\Z[G_{\T}]$ where $CX_{\T}$ is a trivial
$\Z[G_{\T}]$ module. Standard homological algebra now tells us
that $$\eqalign{H_p(S_{\T};CF_{u}) \cong &H_p(G_1\oplus
G_{\T};CX_{\T}\otimes\Z[G_{\T}])\cr \cong &H_p(G_{1}; CX_{\T})\cr}$$
establishing the isomorphism of (4.1)(c) and (d) in homology. A
similar argument also works in cohomology for $C(F_u;\Z)$.

The isomorphism of (4.1)(e) and (f) follows from the Morita
equivalence of the underlying $C^*$-algebras in (2.7) and the
isomorphism $$K_m(C(X_{\T}) \times G_{1}) \cong
\bigoplus_{j=-\infty}^{\infty}H^{m+d+2j}(G_{1};CX_{\T})$$ is one
of the main results of \fh.

The torsion-freedom of these invariants also follows from the
results of \fh, while the vanishing of the (co)homological
invariants outside the range of dimensions stated is immediate
from their identification with the (co)homology of the group
$G_1 \cong \Z^d$.

\sn We make one further reduction of the complexity of the
computation of these invariants. Recall first the construction of
I.2.9 and I.10.1, in particular the equation $F\cap
Q=V\oplus\widetilde\Delta$ splitting $F$ into continuous and
discrete directions, and in which $\pi'(\Z^N)$, the projection of
$\Z^N$ onto $F$ parallel to $E$, is dense. Recall also the map
$\widetilde\mu\colon\Pi_u\lra Q+u$ defined in (I.4.3) for each $u
\in NS$.

\sn{\bf Definition 4.3} The restriction of $\widetilde\mu$ to
$F_u$ is written $\nu\colon F_u\lra F\cap(Q+u)=(F \cap Q) +
\pi'(u)$; this map is $\pi'(\Z^N)$-equivariant and $|\nu^{-1}(v)|
= 1$ precisely when $v \in NS \cap F \cap (Q+u)$ (see Lemma I.9.2).

 Let $\Gamma_{\T}=\{v\in S_{\T}:\pi'(v)\in V\}$ and
$CV_u=\{ f \in CF_u:\nu(supp(f))\subset V+\pi'(u)\}$, where $supp$
refers to the support of the function. This is consistent with
setting $V_u = \{x \in F_u : \nu(x) \in V + \pi'(u)\}$ and taking
$CV_u$ as the continuous integer valued functions on $V_u$ with compact
support.
There is a natural decomposition $CF_u =
CV_u \otimes_{\Z} \Z[\widetilde \Delta]$.

\sn{\bf Lemma 4.4} {\it As a  subgroup of $S_{\T}=G_{\T} + G_1$
(the decomposition of I.10.1), $\Gamma_{\T}$ satisfies
$\Gamma_{\T} =(\Gamma_{\T} \cap G_{\T}) \oplus G_{1}$. Moreover,
$\Gamma_{\T}$ is complemented in $S_{\T}$ by a group
$\Gamma_\Delta$, naturally isomorphic to $\widetilde\Delta$.

With this splitting, the action of $\Z[S_{\T}]=\Z[\Gamma_{\T}]
\otimes\Z[\Gamma_\Delta]$ on $CF_u = CV_u \otimes_{\Z}
\Z[\widetilde \Delta]$ is the obvious one, and hence there is an
isomorphism of homology groups  $H_*(S_{\T}; CF_u)\cong
H_*(\Gamma_{\T}; CV_u)$.}

\sn{\bf Proof} The decompositions and restrictions on $G_{\T}$ and
$G_1$ follow from the definition and the original construction of
(I.10.1). The conclusion in homology is the same
homological argument as used in the previous proof.\qed

\sn  We note that since $S_{\T} \supset \Z^N$ and
$\pi^{\perp}(\Z^N)$ is dense in $Q \cap F$, the group
$\Gamma_{\T}$ acts minimally on $V$ and hence on $V_u$.

\sn{\bf Corollary 4.5} {\it With the data above, $$K_n(C^*(\G\T))
= \bigoplus_{j=-\infty}^{\infty}H_{n+2j}(\Gamma_{\T}
;CV_u).\eqno\bullet$$}

\sn This is, in fact, the most computationally efficient route to
these invariants and, with the exception of Chapter III, the one we
shall use in the remainder of this memoir.

\bigskip

\sn{\sect 5 Homological conditions for self similarity}

\sn To motivate the direction we now move in, we give an immediate
application of these invariants. In this section  we show that the
(co)homological invariants defined in \S4 provide an obstruction
to a pattern arising as a substitution system. This result will be
used in Chapter~IV to show that almost all canonical projection
tilings fail to be self similar.

We adopt the construction of substitution tilings in \ap\
considering only finite type tilings whose
tiles are compact subsets of $\R^d$ homeomorphic to the closed ball.
A substitution procedure as in \ap\ is based on a map which assigns to
each tile of a tiling $\T$
a patch of tiles (a tiling of a compact subset of $\R^d$) which
covers the same set as the tile. Moreover,  the map is
$\R^d$-equivariant in the sense that translationally congruent
tiles are mapped to translationally congruent patches. This is
useful if there is a real constant $\lambda>1$ such that the
tiling which is obtained from $\T$, by first replacing all tiles
with their corresponding patches and then stretching the resulting
tiling by $\lambda$ (keeping the origin fixed), belongs to $M\T$.
This procedure of replacing and stretching can then be applied to
all tilings of $M\T$ thus defining the substitution map
$\sigma:M\T\to M\T$ which is assumed to be injective. As an aside,
Anderson and Putnam show that under this condition $\sigma$ is a
hyperbolic homeomorphism. They
 establish the following property for substitution tilings, i.e.\ tilings
which allow for a substitution map.

\sn{\bf Theorem 5.1   \ap} {\it Suppose that $\T$ is a
finite type substitution tiling of $\R^d$.
Then $M\T$ is the inverse limit of a
stationary sequence
$$Y \bla {\gamma} Y
\bla{\gamma}\cdots$$ of a compact Hausdorff space $Y$ with a continuous map
$\gamma$. \qed}

\sn Let us describe $Y$. Call a collared tile a tile of the tiling $\T$
decorated with the information of what are its adjacent tiles and a
collared prototile the
translational congruence class of a collared tile. As a topological space
we identify the collared prototile with a tile it represents
as a compact subset of $\R^d$, regardless of
its decoration. By the finite
type condition there are only finitely many collared prototiles.
Let $\tilde Y$ be the disjoint union of all collared prototiles.
$Y$ is the quotient of $\tilde Y$ obtained upon
identifying two boundary points $x_i$ of collared prototiles $t_i$, $i=1,2$,
if there are two collared tiles $\hat t_i$ in the tiling,
$\hat t_i$ of class $t_i$,
such that $\hat x_1=\hat x_2$ where $\hat x_i$ is the point on the boundary of
$\hat t_i$ whose position corresponds to that of $x_i$ in $t_i$.
If the tiles are polytopes which
match face to face (i.e.\ two adjacent tiles touch at complete faces)
then $Y$ is a finite CW complex whose highest dimensional cells
are the (interiors of) the tiles represented by the collared
prototiles and whose lower dimensional cells are given by quotients of the
set of their faces of appropriate dimension. The map $\gamma$ above is
induced by the substitution. Going through the details of its construction
in \ap\ one finds that in the case where tiles are polytopes which
match face to face
it maps faces of dimension $l$ to faces of dimension $l$ thus providing
a cellular map. We are then in the situation assumed for the second part of
\ap.

\sn{\bf Corollary 5.2} {\it Suppose that $\T$ is a finite type
tiling whose tiles are polytopes which
match face to face. Then for each $m$, the rationalized \v
Cech cohomology $H^m(M\T) \otimes \Q$ has finite $\Q$-dimension.}

\sn{\bf Proof} Like \ap\ we obtain
from (5.1) and the conclusion that $Y$ is a finite CW complex
and $\gamma$ a cellular map that
$H^m(M\T) = \dirlim H^m(Y)$. So
$H^*(M\T)\otimes\Q=\dirlim (H^*(Y)\otimes\Q)$. Thus the $\Q$
dimension of $H^*(M\T)\otimes\Q$ is bounded by that of
$H^*(Y)\otimes\Q$ and this is finite since $Y$ is a finite CW
complex.\qed

\sn The conclusion of (5.2) applies to much more general
pattern constructions. Note that the only principle used is that
the space $M\T$ of the tiling dynamical system is the inverse
limit of a sequence of maps between uniformly finite CW complexes.
We sketch a generalization whose details can be reconstructed by
combining the ideas to be found in \pr\ and \fortwo.

Finite type substitution tilings are analysed
combinatorially by Priebe in her PhD Thesis \pr\ where the  useful
notion of {\it derivative tiling}, generalized from the
1-dimensional symbolic dynamical concept \du, is developed. We do
not pursue the details here except to note that the derivative of
an almost periodic finite type tiling is almost
periodic and finite type, and that the process of
deriving can be iterated.

Suppose that $\T$ is an almost periodic finite type
tiling. By means of repeated derivatives, and adapting the
analysis of \fortwo\ for periodic lattices, we may build a
Bratteli diagram, $\B$, for $\T$. Its set of vertices at level $t$
is formally the set of translation classes of the tiles in the
$t^{\rm th}$ derivative tiling, and the edges relating two
consecutive levels, $t$ and $t+1$ say, are determined by the way
in which the tiles of the $(t+1)^{\rm th}$ derivative tiling are
built out of the tiles of the $t^{\rm th}$ derivative tiling. The
diagram $\B$ defines a canonical dimension group, $K_0(\B)$.
Adapting the argument of \fortwo\ we can define a surjection
$K_0(\B)\lra H^d(M\T)$ and hence a surjection $K_0(\B)\otimes
\Q\lra H^d(M\T)\otimes\Q$.

In \pr\ it is shown that the repeated derivatives of a finite type
aperiodic substitution tiling have a uniformly bounded number of
translation classes of tiles and are  themselves finite type. In
particular, the number of vertices at each level of its Bratteli
diagram $\B$ is bounded uniformly. Thus $K_0(\B)\otimes\Q$ is
finite dimensional over $\Q$, being a direct limit of uniformly
finite dimensional $\Q$ vector spaces, and so we reprove (5.2) for
the case $m=d$ .

It is worth extracting the full power of this argument since it
applies to a wider class than the substitution tilings.

\sn{\bf Theorem 5.3} {\it Suppose that $\T$ is a finite type
tiling of $\R^d$ whose repeated derivatives have a
uniformly bounded number of translation classes of tiles, then
$H^d(\G\T)\otimes\Q$ is finite dimensional over $\Q$. \qed}

\sn Therefore in Chapter IV, when  we take a pattern, $\T$,
compute its rationalized homology $H_0(\Gamma_{\T}; CV_u) \otimes
\Q$ and find it is infinite dimensional, we know we are treating a
pattern or tiling which is outwith the class specified in Theorem
5.3, and {\it a fortiori} outside the class of finite type
substitution tilings.

\newpage


\headline={\ifnum\pageno=1\hfil\else
{\ifodd\pageno\rightheadline\else\leftheadline\fi}\fi}
\def\rightheadline{\tenrm\hfil{\smallfont
III COHOMOLOGY FOR CODIMENSION 1}\hfil\folio}
\def\leftheadline{\tenrm\folio\hfil{\smallfont
FORREST HUNTON KELLENDONK}\hfil} \voffset=2\baselineskip

\sn{\chap III Approaches to calculation I:}
\medskip
\noindent{\chap Cohomology for codimension 1}
\bigskip

\sn{\sect 1 Introduction}

\sn Our goal in this chapter is to demonstrate the computability
of the invariants introduced in Chapter II and we do so by looking
at the case where $N=d+1$. In this case the lattice $\Delta$ is
always trivial whenever $E \cap \Z^N = 0$.

Recall that when $\Delta=0$ the projection pattern is determined by a small
number of parameters -- the dimensions $d$ and $N$ of the space $E$ and
the lattice $\Z^N\subset \R^N$, the slope of $E$ in $\R^N$ and the
shape of the acceptance domain $K$. We shall restrict ourselves to specific
acceptance domains in later chapters, but the main result of this chapter
(3.1) characterises the invariants of patterns on $\R^d$ arising as projection
from $\R^{d+1}$ for more or less arbitrary acceptance domains $K$.

This chapter thus gives a complete answer to the so-called
codimension 1 patterns ({\it i.e.}, $N-d=1$). After restricting
the shape of the acceptance domain we shall give in Chapter V an
alternative analysis of this case together with formul\ae\ for the
ranks of the invariants in the codimension 2 and 3 cases, when
these are finite.

When $N-d=1$ one of the most important features of $K$ is its
number of path components. To facilitate our computations we
examine in \S2 a general technique which sometimes simplifies the
computations of projection pattern cohomology when the acceptance
domain is disconnected. We prove our main results in \S3.

We note that the case $d=1$ gives the classical Denjoy
counter-example systems whose ordered cohomology is discussed in
\pss, an observation which has also been made by \her. The result
in this chapter, for $d=1$ and $H_*$ finitely generated, can be
deduced from that paper \pss\ a task which has been carried out in
\her.

\sn{\sect 2 Inverse limit acceptance domains}

\sn Suppose that $K$ and $K_i$, $i=1,2,\ldots\,$, are compact
subsets of $E^{\perp}$ each of which is the closure of its
interior, and suppose that $Int K = \cup_i Int K_i$ is a disjoint
union and that $\del K = \cup \del K_i$. Let $K_i^* = \cupij K_j$,
so that $Int K_i^* = \cupij Int K_j$ is a disjoint union and $\del
K_i^* = \cupij \del K_j$.

We define $NS^i = \R^N \setminus (E + \Z^N + \del K_i^*)$ and
$\Sigma^i = K_i^* + E$. So for each $ u \in NS^i$ we have
$\widetilde P_u^i = \Z^N \cap \Sigma^i$ and $P_u^i =
\pi(\widetilde P_u^i)$. From these we construct $\widetilde
\Pi_u^i$, $\Pi_u^i$, $MP_u^i$ and so on, as usual. In fact, in the
following, we shall be interested only in the strip pattern
$\widetilde P^i_u$.

Provided $u$ is non-singular for $K$, and hence is non-singular
for all $K_i^*$, we can take a space $F$ complimentary to $E$ in
$\R^N$ and a corresponding group $G_u$ (Def.~I.10.2)
which will play their usual
roles for all sets of projection data $(E,\R^N,K_i^*)$ and
$(E,\R^N,K)$. For each domain, $K_i^*$, we construct the
corresponding $F^i_u$ {\it etc}. The following lemma follows
easily from the definitions.

\sn{\bf Lemma 2.1} {\it Suppose $j < i$ is fixed throughout the statement
of this lemma. Then $NS^i$ is a dense subset of $NS^j$ and $NS =
\cap_k NS^k$. Moreover, for $u\in NS^i$, we have a natural
continuous $E+\Z^N$ equivariant surjection $\widetilde \Pi_u^i
\lra \widetilde \Pi_u^j$, and a natural continuous $E$ equivariant
surjection $M\widetilde P_u^i \lra M\widetilde P_u^j$; this latter
is described equivalently by the formula $S \mapsto S \cap
\Sigma^j$.

We also have an $\Z^N$-equivariant map $F_u^i \lra F_u^j$, and a
$G_{u}$-equivariant map $X_u^i \lra X_u^j$. All these maps respect
the commutative diagrams of Chapter I and they map many-to-one
only when the image is in (the appropriate embedding of) $NS^j
\setminus NS^i$.}\qed

\sn{\bf Theorem 2.2} {\it With the notation and assumptions above,
we have the following equivariant homeomorphisms.
\item{{\bf (a)}} $\widetilde \Pi_u\cong\invlim\widetilde
\Pi_u^i$,\ \ \ $E+\Z^N$ equivariantly;
\item{{\bf (b)}} $M \widetilde P_u\cong\invlim M \widetilde
P_u^i$,\ \ \  $E$-equivariantly;
\item{{\bf (c)}} $F_u\cong\invlim F_u^i$,\ \ \ $\Z^N$-equivariantly;
\item{{\bf (d)}} $X_u\cong\invlim X_u^i$,\ \ \ $G_{u}$-equivariantly.
}

\sn{\bf Proof} Once again the results are straightforward from the
definitions. The map $M \widetilde P_u \lra M \widetilde P^i_u$ is
again equivalently written $S \mapsto S \cap \Sigma^i$. \qed

\sn The following is now a direct consequence of (4.2), (2.2)(b)
and the continuity of \v Cech cohomology on inverse limits.

\sn{\bf Corollary 2.3} {\it There is a natural equivalence
$H^*(\G\widetilde P_u)= \dirlim H^*(\G\widetilde P_u^i)$. \qed}

\bigskip

\sn{\sect 3 Cohomology of the case d=N--1}

\sn
In this section we determine the cohomology for
projection method patterns when $d=N-1$. It is the only case for
which we have such a complete answer.

Here $E$ is a codimension $1$ subspace of $\R^N$ and so $E \cap
\Z^N=0$ implies that $\Delta =0$.
Therefore, $MP_v = MP_u = M\widetilde P_u$ for all $u,v \in NS$ and
so, given $E$ and $K$, there is only one projection pattern torus
$MP$ to consider, no need to parametrise by $u$, and an
equation $S_{P} = \Z^N$
(I.8.2). With this in mind, we shall avoid further
explicit mention of any particular non-singular point $u$.

Write
$e_1,...,e_N$ for the usual unit vector basis of $\R^N$, which are
also the generators of $\Z^N$. Choose the space $F$ as that
spanned by $e_N$, and so $G_{\T} = \langle e_N \rangle$ and $G_{u}
= \langle e_1,e_2,...,e_{N-1} \rangle$. Recall that we write $K' =
\pi'(K) \subset F$, where $\pi'$ is the skew projection onto $F$
parallel to $E$ and that $\pi'$ maps $K$ homeomorphically to $K'$,
preserving the boundary, $\del K' = \pi'(\del K) \cong \del K$.

Now any compact subset of $E^{\perp} \equiv \R$ which is the
closure of its interior is a countable union of closed disjoint
intervals; and $K$ is such a set. Thus $\del K$ and hence $\del
K'$ is countable. Pick $A=\{p_1,p_2,...\}$, $p_j\in\del K'$, to be
a set of representatives of $\pi'(\Z^N)$ orbits of $\del K'$, a
countable and possibly finite set. Write $k\in\Z_+\cup\infty$ for
the cardinality of $A$.

\sn{\bf Theorem 3.1} {\it If $\T$ is a projection method pattern
with $d = N-1$ and $E\cap\Z^N=0$, then $$H^m(\G\T)= H^m(\torus^N
\setminus k\hbox{ points }) = \left\{\eqalign{ \Z^{{N \choose
m}}\qquad&\hbox{ for }m\leq N-2,\cr \Z^{N+k-1}\qquad&\hbox{ for }m
= N-1,\cr 0\quad\qquad&\hbox{ otherwise. }\cr}\right.$$ An
infinite superscript denotes the countably infinite direct sum of
copies of $\Z$.}

\sn{\bf Proof} We know that $Int K'$ is the union of a countable
number of open intervals, whose closures, $K_i$, are disjoint. We
use the notation and results of \S2, setting $K_i^* = \cupij K_j$
as the finite union of disjoint closed intervals, $\cupij
[s_je_N, r_je_N]$ say. As $M\T = M \widetilde P$, by (2.3) it is
enough to compute the direct limit $\dirlim H^*(M\widetilde P^i)$.

We consider the process of completion giving rise to the space
$M\widetilde P^i$ which we consider as $\widetilde\Pi^i$ (the
completion of the non-singular points $NS^i$) modulo the action of
$\Z^N$. The limit points introduced in $\widetilde\Pi^i$ arise as
the limit of patterns $P^i_{x_n}$ as $x_n$ approaches a singular
point, either from a positive $e_N$ direction, or from a negative
one. To be more precise, suppose that $x_n = x + t_ne_N \in NS$ is
a sequence converging to $x \in \R^N$ in the Euclidean topology.
If $(t_n)$ is an increasing sequence, then
${\displaystyle\mathop{\rm lim}_{n\rightarrow\infty}} \widetilde
P_{x_n}^i$ exists in the $D$ metric and is the point pattern $(x +
\Z^N) \cap ( \cupij (s_je_N, r_je_N]+ E)$.

Likewise, if $(t_n)$ is a decreasing sequence then
${\displaystyle\mathop {\rm lim}_{n\rightarrow\infty}}\widetilde
P_{x_n}^i$ is the point pattern $(x+\Z^N) \cap ( \cupij [s_je_N,
r_je_N) + E)$. These two patterns are the same if and only if $x
\in NS^i$. If $x \nin NS^i$ then these two patterns define the two
$D$ limit points in $\widetilde\Pi^i$ over $x\in\R^N$. Thus the
quotient $M\widetilde P^i \lra \torus^N$ is 1-1 precisely when
mapping to the set $NS^i/\Z^N$, and otherwise it is 2-1; we can
picture the map intuitively as a process of ``closing the gaps''
made by cutting $\torus^N$ along the finite set of hyperplanes
$(\del K_i^* + E)/\Z^N$, {\it c.f.} \le.

We examine the space $M\widetilde P^i$ in more detail. Given $r
>0$, consider the space $M^i_r=\{(S \cap B(r)) \cup\del B(r) : S
\in M \widetilde P_u^i\}$ endowed with the Hausdorff metric $d_r$
among the set of all closed subsets of $B(r)$, the closed ball in
$\R^N$ with centre $0$ and radius $r$. By construction $M^i_r$ is
compact and, for $s\geq r$ and $i\geq j$, there are natural
restriction maps $M^i_s \lra M^j_r$, whose inverse limit for fixed
$i=j$ is $M\widetilde P^i$, and whose inverse limit over all $i$
and $r$ by (2.2) is $M\widetilde P$. The map $M_r^i\lra\torus^N$
given by $\widetilde P_v\mapsto v \mod \Z^N$, $v\in NS^i$, factors
the canonical quotient $M\widetilde P^i\lra\torus^N$.

Define $C^i_r$ as the set $\{v\in \torus^N : (v +\Z^N) \cap (\del
\Sigma^i\cap Int B(r)) \neq \emptyset \}$. As before, $M_r^i \lra
\torus^N$ is 2-1 precisely on those points mapped to $C_r^i$ and
otherwise is 1-1.

The intersection $\del \Sigma^i \cap B(r)$ is, for all $r$ large
enough compared with the diameter of $K_i^*$, equal to a finite
union of codimension $1$ discs, parallel to $E$, and of radius at
least $r-1$, and at most $r$. Each of these discs has centre
$\pi^{\perp}(a)$ for some $a \in \del K_i^*$. Consider this
collection of discs modulo $\Z^N$ and select two, say with centres
$\pi^{\perp}(a)$ and $\pi^{\perp}(b)$, where $a,b \in \del K_i^*$.
Then, for $r$ very large and $a - b \in \pi'(\Z^N)$, these discs
will overlap modulo $\Z^N$. Since there are a finite number of
such pairs in $\del K_i^*$ to consider, we have a universal $r$
such that if $a,b \in \del K_i^*$ and $a-b \in \pi'(\Z^N)$, then
the disc with centre $\pi^{\perp}(a)$ overlaps, modulo $\Z^N$, the
disc with centre $\pi^{\perp}(b)$. If $a-b \nin\pi'(\Z^N)$, then
these discs will not overlap, modulo $\Z^N$. Hence for $r$
sufficiently large, $\del \Sigma^i \cap B(r) \mod \Z^N$ has
precisely $|A_i|$ components.

For the same $r$, $C_r^i \mod \Z^N$ is also a finite union of
discs of radius at least $r-1$ and at most $r$; likewise $C_r^i$
has exactly $|A_i|$ components, in direct correspondence with the
elements of $A_i$.

The description above of the limiting points in $\widetilde \Pi^i$
as we approach $C_r^i$ in a direction transverse to $E$, shows
that $M^i_r$ is homeomorphic to $\torus^N$ with a small open
neighbourhood of $C_r^i$ removed. There is a natural homotopy
equivalence with the space $\torus^N \setminus C_r^i$.

We can now examine what happens as we let first $r$ and then $i$
tend to infinity in this construction. For the above sufficiently
large $r$, the map $M^i_{r+1} \lra M^i_r$ is, up to homotopy, the
injection from $\torus^N \setminus C^i_{r+1}$ to $\torus^N
\setminus C_r^i$, and this is simply, up to homotopy, the identity
from $\torus^N \setminus |A_i|$ points to itself. Hence
$H^*(M_r^i) = H^*(\torus^N \setminus |A_i|$ points$)$ and
$H^*(M_r^i) \lra H^*(M_{r+1})$ is the identity showing that $H^*(M
\widetilde P^i)$ is the cohomology of the torus with $|A_i|$
punctures.

Finally, for each $i$ and for $r$ sufficiently large (depending on
$i$) the map $M_r^{i+1} \lra M_r^i$ is that induced by the
inclusion of $A_i$ in $A_{i+1}$, and this corresponds in the above
description to the adding of a new puncture for each element of
$A_{i+1}\setminus A_i$. In cohomology, the map $H^p(M_r^{i-1})
\lra H^p(M_r^{i})$ is thus the identity for $p \neq N-1$, and in
dimension $d$ gives rise to the direct system of groups and
injections
$\cdots\lra\Z^{N-1+|A_i|}\lra\Z^{N-1+|A_{i+1}|}\lra\cdots$ which
gives the required formula. \qed

\sn We give an alternative proof of this theorem from a different
perspective in Chapter V.

\sn We note that the pattern dynamical  system $(X, G_{1})$ is in
fact a Denjoy example \pss\ generalized to a $\Z^{N-1}$ action and
dislocation along $k$ separate orbits.

\sn{\bf Corollary 3.2} {\it Suppose that $\Gamma$ is a dense
countable subgroup of $\R$ finitely generated by $r$ free
generators. Suppose we Cantorize $\R$ by cutting and splitting
along $k$ $\Gamma$-orbits (as described e.g. in \pss ) to form the
locally Cantor space $R'$ on which $\Gamma$ acts continuously,
freely and minimally. Consider the $\Gamma$-module $C$ of
compactly supported integer valued functions defined on $R'$. Then
$H^*(\Gamma,C) = H^* (\torus^r \setminus k$ points$; \Z)$. \qed}

\newpage

\headline={\ifnum\pageno=1\hfil\else
{\ifodd\pageno\rightheadline\else\leftheadline\fi}\fi}
\def\rightheadline{\tenrm\hfil{\smallfont
IV INFINITELY GENERATED COHOMOLOGY}\hfil\folio}
\def\leftheadline{\tenrm\folio\hfil{\smallfont
FORREST HUNTON KELLENDONK}\hfil} \voffset=2\baselineskip

\def \flag{{\cal F}}

\sn{\chap IV Approaches to calculation II: }
\medskip
\noindent{\chap infinitely generated cohomology}
\bigskip

\sn{\sect 1 Introduction}

\sn In this chapter we restrict ourselves to the classical
projection tilings with canonical acceptance domains $K$ (so $K$
is the projection of the unit cube in $\R^N$). We examine the
natural question of when such tilings arise also as substitution
systems and show that the invariants of chapter 2 are effective
and computable discriminators of such tiling properties.

Much of this chapter is devoted to giving a qualitative
description of the cohomology of canonical projection method
patterns. The main result is formulated in Theorem~2.9. It gives a
purely geometric criterion for infinite generation of pattern
cohomology and for infinite rank of its rationalisation. As a
corollary of this, by the obstruction proved in II.5, we deduce
that almost all canonical projection method patterns fail to be
substitution systems; in fact for vast swathes of initial data
{\it all\/} such patterns fail to be self similar.

The restriction to the canonical acceptance domain allows for a
second, geometric description of the coefficient groups $CF_u$
(I.9.3) and $CV_u$ (II.4.3) in terms of indicator functions on
particular polytopes, we call them {\it $\C_u$-topes}. In \S2 we
introduce this viewpoint, setting up the notation and definitions
sufficient to state the main theorem. From here until the end of
\S5, our aim is to prove Theorem~2.9 by establishing criteria for
the existence of infinite families of linearly independent
$\C_u$-topes. In Sections 3 and 4 we construct such families in
the {\it indecomposable\/} case, and complete the analysis for the
general case in \S5. The final \S6 gives some general classes of
patterns where the conditions of Theorem~2.9 are satisfied, so
proving the generic failure of self-similarity for canonical
projection method patterns.

\sn{\sect 2 The canonical projection tiling}

\sn For the first time in our studies, we narrow our attention to
the classical projection method tilings of \okd\ \db. This section
outlines the simplifications to be found in this case, and
describes the main result of the remainder of this paper; a
sufficient condition for infinitely generated $H^d(\G\T)$. In
Chapter~V we will see that this condition is also neccessary.

Therefore we have data $(K,E,u)$ where $K = \pi^{\perp}([0,1]^N)$
and $u \in NS$. From Section I.8 we see that, but for a few
exceptional cases, we have $S_{\T}=\Z^N$, and if we elect either
to exclude these exceptions (as most authors do) or to include
them only in their most decorated form ($M\T = M\widetilde P_u$),
we make $S_{\T}=\Z^N$ a standing assumption.

With Theorem I.9.4 we have a description of the topology of
$CF_u$: it is generated by intersection and differences of shifts
of a single compact open set $\overline K$, formed by completing
$\pi'K \cap (Q+u) \cap NS$ (Lem.~I.9.6). Topologically $CF_u
\equiv CV_u \times \widetilde \Delta$ and $CV_u$ inherits the
subspace topology and a stabilizing subaction $\Gamma_{\T}$ of
$S_{\T}$ (see II.4.3).

With the choice of $K$ as a canonical acceptance domain above, we
may follow more closely the work of Le \le\ and give other more
geometrical descriptions of the elements of $CF_u$ and $CV_u$.

\sn{\bf Definition 2.1} For each $J \subset \{1,2,...,N\}$, we
construct a subspace $e^J = \langle e_j : j \in J \rangle$ (the
span) of $\R^N$, where $\{e_j\}$ is the standard unit basis of
$\R^N$ or $\Z^N$.

Write $\dim F = n$.

Let $\I = \{ J \subset \{1,2,...,N\} : \dim \pi'(e^J) = n -1 \}$
and define $\I^*$ to be the set of elements of $\I$ minimal with
respect to containment.

Define $\Z^N_{n-1} = \cup\{ e^J+v: |J|=n-1, v \in \Z^N \}$, i.e.
the $n-1$-dimensional skeleton of the regular cubic CW
decomposition of $\R^N$.

\sn The following Lemma gives some combinatorial information about
$\I^*$ and describes the singular points in $F$ - they are formed
by unions of affine subspaces of $F$ of codimension $1$.

\sn{\bf Lemma 2.2} {\it With the construction above,

 i/ $\I^*$ is a sub-collection of the
$n-1$-element subsets of $\{1,2...,N\}$. Also every subspace of
$F$ of the form $\pi'(e^J)$, with $|J| = n-1$, is contained in
$\pi'(e^{J'})$ for some $J' \in \I^*$.

ii/  $\R^N \setminus NS = \pi'^{-1}\pi'(\Z^N_{n-1})$ and
$\pi'(\Z^N_{n-1}) = F \setminus NS$.}

\sn{\bf Proof} i/ Straightforward.

ii/ With the data above, $\pi'(K)$ is a convex polytope in $F$,
with interior, each of whose extreme points is of the form
$\pi'(v)$ where $v \in \{0,1\}^N$. By i/, each of the faces of
$\pi'(K)$ is contained in some $\pi'(e^J+v)$ where $v \in
\{0,1\}^N$ and $J \in \I^*$. Also by i/, we have $\pi'(\Z^N_{n-1})
= \pi'(\cup\{e^J+v : v \in \Z^N, J \in \I^*\} )$.

However, by definition, $F \setminus NS$ is the union of the faces
of those polytopes of the form $\pi'(K+v)$, $v \in \Z^N$. Thus we
deduce immediately that $\pi'(\Z^N_{n-1}) = \pi'(\cup\{e^J+v : v
\in \Z^N, J \in \I^*\} ) \supset F \setminus NS$.

Conversely, a simple geometric argument shows that for each $J \in
\I^*$, there is a face of $\pi'(K)$ which is contained in
$\pi'(e^J+v)$ for some $v \in \Z^N$.

Assuming the claim, we have, for each $J \in \I^*$ and each $v \in
\Z^N$, some suitable shift of $\pi'(K)$, $\pi'(K+w)$ say, $w \in
\Z^N$, one of whose faces, $\Phi$ say, contains the point
$\pi'(v)$ as an extreme point, and $\Phi \subset \pi'(e^J+v)$.
With this same $\Phi$, therefore, we see that $\cup\{ \Phi +
\pi'(w): w \in \Z^N\} \supset \cup\{\pi'(e^J+w): w \in \Z^N\}$.
Thus $\pi'(\Z^N_{n-1}) = \cup\{ \pi'(e^J+v): v \in \Z^N, J \in
\I^*\} \supset F \setminus NS$. So we are done. \qed

\sn{\bf Definition 2.3} Recall the subspace $V$ from (I.2) and
write $\dim V = m$.

Now consider the set $\I^*(V) = \{J \in \I^* : \dim(\pi'(e^J) \cap
V) = m-1\}$.

\sn The following Lemma shows that the sets $\I^*(V)$ can be found
canonically from the sets $\I^*$ and describes the singular points
in $V+\pi'(u)$ preparatory to a description of $V_u$ (See II.4.3).
Once again, the singular points in $V+\pi'(u)$ form a union of
affine hyperplanes in $V+\pi'(u)$.

\sn{\bf Lemma 2.4} {\it i/ $\I^*(V) = \{J \in \I^* :  \pi'(e^J)
\cap V \neq V \}$.

ii/ If $u \in NS$, then $(V + \pi'(u)) \setminus NS = (V +
\pi'(u)) \cap \pi'(\Z^N_{n-1}) = (V + \pi'(u)) \cap (\cup \{
\pi'(e^J+v):v \in \Z^N, J \in \I^*(V)\})$.}

\sn{\bf Proof} Follows directly from definition of $V$ and Lemma
2.2 \qed

\sn With these notations in mind, we are well equipped to describe
the topology of $V_u$. Although $V_u$ is best described as a
subspace of $F_u$ formed from placing $V$ as the affine subspace
$V+\pi'(u)$ in $F$, we prefer to shift the whole construction back
to $V$ by applying a uniform shift by $-\pi'(u)$. There are some
advantages later in having an origin and a vector space structure.

\sn{\bf Definition 2.5} Given $u \in NS$, we use 2.4 i/ to define
a set, $\C_u$, of $m-1$-dimensional affine subspaces of $V$ whose
elements are of the form $(\pi'(e^{J} + v) \cap
(V+\pi'(u)))-\pi'(u)$ where $J \in \I^*(V)$ and $v \in \Z^N$. Such
a space may also be written, $\pi'(e^{J} + v-u) \cap V$.

We say that a subset of $V$ is a {\it $\C_u$-tope}, if

i/ it is compact and is the closure of its interior, and

ii/ it is a polytope, each of whose faces is a subset of some
element of $\C_u$.

We shall also say that a subset, $B$, of $V_u$ is a {\it
$\C_u$-tope} if $B$ is clopen and $\nu(B)-\pi'(u)$ is a
$\C_u$-tope subset of $V$ in the sense above (recall $\nu$ from
(I.9.2) above).

\sn We shall show that the $\C_u$-topes generate the topology on
$V_u$. To help this we describe a third possible topology as it is
found in \le.

For each element, $\alpha$, of $\C_u$, a hyperplane in $V$,
consider the half-spaces $H_{\alpha}^{\pm}$ defined by it. The
sets $H_{\alpha}^{\pm} \cap NS$ can be completed in the $\overline
D$ metric to a closed and open subset of $V_u$ which we write
$\overline H_{\alpha}^{\pm}$. We also call these subsets of $V_u$
{\it half-spaces}.

\sn{\bf Proposition 2.6} {\it The following three collections of
subsets of $V_u$ are the same:

i/ The collection of $\C_u$-topes.

ii/ The collection of compact open subsets of $V_u$.

iii/ The collection of those finite unions and intersections of
half-spaces in $V_u$ which are compact.}

\sn{\bf Proof} Every $\C_u$-tope in $V_u$ is clearly an element of
collection iii/ since every $\C_u$-tope in $V$ is a finite
intersection and union of half spaces in $V$.

Conversely, every open half-space in $V$ is a countable, locally
finite, union of the interiors of $\C_u$-topes (in $V$). Therefore
any pre-compact intersection of half spaces can be formed
equivalently from some union and intersection of a finite
collection of $\C_u$-topes, i.e.\ is a $\C_u$-tope, and each
$\C_u$-tope in $V_u$ is an intersection of clopen half-spaces.

Now compare collection ii/ with iii/ and i/.

Suppose that $\overline H$ is a half-space in $V_u$ defined by the
hyperplane $W \in \C_u$. As noted in 2.5, $W$ has the form
$\pi'(e^{J} + v-u) \cap V$ for some $J \in \I^*(V)$ and $v \in
\Z^N$, and our choice of half-space gives a corresponding choice
of half-space in $V$ on one side or other of $W$. Since $J \in
\I^*(V) \subset \I^*$ (2.4 i/ ), this in turn defines a choice of
half-space in $F$ one side or other of $\pi'(e^{J} + v-u)$. Write
$H$ for the closed half-space of $F$ chosen this way.

Consider the collection of $N$-dimensional cubes in the $\Z^N$
lattice whose image under the projection $\pi'$ is contained
entirely in $H + \pi'(u)$. By Lemma 2.4 ii/, the boundary of $H +
\pi'(u)$ is an $(n-1)$-dimensional hyperplane in $F$, the image
under $\pi'$ of an $(n-1)$-dimensional subspace of $\R^N$
contained in $\Z_{n-1}^N$. Therefore the cubes collected above
line up against this boundary exactly and, projected under $\pi'$,
they cover $H + \pi'(u)$ exactly. In particular, $H + \pi'(u)$ is
covered by $\Gamma_{\T}$-translates of $\pi'([0,1]^N)$.
Restricting this construction to the space $(V+\pi'(u)) \cap NS$,
we find that $\overline H$ is a locally finite countable union of
shifts of $\overline K$. In short, collection iii/ is contained in
collection ii/.

Conversely, we see that $\pi'([0,1]^N) \cap NS = \pi'((0,1)^N)
\cap NS$ so we are sure that the intersection $(V + \pi'(u)) \cap
((\pi'([0,1]^N)) + \pi'(v))$ (for any choice of $v \in \Z^N$) is a
polytope subset of $V+ \pi'(u)$ with interior. The faces of this
polytope are clearly subsets of $NS$, i.e. contained in elements
of $\C_u$ shifted along with $V$ by $\pi'(u)$. Thus every
$\Gamma_{\T}$ translate of $\overline K$ intersects $V_u$ either
non-emptily or as a $\C_u$-tope. Collection ii/ is contained in
collection i/ therefore. \qed

\sn Since $V_u$ is locally compact, we can describe the topology
on $V_u$, defined in II.4.3, equivalently as that generated by
$\C_u$-topes, or as that generated by half-spaces. This gives an
alternative description of $CV_u$.

\sn{\bf Corollary 2.7} {\it Let $V' = (V + \pi'(u)) \cap NS) -
\pi'(u)$. The group $CV_u$ is naturally isomorphic to the group of
integer-valued functions, $V' \lra \Z$, generated by indicator
functions of sets of the form $C \cap V'$ where $C$ is a
$\C_u$-tope. This group isomorphism is a ring and
$\Z[\Gamma_{\T}]$-module isomorphism as well. \qed}

\sn Now we describe an important set of points.

\sn{\bf Definition 2.8} Write $\P$ for the set of points in $V$
which can found as the $0$-dimensional intersection of $m$
elements of $\C_u$.  Note that, under the assumptions on $\C_u$,
$\P$ is a non-empty countable set, invariant under shifts by
$\Gamma_{\T}$.

Say that {\it $\P$ is finitely generated} if $\P$ is the disjoint
union of a finite number of $\Gamma_{\T}$ orbits, {\it infinitely
generated} otherwise.

\sn We may now express the main theorem of this Chapter:

\sn{\bf Theorem 2.9} {\it Given a canonical projection method
pattern, $\T$ and the constructions above. If $\P$ is infinitely
generated, then $H^d(\G\T)$ is infinitely generated and $H^d(\G\T)
\otimes \Q$ is infinite dimensional.}

\sn In Chapter V (Theorem V.2.4) we prove the converse i.e.\ that
$H^d(\G\T)$ infinitely generated is equivalent to $\P$ infinitely
generated.

We complete the proof of Theorem 2.9 in Section 5, but the final
step and basic idea, of the proof can be presented already.

\sn{\bf Proposition 2.10} {\it Suppose that $G$ is a torsion-free
abelian group and that $H$ is an abelian group, with $H/2:=H/2H$
its reduction $\mod 2$. Suppose that there is a group
homomorphism, $\phi : G \lra H/2$, such that $Im \phi$ is
infinitely generated as a subgroup of $H/2$ (equivalently infinite
dimensional as a $\Z/2$ vector space); then $G$ is infinitely
generated as an abelian group, and $G \otimes \Q$ is infinitely
generated as a $\Q$ vector space.}

\sn{\bf Proof} It is sufficient to prove the statement concerning
$G \otimes \Q$.

Suppose that $\phi(s_n)$ is a sequence of independent generators
for $Im \phi$ and suppose that there is some relation
$$\sum_{n=1}^m q_ns_n = 0$$ for $q_n\in\Q$. Since $G$ is
torsion-free, we can assume the $q_n$ are integers and have no
common factor; in particular, they are not all even. Applying the
map $\phi$ then gives a non-trivial relation among the
$\phi(s_n)$, a direct contradiction, as required. \qed

\sn Therefore, to prove 2.9, we shall find a homomorphism from
$H^d(\G\T)\cong H_0(\Gamma_\T,CV_u)$ to an infinite sum of $\Z/2$
whose image is infinitely generated. This is completed in full
generality in Theorem~5.4. In order to construct independent
elements of $H_0(\Gamma_\T,CV_u)$ and its image, we must consider
the geometry of the $\C_u$-topes defined in 2.5.

\bigskip

\sn{\sect 3 Constructing $\C$-topes}

\sn To prove that a group or $\Q$-vector space is infinitely
generated, we must produce independent generators. In the case of
$H_0(\Gamma_\T,CV_u)$, we must find elements of $CV_u$ which
remain independent modulo $\Gamma$-boundaries. In any case, we
must at least produce some elements of $CV_u$, and in this section
we start with constructions of the simplest objects in this space:
indicator functions of $\C_u$-topes.

Rather than refer constantly to the original tiling notation, we
abstract the construction conveniently, basing our development on
a general collection of affine hyperplanes, $\C$, of a vector
space, $V$, with group action, $\Gamma$. Always, the example in
mind is $\C = \C_u$ (2.5), $V$ (I.2) and $\Gamma = \Gamma_{\T}$
(II.4.3), but the construction is potentially more general.
However, the first few definitions and constructions are the
slightest generalization of those of section 2.

\sn{\bf Definition 3.1} Suppose that $V$ is a vector space of
dimension $m$ and that $\Gamma$ is a finitely generated free
abelian group acting minimally by translation on $V$. Thus we
write $w \mapsto w + \gamma$ for the group action by $\gamma \in
\Gamma$, and we may think of $\Gamma$ as a dense subgroup of $V$
without confusion.

Suppose that $\C$ is a countable collection of affine subspaces of
$V$ such that each $W \in \C$ has dimension $m-1$, and such that,
if $ W \in \C$ and $\gamma \in \Gamma$, then $W+\gamma \in \C$.

We suppose that the number of $\Gamma$ orbits in $\C$ is finite.

If $W \in \C$, then we define a unit normal vector, $\lambda(W)$
(with respect to some inner-product). The set $\N( \C) = \{
\lambda(W) : W \in \C\}$ is finite and we suppose that we have
chosen the $\lambda(W)$ consistently so that $- \lambda(W) \nin
\N( \C)$.

We suppose that $\N( \C)$ generates $V$ as a vector space.

\sn We can consider the intersections of elements of $\C$.

\sn{\bf Definition 3.2} Given $0 \leq k \leq m-1$, define
$\C^{(k)}$ to be the collection of $k$-dimensional affine
subspaces of $V$ formed by intersection of elements of $\C$. Thus
$\C$ can be written $\C^{(m-1)}$, and, to be consistent with the
notation of section 2, $\C^{(0)}$ can be written $\P$.

A {\it singular flag}, $\flag$, is a sequence of affine subspaces,
$(\theta_0, \theta_1,...,\theta_{m-1})$ of $V$ such that $\theta_j
\in \C^{(j)}$ for all $0 \leq j \leq m-1$, and $\theta_j \subset
\theta_{j+1}$ for all $0 \leq j \leq m-1$. The set of singular
flags is written $\J_o$.

It is clear that each $\C^{(k)}$ and $\J_o$ supports a canonical
$\Gamma$ action. We write $\J = \J_o/\Gamma$, the set of $\Gamma$
orbits in $\J_o$.

\sn{\bf Definition 3.3} We say that a subset, $C$, of $V$ is a
{\it $\C$-tope}, if

i/ $C$ is compact and is the closure of its interior, and

ii/ $C$ is a polytope, each of whose $(m-1)$-dimensional faces is
a subset of some element of $\C$.

\sn{\bf Definition 3.3} Convex $\C$-topes are finite objects whose
geometry and combinatorics are immediately and intuitively
related. Therefore we think of the definition of face, edge and
vertex, and more generally $k$-dimensional face, in this case, as
the intuitive one. Likewise the idea of incidence of edge on
vertex, face on edge, etc is intuitive.

Suppose that $C \subset V$ is a convex $\C$-tope and that $\flag =
(\theta_j)$ is a singular flag. We say that $\flag$ {\it is
incident on} $C$ if each $\theta_k$ contains a $k$-dimensional
face of $C$.

We write $[\flag]$ for the $\Gamma$-orbit class of $\flag$. We say
that $\flag$ {\it is uniquely incident on} $C$ if $\flag$ is
incident on $C$, but no other $\flag' \in [\flag]$ is incident on
$C$.

\sn The main aim of this section is to build convex $\C$ topes on
which certain singular flags are uniquely incident. This will not
happen in every possible circumstance however, and we develop an
idea of decomposability which will break up the space $V$ into a
direct sum of spaces on which such constructions can be made. In
section 5 we shall show how to recombine these pieces.

\sn{\bf Construction 3.5} Suppose that $A$ is a finite set of
non-zero vectors in $V$ which spans $V$, and no pair of which is
parallel. The example we have in mind is $\N(\C)$, the set of
normals.

A {\it decomposition} of $A$ is a partition $A = A_1 \cup A_2$
such that $V_1 \cap V_2 = 0$ where each $V_j$ is the space spanned
by $A_j$, $j =1,2$.

$A$ is {\it indecomposable} if no such decomposition is possible.
It is not hard to show that every set $A$ has a unique partition
into indecomposable subsets.

Suppose that $B \subset A$ is a basis for $V$. Then, by
stipulating that $B$ is an orthonormal basis, we define an inner
product which we write as square brackets: $[.,.]_B$.

Then we define a finite graph $G(B;A)$ with vertices $B$ and an
edge from $x$ to $y$ whenever there is a $z \in A \setminus B$
such that $[x,z]_B \neq 0$ and $[y,z]_B \neq 0$. We do not allow
loops.

\sn The following is elementary.

\sn{\bf Lemma 3.6} {\it Suppose that $A$ spans $V$, then the
following are equivalent:

i/ $A$ is indecomposable

ii/ for all bases $B \subset A$, $G(B;A)$ is connected

iii/ for some basis $B \subset A$, $G(B;A)$ is connected. \qed}

\sn{\bf Remark 3.7} Note that if $\phi : V \lra V$ is a linear
bijection, then $A$ is indecomposable if and only if $\phi(A)$ is
indecomposable. Therefore, the condition $\N(\C)$ indecomposable
can be stipulated without reference to a particular inner product,
although an inner product must be used to define the normals. We
use this freedom in Theorem 3.10 ahead.

\sn{\bf Construction 3.8} Suppose that $A$ and $B$ are as above,
giving an inner product $[.,.]_B$ to $V$ and defining a graph
$G(B;A)$. Choose $b \in B$ and let $W_b$ be the hyperplane in $V$
orthogonal to $b$ and let $\pi_b$ be the orthogonal projection of
$V$ onto $W_b$.

Consider the sets $\pi_b(A)$ and $\pi_b(B)$. Apart from $0=
\pi_b(b)$ the latter equals $B \setminus b$ precisely, an
othonormal basis for $W_b$. The former contains $B \setminus b$
and other vectors which may be of various lengths and pairs of
which may be parallel. From this set, we form $A_b$ a new set of
vectors in $W_b$ by taking, for each class of parallel elements of
$\pi_b(A)$ a single unit length representative (ignoring $0$). If
the class in question is one which contains some $b' \in B
\setminus b$, then we let $b'$ be the representative chosen. Write
$B_b = B\setminus b$, so that $B_b \subset A_b$ is a basis.

We consider the sets $A_b, B_b$ in $W_b$ and form the graph
$G(B_b;A_b)$ with respect to the inner product $[.,.]_{B_b}$, the
restriction of the inner product $[.,.]_B$ to $W_b$.

If $\C$ is a collection of affine hyperplanes in $V$, then we
define $\C_b$ to be collection of affine hyperplanes in $W_b$ of
the form $W_b \cap W$ where $W \in \C$ is chosen so that $W$ is
not parallel to $W_b$.

\sn The following useful lemma comes straight from the definition.

\sn{\bf Lemma 3.9} {\it i/ The graph $G(B_b;A_b)$ is formed from
$G(B;A)$ by removing the vertex $b$ and all its incident edges
from $G(B;A)$.

ii/ If $A = \N(\C)$, then $A_b = \N(\C_b)$.}

\sn{\bf Theorem 3.10} {\it Suppose that $m = \dim V > 1$ and that
$ 0 \nin A \subset V$ spans $V$, and $A$ has no parallel elements.
Suppose that $B \subset A$ is a basis for $V$ and that $G(B;A)$ is
connected, then there is a closed convex polytope, $C$, of $V$,
with interior, such that

i/ The normal of each face of $C$ is an element of $A$.

ii/  $v$ is a vertex of $C$ which is at the intersection of
exactly  $m$ faces of $C$ and each of these faces is normal to
some element of $B$.

iii/ The vertex $v$ is uniquely defined by property ii/.

(All normals are taken with respect to the inner product
$[.,.]_B$.)}

\sn{\bf Proof} We suppose, without loss of generality, that $v=0$
and proceed by induction on $m = |B| \geq 2$.

Suppose that $|B|=2$. By graph connectedness, we find $a \in A
\setminus B$ with non-zero components in each $B$ coordinate
direction. Thus we can construct easily a triangle, $C$, in $V$
with the required properties.

For larger values of $|B|$, we proceed as follows.

As in Construction 3.8, define, for each $b \in B$, $W_b$ to be
the hyperplane normal to $b$ and let $\pi_b : V \lra W_b$ be the
orthogonal projection.

A simple argument shows that we can find $b_o \in B$ so that
$G(B;A) \setminus b_o$ (i.e. removing the vertex $b_o$ and all
incident edges) is connected. Thus, appealing to the description
of Lemma 3.9 and by induction on $|B|$, we can find a convex
compact polytope subset, $C_1$, of $W_{b_o}$ which obeys
conditions i/ to iii/ above with respect to the basis $B \setminus
b$ and hyperplanes with normals parallel to some $\pi_{b_o}(a) : a
\in A \setminus b$.

Suppose that $W'$ is a hyperplane in $W_{b_o}$ which contains a
face of $C_1$ and has normal $\pi_{b_o}(a)$ for some $a \in A$.
Then $W'$ is in fact of the form $W \cap W_{b_o}$, where $W$ is
the unique hyperplane in $V$, normal to $a$, containing the space
$W'$. Indexing the faces of $C_1$ with $j$ say, each face is
contained in $W'_j$ with normal $\pi_{b_o}(a_j)$. So construct a
hyperplane $W_j$ in $V$ as above, with normal $a_j$. In the case
of the faces of $C_1$, incident on $0$ and normal to $b \in B
\setminus {b_o}$ say, we make sure that we choose $a_j = b$, so
that in this case $W_j = W_{b}$ as defined above.

Now we build $C_2$, a convex closed subset of $V$, as follows. For
each face $j$ of $C_1$, let $H_j$ be the closed half-space defined
by $W_j$, containing the set $C_1$. Let $C_2 = \cap_jH_j$, where
the intersection is indexed over all faces of $C_1$, so that $C_1
= C_2 \cap W_{b_o}$.

$C_2$ is almost certainly unbounded, but it is important to note
that none of its faces is normal to ${b_o}$. Furthermore, by
construction of $C_1$, there is only one edge in $C_2$ which sits
at the intersection of $m-1$ hyperplanes normal to some element of
$B \setminus b_o$, and this edge contains the point $0$ in its
interior and is parallel to $b_o$.

To produce a bounded set, we intersect $C_2$ with a compact convex
polytope $C_4$ whose faces are mostly normal to elements of $B$.
We construct $C_4$ as follows.

For each $b \in B \setminus {b_o}$, let $H_b$ be the half space of
$V$ defined by the hyperplane $W_b$ and containing $C_1$. Now let
$W^+_b$ an affine translation of $W_b$ beyond the other side of
$C_1$. Precisely, we argue as follows: since $W_b$ contains a face
of the convex set $C_1$, either $[b,c]_B \geq 0$ for all $c \in
C_1$, or $[b,c]_B \leq 0$ for all $c \in C_1$; we construct
$W^+_b$ in the first case leaving the second case to symmetry.
Since $C_1$ is compact,  there is an upper bound $r > [b,c]_B$ for
all $c \in C_1$. Let $W^+_b = rb+W_b$. Now let $H^+_b$ be the half
space defined by $W^+_b$ and containing $C_1$.

Let $H_{b_o}$ be a half space of $V$ defined by $W_{b_o}$; it
doesn't matter which one. The intersection, $$C_3 = \bigcap_{b \in
B} H_b \cap \bigcap_{b \in B \setminus b_o} H^+_b$$ is therefore a
semi-infinite rectangular prism whose semi-infinite axis runs
parallel to $b_o$ and whose base contains $C_1$.

We now form an oblique face to truncate this prism. Let $a_o \in A
\setminus b_o$ be chosen so that $[a_o,b_o]_B \neq 0$, as can
indeed be done by the assumed connectedness of $G(B; A)$. Let
$W_{a_o}$ be an affine hyperplane of $V$ with normal $a_o$; we
shall detail its placement soon. Note that, by constuction,
$W_{a_o}$ intersects every line parallel to $b_o$, in particular
we may place $W_{a_o}$ so that it intersects (the interior of) all
the edges of $C_3$ parallel to $b_o$. Let $H_{a_o}$ be the
half-space defined by $W_{a_o}$ and containing $C_1$. Let $C_4 =
H_{a_o} \cap C_3$.

The properties of $C_4$ are summarized: $C_4$ is compact and
convex: all its faces, except exactly one, are normal to some
element of $B$; exactly one of its faces is normal to $b_o$ and
that face contains $C_1$: $0$ is a vertex of $C_4$ and it is at
the intersection of exactly $m$ faces each normal to some element
of $B$: the other vertices of $C_4$ for which this can be said are
in $W_{b_o}$ but are all outside $C_1$.

Let $C = C_2 \cap C_4$. This is clearly a compact convex polytope
in $V$ and $C_1$ is a face of $C$. By construction, the point $0$
is a vertex of $C$ and obeys property ii/. Any other vertex of $C$
with property ii/ must be found in $C_1$ since there are no other
faces of $C_2$ or $C_4$ normal to $b_o$. However, no other vertex
of $C_1$ can have property ii/ by its original definition. \qed

\sn Although the theorem above makes no reference to group
actions, we can adapt it for use in section 4.

\sn{\bf Theorem 3.11} {\it Suppose $\C$ is a collection of
hyperplanes in $V$, $\dim V > 1$, with $\Gamma$ action as in
Def.~3.1, and suppose that $\N(\C)$ is indecomposable. Suppose
that $\flag$ is a singular flag (3.2). Then there is a convex
$\C$-tope on which $\flag$ is uniquely incident. }

\sn{\bf Proof} A singular flag is a descending sequence of
singular spaces and so we can also express it as the sequence of
spaces, $W_1, W_1 \cap W_2, \cap_{1 \leq i \leq 3}W_i,..., \cap_{1
\leq i \leq k} W_i , ..., \cap_{1 \leq i \leq m} W_i$, where $W_i
: 1 \leq i \leq m$ are transverse elements of $\C$. Let $\{v\} =
\cap_{1 \leq i \leq m} W_i$.

Fix some inner product $[.,.]$ in $V$ so that the $W_i$ are
orthogonally transverse and let $B = \{\lambda (W_i) : 1 \leq i
\leq m\}$ , where the normal $\lambda$ is taken with respect to
this inner product. $B$ is an orthonormal basis for $V$. We define
$\N(\C)$ with respect to this inner product (recall Remark 3.7
above) and, to fit it into past notation, let $A = \N(\C)$.

By hypothesis, $A$ is indecomposable. Thus the graph $G(B;A)$ is
connected, by 3.6, and so we may form by 3.10 a convex polytope,
$C_o$, in $V$ with the properties outlined in 3.10.

i/ The normal of each face of $C_o$ is an element of $A$.

ii/  $v$ is a vertex of $C_o$ which is at the intersection of
exactly  $m$ faces of $C_o$ and each of these faces is normal to
some element of $B$.

iii/ The vertex $v$ is uniquely defined by property ii/.

However, we know that the orbit of an element, $W$, of $\C$ is
dense in $V$ in the sense that for every affine hyperplane, $W'$,
of $V$, parallel to $W$, and every $\epsilon > 0$, there is an
$W''$, in the $\Gamma$ orbit class of $W$, such that $W'$ and
$W''$ are separated by a vector of length at most $\epsilon$.
Therefore we may adjust $C_o$ slightly without disturbing the
combinatorial properties of its faces to form a $\C$-tope, $C$,
with the same properties. The vertex $v$ need not be disturbed
atall.

However, then it is clear that $\flag$ is uniquely incident on
$\C$.\qed

\bigskip

\sn{\sect 4 The indecomposable case}

\sn We continue to consider the abstracted situation of section 3
and reintroduce the analysis of 2.6 and 2.7 as a definition.

\sn{\bf Definition 4.1} In the space, $V' = V \setminus \cup\{ W :
W \in \C\}$, define $\A_{\C}$ to be the collection of subsets of
$V'$ of the form $C \cap V'$, where $C$ is a $\C$-tope, and with
the empty set thrown in as well.

We write $CV_\C$ for the ring of integer-valued functions
generated by indicator functions of elements of $\A_{\C}$.

\sn Compare this definition with the construction of the topology
of $V_u$ in section 2.

Recall the set of singular flags, $\J_o$, for $V$ and $\C$ as
above, the $\Gamma$ action on $\J_o$, and the set $\J$ of $\Gamma$
orbits in $\J_o$. This transfers by a coordinatewise action to a
$\Gamma$ action on groups such as $\oplus_{\J_o}\Z/2$, the
$\J_o$-indexed direct sum of $\Z/2$. In this case $(\oplus_{\J_o}
\Z/2)/\Gamma = \oplus_{\J}\Z/2$ canonically.

\sn The main results of this section (4.2) and (4.6) are two
similar technical results. Here is the first.

\sn {\bf Proposition 4.2} {\it With the constructions and notation
of section 3 we suppose that $m = \dim V > 1$ and $\N(\C)$ is
indecomposable. Then, there is a $\Gamma$-equivariant homomorphism
$$\xi_o:CV_\C\to \oplus_{\J_{o}}\Z/2$$ with the following
property: for each singular flag $\flag$ there is an element $e
\in CV_\C$ so that $\xi_o(e)$ has value $1$ at coordinate $\flag$
and value $0$ at all other coordinates $\flag'\in [\flag]$. }

\sn The proof follows directly from section 3. Our aim is to build
a $\Gamma$-equivariant homomorphsism from $CV_{\C}$ to
$\oplus_{\J_o}\Z/2$ with certain further properties. We use the
following lemma to make the construction.

\sn{\bf Lemma 4.3} {\it i/ From 4.1 above, $\A_{\C}$ is an algebra
generated by sets of the form, $C \cap V'$, where $C$ is a convex
$\C$-tope.

ii/ There is a one-to-one correspondance between a group
homomorphisms $\xi : CV_{\C} \lra G$ to an abelian group $G$ and
maps $\xi'$ from convex $\C$-topes to $G$ with the property that
if $C_1$ and $C_2$ are interior disjoint convex $\C$-topes and if
$C_1 \cup C_2$ is convex, then $\xi'(C_1 \cup C_2) = \xi'(C_1) +
\xi'(C_2)$.

iii/ If $\Gamma$ acts homomorphically on $G$ then $\xi$ is
$\Gamma$-equivariant whenever the corresponding $\xi'$ is
$\Gamma$-equivariant.}

\sn{\bf Proof} i/ Every $\C$-tope can be decomposed into a finite
interior-disjoint union of convex $\C$-topes.

ii/ Using the decomposition of part i/ we can define a
homomorphism $\xi : \A_{\C} \lra G$ as follows: Take a $\C$-tope,
$C$, break it into interior-disjoint cover by convex $\C$-topes,
$C_j$. Let $\xi(C \cap V') = \sum \xi'(C_j)$.

To see that this is well-defined, consider a second such
decomposition $C = \cup_i C'_i$. Let $D_{i,j}$ be the closure in
$V$ of $V' \cap C'_i \cap C_j$; this is either empty or a convex
$\C$-tope. With the convention $\xi'(\emptyset) = 0$, the
condition on $\xi'$ extends inductively to show that $\xi'(C_j) =
\sum_i \xi'(D_{i,j})$ and $\xi'(C'_i) = \sum_j \xi'(D_{i,j})$.
Therefore $\sum_j \xi'(C_j) = \sum_i\xi'(C'_i)$ since both are
equal to $\sum_{i,j} \xi'(D_{i,j})$.

It is clear that $\xi$ is additive for interior-disjoint unions in
$\A_{\C}$, and so extends to $CV_{\C}$. It is straightforward to
construct a map $\xi'$ from $\xi$.

iii/ Follows quickly from the construction of ii/. \qed

\sn{\bf Proof of 4.2} It is easy to define a map $\xi'$ on convex
$\C$-topes which reflects the geometric idea of incidence of a
singular flag. Suppose $C$ is a convex $\C$-tope, then define
$\xi'(C)$ to be an element of $\oplus_{\J_o}\Z/2$ whose entry at
the $\flag$ coordinate ($\flag \in \J_o$ of course) is $1$ if and
only if $\flag$ is incident on $C$.

It is clear that $\xi'$ is $\Gamma$-equivariant.

To show that that $\xi'$ is additive, in the sense of Lemma 4.3
ii/, we consider two interior disjoint convex $\C$-topes $C_1$ and
$C_2$, whose union is a convex $\C$-tope, and fix a particular
singular flag, $\flag$. Three cases should be checked, the third
breaking down into two subcases.

Note first the general principle that a singular subspace (of any
dimension, $t$ say) containing a face (of dimension $t$) of $C_1
\cup C_2$ must therefore contain a face (of dimension $t$) of
$C_1$ or a face (of dimension $t$) of $C_2$. The converse, of
course, is not true as faces between $C_1$ and $C_2$ meet and fall
into the interior of a higher dimensional faces of $C_1 \cup C_2$;
and this is generally the only way that faces can be removed from
consideration. So we have a second general principle that a
singular subspace contains a face of $C_1 \cup C_2$ when it
contains a face of $C_1$ and no face of $C_2$ (matching dimensions
always).

Case a/: $\flag$ is incident neither on $C_1$ nor on $C_2$. By our
first general principle, it is immediate that $\flag$ is not
incident on $C_1 \cup C_2$.

Case b/: $\flag$ is incident on precisely one of $C_1$ or $C_2$.
The second general principle above shows then that $\flag$ is
incident on $C_1 \cup C_2$.

Case c/: $\flag = (\theta_j)_{0 \leq j < N}$ is incident on both
$C_1$ and $C_2$. Suppose that $\theta_t$ is the singular subspace
(of dimension $t$) containing a face $F_1$ (of dimension $t$) of
$C_1$ and a face $F_2$ (of dimension $t$) of $C_2$. We analyse two
possibilities: i/ $F_1 = F_2$: ii/ $F_1$ and $F_2$ are interior
disjoint (as subsets of $\theta_t$). Note that we have used the
convexity of $C_1 \cup C_2$ and other assumed properties to make
this dichotomy.

Case ci/: Here $F_1$ is no longer a face of $C_1\cup C_2$ and
$\theta_t$ does not contain a face of $C_1\cup C_2$. Thus $\flag$
is not incident on $C_1\cup C_2$.

Case cii/: Here $\theta_t$ contains the face $F_1 \cup F_2$ of
$C_1\cup C_2$. However, consider $\theta_{t-1}$ (note that $t \geq
1$ automatically in case ii/) which contains a face of both $C_1$
and $C_2$ of dimension $t-1$. In fact, by convexity, this face
must be $F_1 \cap F_2$. However, this set is not a $t-1$
dimensional face of $C_1 \cup C_2$, having been absorbed into the
interior of the $t$-dimensional face $F_1 \cup F_2$.

Either way, in case c/ $\flag$ is not incident on $C_1 \cup C_2$.

Combining all these three cases gives the additivity mod $2$
required of $\xi'$. Therefore, we define a $\Gamma$-equivariant
homomorphism $\xi_o : CV_{\C} \lra \oplus_{\J_o}\Z/2$ as required.

The element $e$ is provided by the construction of 3.11. Consider
a convex $\C$-tope, $C$ produced by Theorem 3.11 from the singular
flag $\flag$. The indicator function of $C$ is an element of
$CV_{\C}$ which we will write $e$. The properties claimed of $e$
in the theorem follow automatically from unique incidence of
$\flag$ on $C$. \qed

\sn $\Gamma$-equivariance of the homomorphism from Proposition~4.2
allows us to build the commuting diagram $$\matrix{CV_\C &\bra{
\xi_{o}}& \oplus_{\J_{o}}\Z/2 \cr&&\cr\downarrow^{q}&&
\downarrow^{q} \cr&&\cr H_0(\Gamma;CV_\C) &\bra{\xi_{*}}&
\oplus_{\J_{}}\Z/2 \cr}$$ to define $\xi_*$, where the maps $q$
quotient by the action of $\Gamma$. In the indecomposable case
therefore, Theorem~5.4 below is a direct consequence of this
observation and the reader who is only interested in the
indecomposable case can jump directly to that theorem replacing
the first words "Given the data above" by the first sentence of
Proposition 4.2. In the decomposable case, however, we need a lot
more work to establish such a diagram.

\sn We prove an analogous technical result for the case $\dim V =
1$. It seems necessary and slightly surprising that its proof is
not as directly geometric.

For the remainder of the section, therefore, we assume that $\dim
V = 1$, but, for technical reasons, we also relax the assumption
that $\Gamma$ is a subgroup of $V$. Rather $\Gamma$ is a free
abelian group acting minimally by translation, but some of these
translations may be $0$. In this case we decompose $\Gamma =
\Gamma_0 \oplus \Gamma_1$ where $\Gamma_1$ acts minimally and
freely and $\Gamma_0$ fixes every point in $V$.

We note that in the case $\dim V = 1$, the flags are simply single
points from $\P$ (sections 2 and 3). Thus there is a canonical
correspondence between $\J_o$ and $\P$, which is clearly
$\Gamma$-equivariant.

\sn{\bf Construction 4.4} Consider the direct sum group $L =
\oplus_{\J}\Z$ whose elements can be considered as collections of
integers indexed by elements of $\J$, or equivalently by
$\Gamma$-orbit classes in $\P$. Let $h_j : j=1,2,..,k$ (with the
obvious adaption for $k = \infty$) be the canonical free
generating set for $L$.

Let $\P_1,...\P_k$ be the $\Gamma$-orbit classes of $\P$ and
choose $x_0\in \P_1$ and, for each $1 \leq j \leq k$, also choose
$x_j \in \P_j$ such that $x_j > x_0$ (recall that $\P \subset V =
\R$ in the case $\dim V =1$).

According to definition 4.1 above, $V' = V \setminus \P$, and the
sets $I_j = [x_0,x_j] \cap V'$ are elements of $\A_{\C}$. Thus the
functions $f_j$ which indicate respective $I_j$, are elements of
$CV_{\C}$.

This allows us to define a homomorphism $\beta : L \lra CV_{\C}$,
defined $\beta(h_j) = f_j$ for each $1 \leq j \leq k$. Consider
$\oplus_{\Gamma_1}L$, the $\Gamma_1$ indexed direct sum of copies
of $L$, whose elements we shall consider as $\Gamma_1$-indexed
elements of $L$, $(g_{\gamma})_{\gamma \in \Gamma_1}$, all but a
finite number of which are $0$. Note that, by correponding the
coordinate indices $\Gamma$-equivariantly, $\oplus_{\Gamma_1}L$ is
$\Gamma$-equivariantly isomorphic to $\oplus_{\J_o}\Z$. (The
$\Gamma_0$ component of the action acts trivially).

Now consider the homomorphism $\beta^* : \oplus_{\Gamma_1} L \lra
CV_{\C}$, defined as $$\beta^*((g_{\gamma})_{\gamma \in \Gamma_1})
= \sum_{\gamma \in \Gamma_1} \gamma \beta (g_{\gamma})$$

\sn{\bf Lemma 4.5} {\it With the construction above, $\beta^*$ is
$\Gamma$-equivariant.

i/ $\beta^*: \oplus_{\Gamma_1} L \lra CV_{\C}$ is injective.

ii/ The image of $\beta^*$ is complemented in $CV_{\C}$.}

\sn{\bf Proof} The proof relies on a construction whose proper
generalization is made in Chapter 5. In the case $\dim V = 1$,
however, it is easy enough to describe directly.

Consider those subsets, $I_{a,b}$, of $\A_{\C}$ formed from a
$\C$-tope interval $[a,b] \cap V'$ (much as the sets $I_j$ were
formed above). In this case $a,b \in \P$ necessarily. Therefore to
this set, we associate an element of $\oplus_{\P}\Z$ namely one
which is $0$ at every coordinate except the $a$ coordinate, where
it is $1$, and at the $b$ coordinate, where it is $-1$. We write
this element $\delta(I_{a,b}) = 1_a-1_b$ in an obvious notation.
This function is additive in the sense of Lemma 4.3 and so extends
to a homomorphism $\delta: CV_{\C} \lra \oplus_{\P}\Z$. Moreover,
$\delta$ is clearly $\Gamma$-equivariant on such intervals, and
hence it is equivariant when extended to a group homomorphism.

Recall the elements, $f_j$, of $CV_{\C}$ defined before. Consider
an equation of the form $\sum_{1 \leq j \leq k, \gamma \in
\Gamma_1} t_{j,\gamma}\gamma f_j = 0$ (where $t_{j,\gamma} \in \Z$
equal $0$ for all but a finite set of indices). If we apply
$\delta$ to this, we find $\sum_{1 \leq j \leq k, \gamma \in
\Gamma_1} t_{j,\gamma}\gamma \delta(f_j) = 0$, a sum in
$\oplus_{\P}\Z$, and we can now use the equation $\delta(f_j) =
1_{x_0}-1_{x_j}$.

Decompose $\oplus_{\P}\Z = \oplus_j \oplus_{\P_j}\Z$ and examine
what happens on each $\oplus_{\P_j}\Z$.

Consider first an index $j \geq 2$. In the component
$\oplus_{\P_j}\Z$ for such a $j$, the sum reduces to  $\sum_{
\gamma \in \Gamma_1} t_{j,\gamma}\gamma 1_{x_j} = 0$ whence
$t_{j,\gamma}=0$ for all $\gamma \in \Gamma_1$ since $\Gamma_1$
acts freely on $\P_j$.

The case remaining is $j=1$ for which the sum reduces to $$\Bigl(
\sum_{1 \leq j \leq k, \gamma \in \Gamma_1} t_{j,\gamma}\gamma
1_{x_0} \Bigr) - \Bigl( \sum_{ \gamma \in \Gamma_1}
t_{1,\gamma}\gamma 1_{x_1} \Bigr) = 0$$ however, by the last
paragraph, this reduces to $\sum_{ \gamma \in \Gamma_1}
t_{1,\gamma}\gamma (1_{x_0}-1_{x_1}) = 0$. But $x_1 =
\gamma_1(x_0)$ for some $\gamma_1 \in \Gamma_1$ and so we find
$\gamma_1f = f$ where $f = \sum_{ \gamma \in \Gamma_1}
t_{1,\gamma}\gamma 1_{x_0}$. This implies that $f=0$ ($\Gamma_1$
acts freely on $\oplus_{\P}\Z$) and so we have $t_{1,\gamma}=0$
for all $\gamma \in \Gamma_1$ as well.

In conclusion, the equation $\sum_{1 \leq j \leq k, \gamma \in
\Gamma_1} t_{j,\gamma}\gamma f_j = 0$ implies $t_{j,\gamma}=0$ for
all $j$ and $\gamma$, and so i/ follows.

To show ii/ we note first that $CV_{\C}$ is isomorphic to a
countable direct sum of $C(X;\Z) \equiv \oplus_{\infty}\Z$ ($X$
Cantor) and so is itself free abelian. Therefore to prove
complimentarity of the image of $\beta^*$ it's enough to show that
$CV_{\C}/Im\beta^*$ is torsion free.

Therefore, we argue to the contrary and suppose that we have an
element $g = (g_{\gamma})_{\gamma \in \Gamma_1} \in
\oplus_{\Gamma_1}L$ for which $\beta^* (g) = tf$ for some $t \in
\Z$, $t\geq 2$, and $f \in CV_{\C}$. This equation may be written
$\sum_{1 \leq j \leq k, \gamma \in \Gamma_1} t_{j,\gamma}\gamma
f_j = tf$, where $\beta(g_{\gamma}) = \sum_{1 \leq j \leq k}
t_{j,\gamma} f_j$.

The analysis we have just completed in part i/ can be performed
equally well modulo $t$. Therefore each $t_{j,\gamma}= 0 \mod t$
and we deduce that each $\beta(g_{\gamma}) = tf_{\gamma}$ for some
$f_{\gamma} \in \beta(L)$. In particular $f = \sum_{\gamma \in
\Gamma_1} \gamma f_{\gamma} \in \beta^*(\oplus_{\Gamma_1}L)$. Thus
every element of $\beta^*(\oplus_{\Gamma_1}L)$ which can be
divided in $CV_{\C}$ can also be divided in
$\beta^*(\oplus_{\Gamma_1}L)$, i.e. the quotient is torsion-free
as required. \qed

\sn   This gives the key to the second and final technical result
of this section.

\sn {\bf Proposition 4.6} {\it Assume the constructions and
notation of 4.4 and section 3: in particular we suppose that $\dim
V = 1$. Then there is a surjective $\Gamma$-equivariant
homomorphism $\xi_o:CV_\C \to \oplus_{\J_{o}}\Z/2$. }

\sn{\bf Proof} Lemma 4.5 builds an injective $\Gamma$-equivariant
homomorphism $\beta^*: \oplus_{\J_o}\Z \lra CV_{\C}$ whose image
is complemented in $CV_{\C}$. Therefore we have automatically a
surjective $\Gamma$-equivariant homomorphism $\xi: CV_{\C} \lra
\oplus_{\J_o}\Z$ reversing this (i.e. $\xi \beta^* = $ identity).
Extending $\xi$ by the reduction mod $2$, $\oplus_{\J_o}\Z\lra
\oplus_{\J_o}\Z/2$, gives us the homomorphism $\xi_o: CV_{\C} \lra
\oplus_{\J_o}\Z/2$. \qed

 \bigskip

\sn{\sect 5 The decomposable case}

\sn Proposition 4.2 is very close to a proof of Theorem 2.9 in the
case that $\N(\C)$ is indecomposable. Rather than finish the
argument for that case, we clear up the decomposable case and so
complete the most general proof. Once again we exploit the
generalities of section~3.

\sn{\bf Construction 5.1} Suppose that $\N(\C)$ is decomposable
and that $\N_1, \N_2,..., \N_k$ are its components. Let $V_j$ be
the vector space spanned by $\N_j$, defined for each $1 \leq j
\leq k$. Without loss of generality (see Remark 3.7) we can impose
an inner product on $V$ which will make the $V_j$ mutually
orthogonal and we do this from the start.

Therefore $V = \sum_j V_j = \oplus_jV_j$ is an orthogonal direct
sum. Let $\pi_j$ be the orthogonal projection onto $V_j$ with
kernel $\sum\{V_i : i \neq j\}$.

Let $\C_j$ be the set of those elements of $\C$ whose normal is
contained in $V_j$. Thus, for each $j$, $\N_j = \N(\C_j)$, an
indecomposable subset of $V_j$. The $\Gamma$ action on $\C$ leaves
each of the sets $\C_j$ invariant.

The natural $\Gamma$ action on $V_j$ is more complicated than mere
restriction. We consider $\pi_j(\Gamma)$ as a subgroup of $V_j$
and let $\Gamma$ act by translation: $\gamma(v) = v +
\pi_j(\gamma)$. We call this the {\it projected action}.

On each $V_j$ we construct the sets of singular flags, $\J_{oj}$
and $\J_j=\J_{oj}/\Gamma$, and the spaces $CV_{j\, \C_j^*}$
according to section 3, using the projected action of $\Gamma$ on
$V_j$ and the singular hyperplanes $\C^*_j = \{W \cap V_j: W \in
\C_j\}$. Note again that, for each $j$, $\N(\C^*_j) = \N(\C_j)$ is
indecomposable as a spanning subset of $V_j$. Moreover, $\P_j$,
the intersection point set defined using $\C^*_j$ in $V_j$, is
equal to $\pi_j(\P)$ which in turn is equal to $\P \cap V_j$.

We abbreviate the $\Z[\Gamma]$ module $CV_{j\, \C^*_j}$ as $CV_j$
in what follows.

\sn{\bf Lemma 5.2} {\it Given the constructions above, we have the
following relations:

i/ There is a canonical bijection $\C \leftrightarrow \cup_j
\C^*_j$.

ii/ There is a canonical $\Gamma$-equivariant bijection $\Pi_j
\P_j \leftrightarrow \P$ (direct product of sets), where $\Gamma$
acts on the direct product diagonally by its projected actions.

iii/ $V \equiv \oplus_j V_j$ is a canonical orthogonal direct sum
on which the usual $\Gamma$ action on $V$ is retrieved as the
diagonal action of the projected actions of $\Gamma$.

iv/ $CV \equiv \otimes_j CV_j$ canonically as a
$\Z[\Gamma]$-module, where the module action is given by the
diagonal action of the projected actions.

v/ There is a canonical $\Gamma$-equivariant surjection $\sigma:
\J_o \lra \Pi_j \J_{oj}$ (direct product of sets) and a natural
$\Gamma$-equivariant injection $\tau: \Pi_j \J_{oj} \lra \J_o$
whose composition $ \sigma \tau $ is the identity.}

\sn{\bf Proof} The first four parts all follow from the fact that,
for each $j$, the singular spaces in $V$ break into two classes.
Those spaces which are parallel to $V_j$ do not impinge on its
singular geometry at all. And those singular spaces, $W$, which
are not parallel, intersect $V_j$ in exactly the same way as they
project to $V_j$, i.e. $\pi_j(W) = W \cap V_j$.

Part v/ is more involved. Consider a singular flag, $\flag$, in
$V$ with the singular plane set $\C$. The bijection of part i/ can
be detailed: each element of $\C$ is of the form $V_1 + V_2 + ...+
V_{j-1} + W_j + V_{j+1} + ..+V_k$, an orthogonal sum, where $W_j
\in \C^*_j$. Thus any intersection of elements of $\C$ can be
written as an orthogonal sum $W_1+W_2+...+W_k$, where, for each
$j$, $W_j$ is an intersection of elements from $\C^*_j$.

Now consider the sequence of singular spaces listed in $\flag$ and
their orthogonal decomposition as above. As we read from point to
hyperplane, their dimension rises by exactly one and so the
orthogonal summands rise in dimension, but only in one of the
summands and only by one dimension. Suppose in the orthogonal
decomposition of this increasing sequence, we ignore all
directions but the $j$th say. Then we find a nested sequence of
singular spaces in $V_j$ whose dimension rises by at most $1$ at
each step, either going all the way up to $V_j$ or (for precisely
one value of $j$) stopping one dimension short. Extract the
sequence as a strictly increasing subsequences, neglecting the
last term if it happens to be $V_j$, and call the result
$\flag_j$. For each $j$, $\flag_j$ will be a singular flag in
$V_j$ with respect to the singular planes $\C^*_j$. The map $\flag
\mapsto (\flag_1,\flag_2,...,\flag_k)$ is $\sigma$.

To show it's onto, we construct $\tau$. Given $\flag_j =
(W_{j,0},...,W_{j,m_j-1})$ (where dimension of $V_j$ is $m_j$)
consider the sequence: $W_{1,0} + W_{2,0} + ... + W_{k,0}$,
$W_{1,1} + W_{2,0} + ... + W_{k,0}$, $W_{1,2} + W_{2,0} + ... +
W_{k,0}$,..., $W_{1,m_1-1} + W_{2,0} + ... + W_{k,0}$, $V_1 +
W_{2,0} + ... + W_{k,0}$, $V_1 + W_{2,1} + ... + W_{k,0}$, ...,
$V_1 + W_{2,m_2-1} + ... + W_{k,0}$, $V_1 + V_2 + W_{3,0}+ ... +
W_{k,0}$,..., $V_1+V_2+...V_{k-1} + W_{k, m_k-2}$,
$V_1+V_2+...V_{k-1} + W_{k, m_k-1}$, where all sums are orthogonal
sums in $V$. This is a singular flag, in $V$ with respect to $\C$,
which we shall call $\tau(\flag_1,\flag_2,...,\flag_k)$.

It is immediate from definition that $\sigma$ and $\tau$ are
$\Gamma$-equivariant and that $\sigma \tau$ is the identity. \qed

\sn Note that although $\P$ may have an infinite number of
$\Gamma$-orbits, each of the $\P_j$ may have only finitely many
$\Gamma$ orbits under the projected action. This subtle problem
cuts out the possibility of an easy argument by induction on the
number of indecomposable components of $\N(\C)$.

Applying $\tau$ coordinate-wise, we find

\sn{\bf Corollary 5.3} {\it There is a $\Gamma$-equivariant
homomorphism $\tau^+ : \otimes_j \oplus_{\J_{oj}} \Z/2 \lra
\oplus_{\J_o} \Z/2$. \qed}

\sn Now we are ready for the main result of this section.

\sn{\bf Theorem 5.4} {\it Given the data above, there is a
homomorphism $\xi_*: H_0(\Gamma; CV_\C) \lra  \oplus_{\J}\Z/2$
such that for all $v \in \P$, there is a singular flag $\flag =
(\theta_j)_{1 \leq j < m} \in \J_o$ so that $\theta_0 = \{v\}$,
and an element $e_v \in H_0(\Gamma, CV_\C)$ such that $\xi_*(e_v)$
has value $1$ at coordinate $[\flag]$ (i.e. the $\Gamma$-orbit
class of $\flag$, an element of $\J$ (3.3)).}

\sn{\bf Proof} Suppose that $\N(\C) = \cup_j \N(\C_j)$ is an
indecomposable partition, forming the spaces $V_j$, of dimension
$m_j = \dim V_j$, etc as above.

Note that on each $V_j$, $\Gamma$ acts by translation by elements
of $\pi_j(\Gamma)$ and that action may or may not be free.
Proposition 4.6 applies in this case by hypothesis. Also the proof
of Proposition 4.2 does not depend on the freedom of the $\Gamma$
action. No complication arises therefore if we stick with the
projected $\Gamma$ action on each $V_j$ even though the action may
not be free.

For each $j$, $\pi_j(\Gamma)$ is the projection of a dense subset
of $V$ and so itself is dense in $V_j$. Therefore the $\Gamma$
action on each $V_j$ is minimal and $rank (\pi_j(\Gamma)) > 1$.

Given $v \in \P$, we find $\pi_j(v) = v_j \in \P_j$. By
Propositions 4.2 and 4.6 we find for each $j$ a homomorphism $$
\xi_{oj} : CV_j \lra \oplus_{\J_{oj}}\Z/2,$$
a singular flag $\flag_j$ in $V_j$ such that $\flag_{j0} =
\{v_j\}$, and an element $e_{v,j} \in CV_j$ so that
$\xi_o(e_{v,j})$ has value $1$ in the coordinate $\flag_j$ and $0$
in the coordinates of all other flags from the $\Gamma$-orbit
$[\flag_j]$ of $\flag$ (if $\dim V_j=1$ flags are just points and
the latter properties follow from surjectivity). The homomorphisms
$\xi_{oj}$ are independent of the choice of $v$.

The equation $CV_\C = \otimes_j CV_j$ allows us to build the
homomorphism $$\otimes \xi_{oj} : CV_\C \lra \otimes_j
\oplus_{\J_{oj}}\Z/2$$ and using $\tau^+$ of Cor 5.3, we can
continue this homomorphism to $ \oplus_{\J_{o}} \Z/2$. This
homomorphism is clearly $\Gamma$-equivariant and so we deduce a
quotiented homomorphism: $\xi_* : H_0(\Gamma; CV_\C) \lra
\oplus_{\J} \Z/2$ which completes a square as required
$$\matrix{CV_{\C} &\bra{ \tau^+(\otimes_j \xi_{oj})}&
\oplus_{\J_{o}}\Z/2 \cr&&\cr\downarrow^{q}&& \downarrow^{q}
\cr&&\cr H_0(\Gamma;CV_\C) &\bra{\xi_{*}}& \oplus_{\J}\Z/2 \cr}.$$
Given $v$, the element $e_v = \otimes_j e_{v,j} \in CV_\C$ is
mapped by $\tau^+(\otimes_j \xi_{oj})$ to an element with value
$1$ at the coordinate $\flag = \tau(\flag_1,\flag_2,...,\flag_k)$
and with value $0$ at coordinates $\flag'=
\tau(\flag_1',\flag_2',...,\flag_k')$ where $\flag_j'\in
[\flag_j]$ but $\flag\neq\flag'$. In particular, $\tau^+(\otimes_j
\xi_{oj})(e_v)$ has value $0$ at all coordinates of the
$\Gamma$-orbit $[\flag]$ of $\flag$ except at $\flag$ itself.
Hence $q \tau^+(\otimes_j \xi_{oj}) : CV_\C \lra \oplus_{\J} \Z/2
$ maps $e_v$ to an element with value $1$ at the coordinate
$[\flag]$. The zero dimensional element of $\flag$ is $\{v\}$ by
construction. \qed

\sn{\bf Proof of Theorem 2.9}  Suppose that $\P$ is infinitely
generated and that $v_1,v_2,...$ are representatives of distinct
$\Gamma$-orbit classes. By Theorem 5.4, we find a homomorphism
$\xi_*: H_0(\Gamma; CV_\C) \lra  \oplus_{\J}\Z/2$, singular flags
$\flag_j \in \J_o$ so that the zero dimensional element of
$\flag_{j}$ is $\{v_j\}$, and elements $e_j
\in
H_0(\Gamma, CV_\C)$ such that $\xi_*(e_j)$ has value $1$ at
coordinate $[\flag_j]$.

By taking a subsequence if necessary we may assume that for $i <
j$, all the values $\xi_*(e_i)$ have value $0$ at the $[\flag_j]$
coordinate. In particular, we have ensured that the set $\{
\xi_*(e_j) : j \geq 1 \}$ is $\Z/2$ independent in
$\oplus_{\J}\Z/2$. By Proposition 2.10 we have $H_0(\Gamma, CV_\C)
\otimes \Q$ infinitely generated.

Now to prove Theorem 2.9, we apply this analysis to the
construction of section 2, using the same $V$, and setting $\C =
\C_u$ and $\Gamma = \Gamma_{\T}$. \qed

\bigskip

\sn{\sect 6  Conditions for infinitely generated cohomology}

\sn To apply Theorem 2.9 we must be able to count the orbits in
$\P$. This is a geometric exercise, and each case will have its
own peculiarities. We present in this section elementary general
conditions which are sufficient to give infinite orbits in $\P$.

Recall the general set-up from IV.3, and the construction of $\P$
as points which are the proper intersection of $m=\dim V$
hyperplanes picked from $\C$.

\sn{\bf Definition 6.1} Suppose that $W_1,..,W_m$ is a set of
hyperplanes chosen from $\C$, intersecting in a single point $p$.
For each subset $A$ of $\{1,\cdots,m\}$ the intersection
$$W_A:=\bigcap_{i\in A} W_i$$ has dimension $m-|A|$. We write
$A^c=\{1,\cdots,m\}\setminus A$ and define $\Gamma_A = \{x \in V |
\exists \gamma \in \Gamma: \{x + p\} = W_A \cap (W_{A^c} +
\gamma)\}$. Finally let $\Gamma^A\subset\Gamma$ be the stabilizer
of $W_A$.

\sn We think of $\Gamma_A$ as the projection of $\Gamma$ onto
$W_A$ along $W_{A^c}$. The following is straightforward from the
definitions.

\sn{\bf Lemma 6.2} {\it With the notation above:

i/ $\Gamma_A$ is a group with $\Gamma^A$ as a subgroup.

ii/ $\Gamma_A + p = \P \cap W_A$.

iii/ If $q \in \P \cap W_A$, then $(\Gamma + q) \cap W_A =
\Gamma^A + q$.

iv/ If $A_1,A_2\subset \{1,\cdots,m\}$ are disjoint, then
 $\Gamma^{A_1\cup A_2}= \Gamma^{A_1}\cap \Gamma^{A_2}$.}

\sn This gives immediately an easy way to determine whether we
have an infinite number of orbits in $\P$.

\sn{\bf Proposition 6.3} {\it If, for some choice of $W_1,...,W_m$
and $A\subset\{1,\cdots,m\}$, the stabilizer $\Gamma^A$ has
infinite index in $\Gamma_A$ (equivalently, if $\rk \Gamma^A<\rk
\Gamma_A$), then $\P$ is infinitely generated.}

\sn{\bf Proof} By Lemma 6.2 ii/ and iii/ the orbits in $\P$ which
intersect $W_A$ are enumerated precisely by the cosets of
$\Gamma^A$ in $\Gamma_A$. This is infinite by assumption. \qed

\sn We can pursue the construction above a little further to get
even a sharper condition for infinitely generated $\P$. Given
$A_1,A_2\subset\{1,\cdots,m\}$, note that $A_1 \cup A_2 =
\{1,2,..,m\}$ implies $W_{A_1}\cap W_{A_2} = \{p\}$ and therefore
$\Gamma_{A_1} \cap \Gamma_{A_2} = \{0\}$. This together with
Lemma~6.2 gives the following result.

\sn{\bf Lemma 6.4} {\it For every choice of $W_1,...,W_m$ as above
and for every pair of sets, $A_1,A_2\subset\{1,\cdots,m\}$ we
have,

i/ if $A_1\cap A_2=\emptyset$, then $$\rk \Gamma^{A_1} +
\rk\Gamma^{A_2} -\rk \Gamma^{A_1\cup A_2} \leq \rk\Gamma \leq \rk
\Gamma_{A_1} + \rk\Gamma_{A_2} -\rk (\Gamma_{A_1} \cap
\Gamma_{A_2}),$$

ii/ if $A_1 \cup A_2 = \{1,2,..,m\}$, then $$\Gamma^{A_1} +
\Gamma^{A_2} \subset \Gamma^{A_1\cap A_2} \subset \Gamma_{A_1\cap
A_2} \subset \Gamma_{A_1} + \Gamma_{A_2}$$ both sums being
direct.} \qed

\sn{\bf Corollary 6.5} {\it If $\P$ is finitely generated, then
for every choice of $W_1,...,W_m$ and for every pair of sets,
$A_1,A_2\subset\{1,\cdots,m\}$ we have,

i/ if $A_1\cap A_2=\emptyset$, then $$\rk \Gamma^{A_1} +
\rk\Gamma^{A_2} -\rk \Gamma^{A_1\cup A_2} = \rk\Gamma,$$

ii/ if $A_1 \cup A_2 = \{1,2,..,m\}$, then $$\rk\Gamma^{A_1} +\rk
\Gamma^{A_2} =\rk \Gamma^{A_1\cap A_2}.$$}

\bew\ Suppose $A_1\cap A_2=\emptyset$. Then $\Gamma^{A_1\cup A_2}=
\Gamma^{A_1} \cap \Gamma^{A_2} \subset \Gamma_{A_1} \cap
\Gamma_{A_2}$ which, by Proposition~6.3, implies
$\rk\Gamma_{A_1\cup A_2}\leq  \rk(\Gamma_{A_1} \cap
\Gamma_{A_2})$. Hence, by Lemma 6.4 $\rk\Gamma \leq \rk
\Gamma_{A_1} + \rk\Gamma_{A_2} - \rk \Gamma_{A_1\cup A_2}$ and now
Proposition~6.3 and Lemma~6.4 i/ allow to conclude i/. ii/ follows
directly from Proposition~6.3 and Lemma~6.4 ii/. \qed

\sn{\bf Definition 6.6} Let us now look at a situation in which we
pick more hyperplanes, $\W = \{W_1,\cdots,W_f\}$, $f>m=\dim V$,
from $\C$ demanding that the set $\N(\W)$ of normals of the planes
(defined as in Definition 3.1 but for the subset $\W\subset\C$) is
indecomposable in the sense of 3.5. This requires that $\N(\W)$
spans $V$ and by (3.6) is equivalent to the fact that all graphs
$G(B;\N(\W))$ with $B\subset \N(\W)$ a basis for $V$ are
connected. Let us denote by $\I_l$ the collection of subsets
$A\subset \{1,\cdots,f\}$ of $m-l$ elements such that $W_A$ has
dimension $l$ (compare with 2.1). Note that any $A\in \I_0$
defines a basis $B_A$ by the normals to all $W_i$, $i\in A$.
\sn{\bf Theorem 6.7} {\it With $\W$ as above suppose that $\P$ is
finitely generated. Then, for all $A\in \I_l$ $$\rk\Gamma^A =
l\frac{\rk\Gamma}{\dim V}.$$ In particular $\dim V$ divides
$\rk\Gamma$.}

\bew\ Fix $i$ and choose an $A\in \I_0$ which contains $i$ (this
is clearly possible). Applying Corollary~6.5~i/ iteratively we
obtain $$\sum_{i\in A}\rk\Gamma^{\{i\}} = (m-1)\rk\Gamma.
\eqno(6.1)$$ By hypothesis there exists a $j\in A$ such that the
vertices of $G(B_A,\N(\W))$ corresponding to the normals of $W_i$
and $W_j$ are linked. Hence there is a $k \notin A$ such that the
normal of $W_k$ has nonvanishing scalar product $[\cdot,\cdot]$
with the normals of $W_i$ and $W_j$. This implies that, both,
$(A\backslash\{i\})\cup\{k\}$ and $(A\backslash\{j\})\cup\{k\}$
belong to $\I_0$. Hence we can apply (6.1) to both sets to
conclude $\nu_i = \nu_j$. By assumption $G(B_A;\N(\W))$ is
connected so that a repetition of the argument shows that $\nu_i$
does not depend on the choice of $i\in A$. This proves the theorem
for $l=m-1$ ($l=m$ is clear). The statements for $l<m-1$ follow
now inductively from Corollary~6.5~i/.\qed

\sn Theorem~6.7 contains as a special case a condition for
canonical projection method patterns as to whether they have
infinitely generated cohomology which can quickly be checked.

\sn{\bf Corollary 6.8} {\it Suppose that $\T$ is a tiling in
$\R^d$, topologically conjugate (I.4.5) to a canonical projection
method pattern with data $(E,u)$, and suppose that $E \cap \Z^N =
0$. If $N - \rk \Delta$ is not divisible by $N - \rk \Delta-d$
then $H_0(\G\T_u) \otimes \Q$ is infinite dimensional, and so
$\T_u$ is not a substitution tiling.}

\sn{\bf Proof} In the tiling case $\rk \Gamma_{\T_u} = N - \rk
\Delta$ and $\dim V = N-d-\rk \Delta$. \qed

\sn{\bf Examples 6.9} (Cf.\ Examples I.2.7) The Octagonal tiling,
a canonical projection tiling with $N=4$, $\Delta = 0$ and $d=2$,
and the Penrose tiling, a canonical projection tiling with $N=5$,
$\Delta \cong \Z$ and $d=2$, are both substitutional and hence
have finitely generated cohomology. Theorem 6.7 therefore tells us
that the stabilizers of the hyperplanes (here lines) have rank
two.

\sn It is clear that in the generic placement of planes the ranks
of the stabilizers of the intersections of hyperplanes will have
smaller rank than compatible with the last theorem. Thus we
deduce:

\sn{\bf Theorem 6.10} {\it Suppose that $\T$ is a tiling in
$\R^d$, topologically conjugate to a canonical projection method
pattern with data $(E,u)$ and $N > d+1$, and suppose that $E$ is
in generic position. Then $H_0(\G\T) \otimes \Q$ is infinite
dimensional, and $\T$ is not a substitution tiling. \qed}

\newpage

\headline={\ifnum\pageno=1\hfil\else
{\ifodd\pageno\rightheadline\else\leftheadline\fi}\fi}
\def\rightheadline{\tenrm\hfil{\smallfont
V COHOMOLOGY FOR SMALL CODIMENSION}\hfil\folio}
\def\leftheadline{\tenrm\folio\hfil{\smallfont
FORREST HUNTON KELLENDONK}\hfil} \voffset=2\baselineskip

\sn{\chap V Approaches to calculation III:}
\medskip
\noindent{\chap cohomology for small codimension}
\bigskip

\sn{\sect 1 Introduction}

\sn The last chapter was devoted to the case in which the
cohomology groups of a canonical projection tiling were infinitely
generated. Now we turn to the opposite case. In particular we
shall assume that $\P$ has only finitely many $\Gamma$ orbits,
i.e. is finitely generated (IV.2.8). In that case we find that the
cohomology groups of canonical projection tilings are finitely
generated free abelian groups and we can provide explicit formulae
for their ranks if the codimension is smaller or equal to $3$. We
obtain a formula for their Euler characteristic even in any
codimension. Although we saw that finitely generated $\P$ is
non-generic this seems to be the case of interest for quasicrystal
physics. In particular we will present calculations for the
Ammann-Kramer tiling as an example. It is a three dimensional
analog of the Penrose tilings and often used to model icosahedral
quasicrystals.

The material presented here extends
\FHKphys\ to the codimension $3$ case.
This is important, because all known tilings
which are used to model icosahedral quasicrystals are obtained from
projection out of a $6$-dimensional periodic structure with codimension $3$.
Unlike in \FHKphys\ we use here a spectral sequence derived from a double
complex to prove our formulae.

Perhaps the first use of spectral sequences in calculation of tiling
cohomology or $K$-theory is found in \bcl . However, we note that the spectral
sequence we use here differs significantly from the spectral sequence found in
\bcl. While the sequence of \bcl\ produces an isomorphism of the
$K$-theory of the groupoid $C^*$-algebra of the tiling
with its cohomology in $2$ dimensions (as is generalized in Chapter II), our
sequence represents a geometric decomposition of the tiling cohomology itself
which is a powerful calculating tool.

In Section~2 we recall the set-up and state the main results
(Thms.~2.4, 2.5, 2.8 for arbitrary codimension and Thm.~2.7 for
codimension $3$). In Section~3 we explain one part of the double
complex, namely the one which is related to the structure of the
set of singular points, and in Section~4 we recall the other half
which is simply group homology. We put both together in Section~5
where we employ the machinery of spectral sequences to proof our
results. In principle there is no obstruction to pushing the
calculation further to higher codimension, it is rather the
complexity which becomes overwhelming to do this in practice. In
the final section we sketch how our formulae apply to the
Ammann-Kramer tiling.

\newpage

\sn{\sect 2 Set up and statement of the results}

\sn In the second chapter we defined the cohomology of a tiling as
the cohomology of one of its groupoids (it turned out that it does
not matter which one we take) and interpreted it in various forms,
in particular as a standard dynamical invariant of one of the
systems. This dynamical system is of the following type: Consider
a dense lattice $\Gamma$ of rank $N$ in a euclidian space $V$ as
it arises generically if one takes $N>\dim V$ vectors of $V$ and
considers the lattice they generate. Let $K$ be a compact set
which is the closure of its interior and consider the orbit
$S=\partial K+\Gamma\subset V$ of $\partial K$ under $\Gamma$. We
showed in Chapters I and II how such a situation arises for
projection method patterns and how in this case the rather simple
dynamical system $(V,\Gamma)$ ($\Gamma$ acting by translation)
extends to a dynamical system $(\overline{V},\Gamma)$ which
coincides with the old one on the dense $G_\delta$-set
$V\backslash S$. $\overline{V}$ is locally a Cantor set and
obtained from $V$ upon disconnecting it along the points of $S$.
We are interested in calculating the homology groups
$H_p(\Gamma,C\overline{V})$ of the group $\Gamma$ with
coefficients in the compactly supported integer valued continuous
functions over $\overline{V}$. Already for the results of
Chapter~IV we specialized to the situation in which $S$ is a union
of a collection $\C$ of hyperplanes, the collection consisting of
finitely many $\Gamma_\T$-orbits. We denoted there $\overline{V}$
by $V_\C$. We will restrict our attention to this case here too.
In the context of projection method patterns (Def.\ I.4.4), where
$\Gamma=\Gamma_\T$ and $V=V+\pi'(u)$ (II.4.3), this means that we
consider only polytopal acceptance domain $K$ and such that the
orbit of a face under the action of $\Gamma$ contains the
hyperplane it spans (this is hypothesis H3 in \FHKphys). The
$\Gamma_\T$-orbits of these hyperplanes constitute the set $\C$.
In particular, we rule out fractal acceptance domain although this
case might be important for quasicrystal physics.

\sn {\bf Set-up 2.1} All our results below depend only on the
context described by the data $(V,\Gamma,\W)$, a dense lattice
$\Gamma$ of finite rank in a Euclidean space $V$ with a finite
family $\W=\{W_i\}_{i=1,\cdots,f}$ of (affine) hyperplanes whose
normals span $V$. Thus they are not specific for the tilings
considered here but can be applied e.g.\ also to the situation of
\FHKphys. Our aim is to analyse the homology groups
$H_p(\Gamma,CV_\C)$ where $\C:= \{W+\gamma : W\in \W,
\gamma\in\Gamma\}$.\bigskip

We extend the notation of IV.6.6.

\sn{\bf Definitions 2.2} Given a finite family
$\W=\{W_i\}_{i=1,\cdots,f}$ of hyperplanes, recall that $J_l$ is
the collection of subsets $A\subset  \{1,\cdots,f\}$ of $\dpe-l$
elements such that $W_A$ has dimension $l$. On $W_A$ we have an
action of $\Gamma^{|A|}$: $$\x\cdot W_A:= \bigcap_{i\in
A}(W_i-x_i).$$ We call $\x\cdot W_A$ a {\it singular space} or
singular $l$-space if we want to  specify its dimension $l$. Let
$\JJ_l$ be the quotient $\Gamma^{|A|}\times J_l/\sim$ with
$(\x,A)\sim (\x',A')$ if $\x\cdot W_A =\x'\cdot W_{A'}$. It is in
one to one correspondence to
 the set of singular $l$-spaces.
We denote equivalence classes by $[\x,A]$ and the space $\x\cdot
W_A$ also by $W_{[\vec{x},A]}$. On $\JJ_l$ we have an action of
$\Gamma$: $y\cdot [\vec{x},A] = [y\cdot\vec{x},A]$ where
$(y\cdot\vec{x})_i=y+x_i$. We denote the orbit space
$\JJ_l/\Gamma$ by $I_l$. Since this action coincides with the
geometric action of $\Gamma$ (by translation) on the singular
$l$-spaces we can use the elements of $I_l$ to label the
$\Gamma$-orbits of singular $l$-spaces. Note that the map
$\{1,\cdots,f\}\to I_{\dpe-1}$ which assigns to $i$ the orbit of
$[0,\{i\}]$ is surjective but not necessarily injective. The
stabilizer of a singular $l$-space depends only on its orbit
class, if the label of its orbit is $\Theta\in I_l$ we denote the
stabilizer by $\Gamma^\Theta$.

Fix $\hat\Theta\in \JJ_{l+k}$, $l+k<\dpe$ and let
$\JJ_l^{\hat\Theta}:=\{\hat\Psi\in\JJ_l|W_{\hat\Psi}\subset
W_{\hat\Theta}\}$. Then $\Gamma^\Theta$ ($\Theta$ the orbit class
of $\hat\Theta$) acts on $\JJ_l^{\hat\Theta}$ (diagonally, by the
same formula as above) and we let
$I_l^{\hat\Theta}=\JJ_l^{\hat\Theta}/\Gamma^\Theta$, the orbit
space. It labels the $\Gamma^\Theta$-orbits of singular $l$-spaces
in the $l+k$-dimensional space $W_{\hat\Theta}$. We can naturally
identify $I_l^{\hat\Theta}$ with $I_l^{\hat\Theta'}$ if
$\hat\Theta$ and $\hat\Theta'$ belong to the same $\Gamma$-orbit
and so we define $I_l^{\Theta}$, for the class $\Theta\in
I_{l+k}$. $I_l^{\Theta}$ is the subset of $I_l$ labelling those
orbits of singular $l$-spaces which have a representative that
lies in singular space whose label is $\Theta$. Finally we denote
$$ L_l = |I_l|,\quad L^\Theta_l = |I^\Theta_l| $$ where $\Theta\in
I_{l+k}$, $l+k<\dpe$.

Note that $\JJ_0$ can be identified with $\P$ (Def.~IV.2.8). So
our assumption of this chapter is that $L_0$ is finite. The proof
of the following lemma is straightforward.

\sn{\bf Lemma 2.3} {\it If $L_0$ is finite then $L^\Theta_l$ is
finite for all $\Theta$ and $l$.}\qed\bigskip

The proof of the following five theorems will be given in
Section~5. The first one is the converse of Theorem~IV.2.9.

\sn{\bf Theorem 2.4} {\it Given data $(V,\Gamma,\W)$ as in (2.1)
with $L_0$ finite. Then $H_p(\Gamma,CV_\C)\otimes\Q$ has finite
rank over the rational numbers $\Q$.}\bigskip

We denote $$ D_p = \rk H_p(\Gamma,CV_\C)\otimes\Q.$$

\sn{\bf Theorem 2.5} {\it Given data $(V,\Gamma,\W)$ as in (2.1)
with $L_0$ finite. Then $H_p(\Gamma,CV_\C)$ is free abelian. In
particular it is uniquely determined by its rank which is
$D_p$.}\bigskip

Note that Theorems~2.4 and 2.5 exclude the possibility that
$H(\Gamma,CV_\C)$ would contain e.g.\ the dyadic numbers as a
summand, a case which occurs frequently in cohomology groups of
substitution tilings.

Recall that, if the normals of the hyperplanes $\W$ form an
indecomposable set in the sense of IV.3.4 and $L_0$ is finite then
Theorem~IV.6.7 implies that the rank of the stabilizer
$\Gamma^\Theta$ depends only on the dimension of the plane it
stabilizes, i.e.\ $$\rk\Gamma^\Theta=\nu \dim\Theta$$ where $\nu =
\frac{\rk\Gamma}{\dim V}$ and $\dim\Theta=l$ provided $\Theta\in
I_l$.

Recall Theorem~III.3.1 applied to the present situation where we
have data $(V,\Gamma,\W)$ with $L_0$ finite and $\dim V=1$. Then $
H_p(\Gamma,CV_\C) = \Z^{D_p}$ with
 $$ D_p = {\left(\nu \atop p+1
\right)},\quad p>0, \eqno(2.1) $$ $$ D_0  = (\nu -1) + L_0.
\eqno(2.2) $$ If $\{M_i:i\in I\}$ is a family of submodules of
some bigger module we denote by $\erz{M_i:i\in I}$ their span. For
a finitely generated lattice $G$ we let $\Lambda G$ be the
exterior ring (which is a $\Z$-module) generated by it. In
\FHKphys\ we obtained the following theorem.

\sn{\bf Theorem 2.6} {\it Given data $(V,\Gamma,\W)$ as in (2.1)
with $\dim V=2$. Suppose that $L_0$ is finite and that the normals
of the hyperplanes $\W$ form an indecomposable set in the sense of
IV.3.4. Then $$ D_p = \left(2\nu \atop p+2 \right) + L_1\left(\nu
\atop p+1 \right)-r_{p+1}-r_p,\quad p>0, $$ $$ D_0 = \left(2\nu
\atop 2 \right)-2\nu+1 + L_1(\nu -1) + e -r_1 $$ where
$$r_p=\rk\erz{\Lambda_{p+1}\Gamma^{\alpha}:\alpha\in I_1},$$ and
the Euler characteristic is $$e:=\sum_{p} (-1)^p D_{p} =
-L_0+\sum_{\alpha\in I_1}L^\alpha_0.$$}

\sn The main result of this chapter is an extension of this result to
codimension $3$.

\sn{\bf Theorem 2.7} {\it Given data $(V,\Gamma,\W)$ as in (2.1)
with $\dim V=3$. Suppose that $L_0$ is finite and that the normals
of the hyperplanes $\W$ form an indecomposable set in the sense of
IV.3.4. Then, for $p>0$, $$ D_p
 = \left( 3\nu\atop p+3\right)
+L_2 \left( 2\nu\atop p+2\right) +\tilde L_1\left( \nu\atop
p+1\right) - R_p - R_{p+1}, $$ $$ D_0
 = \sum_{j=0}^3 (-1)^j
\left( 3\nu\atop 3-j\right)+ L_2\sum_{j=0}^2 (-1)^j \left(
2\nu\atop 2-j\right) $$ $$ + \tilde L_1\sum_{j=0}^1 (-1)^j
\left(\nu\atop 1-j\right) + e  - R_1 $$ where $\tilde L_1 =
-L_1+\sum_{\alpha\in I_2}L_1^\alpha$,
$$R_p=\rk\erz{\Lambda_{p+2}\Gamma^{\alpha}:\alpha\in I_2}
-\rk\erz{\Lambda_{p+1}\Gamma^{\Theta}:\Theta\in I_1}
+\sum_{\alpha\in I_2}
\rk\erz{\Lambda_{p+1}\Gamma^{\Theta}:\Theta\in I_1^\alpha} ,$$ and
the Euler characteristic is $$ e:=\sum_{p} (-1)^p D_{p} =
L_0-\sum_{\alpha\in I_2} L^\alpha_0 +\sum_{\alpha\in
I_2}\sum_{\Theta\in I_1^\alpha}L^\Theta_0 -\sum_{\Theta\in
I_1}L^\Theta_0. $$}

\sn The last two theorems make combinatorial patterns for higher
codimensional tilings apparent. We are able to present one such
for the Euler characteristic of the general case. To set this up
define a {\it singular sequence} to be a (finite) sequence $c =
\Theta_1, \Theta_2, ..., \Theta_k$ of $\Gamma$-orbits of singular
spaces strictly ascending  in the sense that $\Theta_j \in
I^{\Theta_{j+1}}_{\dim \Theta_j}$, $\dim \Theta_j<\dim
\Theta_{j+1}$, and $\dim \Theta_1=0$. An example of such a
sequence is the $\Gamma$-orbit of a singular flag, but the
dimension in the sequence can also jump by more than one. The
length of the chain $c$ is $k$, written $|c| = k$.

\sn{\bf Theorem 2.8} {\it Given data $(V,\Gamma,\W)$ with $L_0$
finite. Then the Euler Characteristic equals $$ e:=\sum_{p} (-1)^p
D_{p} = \sum (-1)^{|c|+\dim V}$$ where the sum is over all
singular chains $c$.}

\bigskip

\sn{\sect 3 Complexes defined by the singular spaces}

\sn Let $\C'$ be an arbitrary countable collection of affine hyperplanes
of $V'$, a linear space, and define $\C'$-topes as in IV.3.2:
compact polytopes which are the closure of their interior and
whose boundary faces belong to hyperplanes
from $\C'$.

\sn{\bf Definition 3.1} For $n$ at most the dimension of $V'$ let
$C_{\C'}^n$ be the $\Z$-module generated by the $n$-dimensional
faces of convex $\C'$-topes satisfying the relations
$$[U_1]+[U_2]=[U_1\cup U_2]$$ for any two $n$-dimensional faces
$U_1,U_2$, for which $ U_1\cup U_2$ is as well a face and $U_1\cap
U_2$ has no interior (i.e.\ nonzero codimension in $U_1$). (The
above relations then imply $[U_1]+[U_2]=[U_1\cup U_2]+[U_1\cap
U_2]$ if $U_1\cap U_2$ has interior.)

If we take $\C'=\C:=\{W+x:W\in\W,x\in\Gamma\}$, our collection of singular
planes, then $C^n:=C_\C^n$
carries an obvious $\Gamma$-action,
namely $x\cdot [U]=[U+x]$. It is therefore an $\Gamma$-module.
Recall Definition~IV.3.2
in which we defined $V_\C$ by identifying $CV_\C$
with $C^{\dpe}$ as $\Gamma$-module.
An isomorphism between
$C^{\dpe}$ and $CV_\C$ is  given by assigning to
$[U]$ the indicator function on the connected component containing
$U\backslash S$ (which is a clopen set). Moreover, $C^0$ is
a free $\Z$-module, its generators are in one to one correspondence
to the elements of $\P$.
The following proposition is from \FHKphys.

\sn{\bf Proposition 3.2} {\it There exist $\Gamma$-equivariant
module maps $\delta$ and $\epsilon$ such that $$ 0\to C^\dpe
\bra{\delta} C^{\dpe-1}\bra{\delta}\cdots C^0\bra{\epsilon}\Z\to 0
$$ is an exact sequence of $\Gamma$-modules and $\epsilon[U]=1$
for all vertices $U$ of $\C$-topes.}

\bew\ For a subset $R$ of $\Gamma$ let $\C_R:=\{W+r:W\in\W,r\in
R\}$ and $S_R=\{x\in W:W\in\C_R\}$. Let $\cal R$ be the set of
subsets $R\subset\Gamma$ such that all connected components of
$V\backslash S_R$ are bounded and have interior. $\cal R$ is
closed under union and hence forms an upper directed system under
inclusion. Let $V_R$ be the disjoint union of the closures of the
connected components of $V\backslash S_R$. $R\subset R'$ gives
rise to a natural surjection $V_{R'}\to V_R$ which is the
continuous extension of the inclusion $V\backslash S_{R'}\subset
V\backslash S_R$ and $V_\C$ is the projective limit of the $V_R$.
For any $R\in\cal R$, the $\C_R$-topes define a regular polytopal
CW-complex $$ 0\to C_{\C_R}^\dpe\bra{\delta_R}C_{\C_R}^{\dpe-1}
\bra{\delta_R}\cdots C_{\C_R}^0\to 0, \eqno(3.1) $$ with boundary
operators $\delta_R$ depending on the choices of orientations for
the $n$-cells ($n>0$) \massey. Moreover, this complex is acyclic
($V$ is contractible), i.e.\ upon replacing $C_{\C_R}^0\to 0$ by
$C_{\C_R}^0\bra{\epsilon_R}\Z\to 0$ where $\epsilon_R[U]=1$, (3.1)
becomes an exact sequence. Let us constrain the orientation of the
$n$-cells in the following way: Each $n$-cell belongs to a unique
singular $n$-space $W_{\hat\Theta}$, $\hat\Theta\in\JJ_n$. We
choose its orientation such that it depends only on the class of
$\hat\Theta$ in $I_l$ (i.e.\ we choose an orientation for all
parallel $W_{\hat\Theta}$ and then the cell inherits it as a
subset). By the same principle, all cells of maximal dimensions
are supposed to have the same orientation. Then the  cochains and
boundary operators $\delta_R$ share two crucial properties: first,
if $R\subset R'$ for $R,R'\in\cal R$, then we may identify
$C_{\C_R}^n$ with a submodule of $C_{\C_{R'}}^n$ and under this
identification $\delta_R(x)=\delta_{R'}(x)$ for all $x\in
C_{\C_R}^n$, and second,
if $U$ and $U+x$ are $\C_R$-topes then
$\delta_R[U+x]=\delta_R[U]+x$.
The first property implies that the directed system $\cal R$ gives rise
to a directed system of acyclic cochain complexes, and hence its
direct limit is an acyclic complex, and the second implies, together with
the fact that for all $x\in\Gamma$ and $R\in\cal R$ also $R+x\in
\cal R$, that this complex becomes a complex of $\Gamma$-modules.
The statement now follows since $C_{\C}^n$ is the direct limit of
the $C_{\C_R}^n$, $R\in\cal R$.\qed

\sn Let $C^l_{[\vec{x},A]}$ be the restriction of $C^l$ to
polytopes which belong to $W_{[\vec{x},A]}$. We can naturally
identify $C^l_{[\x,A]}$ with  $C^l_{[\x',A']}$ if
$y\cdot[\x,A]=[\x',A']$ for some $y\in\Gamma$. Moreover, in this
case the complexes obtained by restriction of $\delta$, $$0\to
C^{l}_{[\x,A]}\bra{\delta^{[\x,A]}} C^{l-1}_{[\x,A]} \cdots
C^{0}_{[\x,A]}\to 0,$$ are isomorphic for all $[\x,A]$ of the same
orbit class in $I_l$. With this in mind we define
$((C_\Theta^l)_{l\in\Z},\delta^\Theta)$ as a complex isomorphic to
the above complex where $\Theta\in I_l$ is the class of $[\x,A]$.
Thus every $\Theta\in I_l$ defines a new acyclic complex.

\sn{\bf Lemma 3.3}{\it $$C^l\cong \bigoplus_{\Theta\in I_l}
C^l_\Theta \otimes \Z[\Gamma/\Gamma^\Theta]$$ and if $\Theta\in
I_{l+k}$, $l+k<\dpe$, then $$C^l_\Theta\cong \bigoplus_{\Psi \in
I_l^{\Theta}} C^l_{\Psi}\otimes
\Z[\Gamma^\Theta/\Gamma^{\Psi}].$$} \bew\ First observe that
$$C^l\cong \bigoplus_{[\x,A]\in \JJ_l} C^l_{[\x,A]},$$ because an
$l$-face of a $\C$-tope belongs a unique singular $l$-space. From
the definition of $C_\Theta$ and the observation that
$y\cdot[\x,A]=[\x,A]$ whenever $y$ leaves $W_A$ invariant the
first statement of the lemma follows. The proof for the second is
similar.\qed

\bigskip

\sn{\sect 4 Group homology}

\sn Recall that
the homology of a discrete group $\Gamma'$ with coefficients in a
$\Gamma'$-module $M$
is defined as the homology of a complex which is obtained from
a resolution of $\Z$ by (projective) $\Gamma'$-modules upon application of
the functor $\otimes_{\Gamma'} M$. We use here a free
and finite resolution of $\Gamma'=\Gamma\cong \Z^N$ by its exterior module
so that the
complex looks like
$$0\to \Lambda_N\Gamma\otimes M\bra{\partial}\cdots
\Lambda_0\Gamma\otimes M\to 0$$
with boundary operator $\partial$ given by
$$\partial(e_{i_1}\cdots e_{i_k}\otimes m) =
\sum_{j=1}^k e_{i_1}\cdots\hat e_{i_j}\cdots
e_{i_k}(e_{i_j}\cdot m-m)$$
denoting by $\hat e_{i_j}$ that it has to be left out and
by $g\cdot$ the action of $g\in\Gamma$ on $M$.
In particular, $H_k(\Gamma,\Z\Gamma)$ is trivial for all $k>0$ and
isomorphic to $\Z$ for $k=0$.

Suppose that we can split $\Gamma=\Gamma'\oplus \Gamma''$ and let
us compute $H_k(\Gamma,\Z \Gamma'')$ where $\Z \Gamma''$ is the
free $\Z$-module generated by $\Gamma''$ which becomes an
$\Gamma$-module under the action of $\Gamma$ given by $(g\oplus
h)\cdot h' = h+h'$. Then we can identify $$ \Lambda_k\Gamma\otimes
\Z \Gamma'' \cong \bigoplus_{i+j=k} \Lambda_i \Gamma' \otimes
\Lambda_j \Gamma''\otimes\Z \Gamma'' \eqno(4.1) $$ and under this
identification $\partial\otimes 1$ becomes $(-1)^{deg}\otimes
\partial'$ where $\partial'$ is the boundary operator for the
homology of $\Gamma''$. It follows that $$H_k(\Gamma,\Z \Gamma'')
\cong \bigoplus_{i+j=k} \Lambda_i \Gamma' \otimes H_j( \Gamma'',\Z
\Gamma'') \cong \Lambda_k \Gamma'.$$ As a special case,
$H_k(\Gamma,\Z)\cong \Lambda_k\Gamma\cong \Z^{\left(N\atop
k\right)}$. Now let $\epsilon:\Z \Gamma''\to \Z$ be the sum of the
coefficients, i.e.\ $\epsilon[h]=1$ for all $h\in \Gamma''$. We
shall later need the following lemma: \sn{\bf Lemma 4.1}{\it Under
the identifications $H(\Gamma,\Z \Gamma'')\cong \Lambda \Gamma'$
and $H(\Gamma,\Z)\cong \Lambda \Gamma$ the induced map
$\epsilon_k:H_k(\Gamma,\Z \Gamma'')\to H_k(\Gamma,\Z)$ becomes the
inclusion $\Lambda_k \Gamma'\hookrightarrow\Lambda_k\Gamma$.}
\bew\ Using the decomposition (4.1) it is not difficult to see
that the induced map $\epsilon_k: \bigoplus_{i+j=k} \Lambda_i
\Gamma' \otimes H_j( \Gamma'',\Z \Gamma'')\to \bigoplus_{i+j=k}
\Lambda_i \Gamma' \otimes H_j( \Gamma'',\Z)$ preserves the
bidegree and must be the identity on the first factors of the
tensor product. Since $H_k( \Gamma'',\Z \Gamma'')$ is trivial
whenever $k\neq 0$ and one dimensional for $k=0$, $\epsilon_k$ can
be determined by evaluating $\epsilon_0$ on the generator of $H_0(
\Gamma'',\Z \Gamma'')$ and one readily checks that this gives a
generator of $H_0( \Gamma'',\Z)$ as well. \qed

\sn
Finally, we note two immediate corollaries of Lemma~3.3.
\sn{\bf Corollary 4.2} {\it
$$H(\Gamma,C^l)\cong \bigoplus_{\Theta\in I_{l}}
H(\Gamma^\Theta,C^l_\Theta).$$}\qed

\sn{\bf Corollary 4.3} {\it
If $\Theta\in I_{l+k}$, $l+k<\dpe$, then
$$H(\Gamma^\Theta,C^l_\Theta)\cong
\bigoplus_{\Psi\in I_l^\Theta}
H(\Gamma^{\Psi},C^l_{\Psi}).$$}\qed

\sn{\sect 5 The spectral sequences}

\sn In the last two sections we described two complexes. They can
be combined to yield a double complex and then spectral sequence
techniques can be applied to obtain the information we want. We do
not explain these techniques in detail but only set up the
notation, see, for example, \br\ for an introduction. Consider a
double complex $(E^0,\partial,\delta)$ which is a bigraded module
$E^0=(E^0_{pq})_{p,q\in \Z}$ with two commuting differential
operators $\partial$ and $\delta$ of bidegree $(-1,0)$ and
$(0,-1)$, respectively. The associated total complex is $TE^0 =
((TE^0)_k)_{k\in\Z}$ with $(TE^0)_k = \bigoplus_{p+q=k} E^0_{pq}$
and differential $\delta+\partial$. Furthermore one can form two
spectral sequences of bigraded modules \br. The first one,
$(E^k)_{k\in \Nat_0}$, is equipped with differential operators
$\df^k$ of bidegree $(k-1,-k)$, i.e.\ $$ \df^k_{pq} : E^k_{pq}
\longrightarrow E^k_{p\!+\!k\!-\!1\,q\!-\!k}, \eqno (5.1) $$ such
that $E^{k+1} = H_{\df^k}(E^k)$. (Here and below we use also the
notation $Z_{d}({C}^k)$ and $H_{d}({C}^k)$ for the degree $k$
cycles and homology resp.\ of a complex $(C,d)$.) $E^0$ is the
original module, $\df^0 = \partial$, $\df^1 = \delta_*$ and the
higher differentials become more and more subtle. The second
sequence, $(\tilde E^k)_{k\in \Nat_0}$, is obtained in the same
way except that one interchanges the role of the two
differentials, i.e.\ $\tilde E^0=E^0$, $\tilde{\df}^0 =\delta$,
$\tilde{\df}^1 = \partial_*$. The two spectral sequences may look
rather different but their "limits" are related to the homology of
the total complex. This relation may be quite subtle but in our
case we need only the following. Suppose that $E^0_{pq}$ is
non-trivial only for finitely many $p,q$. Then the higher
differentials vanish for both sequences from some $k$ on so that
the modules stabilize and we have well defined limit modules
$E^\infty$, $\tilde E^\infty$. If the modules are moreover vector
spaces then $$\bigoplus_{p+q=k} E^\infty_{pq} \cong
\bigoplus_{p+q=k} \tilde E^\infty_{pq} \cong
H_{\partial+\delta}((TE^0)_k).$$ If the modules are only over some
ring such as the integer numbers we may at least conclude that
finitely generated $E^\infty$ implies that $\tilde E^\infty$ and
$H_{\partial+\delta}(TE^0)$ are finitely generated as well.

Inserting the complex from Proposition~3.1 for $M$ in the last
section we get the double complex $(\Lambda_p\Gamma\otimes
C^q,\partial_{pq},\delta_{pq})$ where
$\partial_{pq}=\partial_p\otimes 1$ and $\delta_{pq} =
(-1)^p\otimes \delta_q$.

\sn{\bf Proposition 5.1} {\it Consider the spectral sequence
$\tilde E^k$ derived from the double complex
$(E_0^0=\Lambda_p\Gamma\otimes C^q,\partial_{pq},\delta_{pq})$
when starting with homology in $\delta$. It satisfies $$
\bigoplus_{p+q=k}\tilde E^\infty_{pq} \cong \Lambda_k\Gamma. $$ In
particular, $\rk \left(H_{\delta+\partial}(
\bigoplus_{p+q=k}\Lambda_p\Gamma\otimes C^q)\otimes\Q\right) =
\left( \rk \Gamma \atop k\right)$. } \bew\ The first page of that
spectral sequence is $$ \tilde E^1_{pq} =
H_\delta(\Lambda_p\Gamma\otimes C^q) \cong
\eos{\Lambda_p\Gamma}{q=0}{0}{q>0.}$$ Under the isomorphism
$\tilde E^1_{p0}\cong \Lambda_p\Gamma$, $\tilde{\df}^1$ becomes
trivial and all other differentials $\tilde{\df}^k$ vanish because
they have bidegree $(-k,k-1)$. Hence $\tilde E^\infty=\tilde E^1$
from which the first statement follows. The second statement is
then clear. \qed

 \sn
The more difficult spectral sequence will occupy us for the rest
of this chapter. Its first page is obtained when starting with
homology in $\partial$, i.e.\ $${E}^1_{pq} =
H_\partial(E^0_{pq})=H_p(\Gamma,C^q),$$ and the first differential
is $\df^1 = \delta_*$. We realize that ${E}^1_{p\,\dpe}$ is what
we want to compute. For that to carry out we have to determine the
ranks of the higher differentials which is quite involved. On the
other hand the proofs of Theorems~2.4, 2.5, and 2.8 involve only
the general fact that the higher differentials $\df^k$ have
bidegree $(k-1,-k)$ and begin with that. All three proofs are by
induction on the dimension of $V$ the one-dimensional case
following from Eqs.~(2.1),(2.2).

\sn {\bf Proof of Theorem 2.4} Theorem~2.4 is a rational result.
So we work with a rationalized version of the above. Corollary~4.2
yields $${E}^1_{pq}\cong \bigoplus_{\Theta\in
I_q}H_p(\Gamma^\Theta,C^q_\Theta).$$ Hence by Lemma~2 and
induction on the dimension of $V$ we find that
${E}^1_{pq}\otimes\Q$ is a finite dimensional vector space,
provided $q<\dim V$. This implies that ${E}^k_{pq}\otimes\Q$ is
finite dimensional, provided $q<\dim V$, for arbitrary $k$, and is
trivial for $q>\dim V$ or $p>\rk\Gamma$. On the other hand by
Proposition~5.1 $E^\infty_{pq}\otimes\Q$ has to be finite
dimensional for all $p,q$. By construction $E^\infty_{p\,\dpe}$ is
obtained upon taking successively the kernels of the higher
differentials. So if ${E}^1_{p\,\dpe}\otimes\Q$ was infinite, one
of the ranks of the higher differentials would have to be
infinite. But this would contradict the finite dimensionality of
${E}^k_{pq}\otimes\Q$ for $q<\dim V$.\qed

\sn {\bf Proof of Theorem 2.5} It is a result of \fh\ that
$H_p(\Gamma,CV_\C)$ is torsion free. Therefore Theorem~2.5 follows
if we can show that $H_p(\Gamma,CV_\C)$ is finitely generated. Now
we use Lemma~2.2 and Corollary~4.2 to see inductively that
${E}^1_{pq}$ is finitely generated, provided $q<\dim V$. This
implies that ${E}^k_{pq}$ is finitely generated, provided $q<\dim
V$, for arbitrary $k$. By Eq.~(5.1) we get sequences $$
E^{k+1}_{p\,\dpe}\hookrightarrow E^{k}_{p\,\dpe} \bra{\df^k}
E^{k}_{p+k-1\,\dpe-k}$$ which are exact at $E^{k}_{p\,\dpe}$ and
whose right hand side modules are finitely generated if $k\geq1$.
Hence $E^{k}_{p\,\dpe}$ is finitely generated if
$E^{k+1}_{p\,\dpe}$ is finitely generated. But
$E^{\infty}_{p\,\dpe}$ is finitely generated by Proposition~5.1.
Again inductively we conclude therefore that $E^{1}_{p\,\dpe}$ is
finitely generated.\qed

\sn {\bf Proof of Theorem 2.8} By Theorem~2.4 the dimensions of
the rational vector spaces $E^k_{pq}\otimes\Q$ are finite. From
Eq.~(5.1) it follows therefore that $\sum_{p,q}(-1)^{p+q}\dim
E^k_{pq}\otimes\Q$ is independent of $k$. Choosing $k=\infty$ we
see from Proposition~5.1 that this sum vanishes. With $e_q =
\sum_p (-1)^p\dim H_p(\Gamma,C^q)\otimes\Q$ we thus get $$e =
e_\dpe = -\sum_{q=0}^{\dpe-1}(-1)^{q+\dpe}e_q.$$ But $\dim
H_p(\Gamma,C^q)\otimes\Q = \sum_{\Theta\in I_q}\dim
H_p(\Gamma^\Theta,C_\Theta^q)\otimes\Q$ by Corollary~4.2 and we
can use induction on $\dpe$ to see that $e_q = \sum_{\Theta\in
I_q}\sum (-1)^{q+|c|-1}$ where the second sum is over all singular
sequences $c$ whose last element is $\Theta$. Inserting this into
the above formula for $e$ gives directly the result.\qed

\sn The strategy to prove the remaining Theorems~2.6 and 2.7 is as
follows. Suppose that we know $\d_{pq}=\dim {E}^1_{pq}\otimes\Q$
for $q<\dpe$. Then we can determine $\d_{p\,\dpe}$ from these
data, the dimensions of the rationalized total homology groups
(Proposition~5.1), and the ranks of the differentials. To carry
that out we have to consider only modules over $\Q$. We will
therefore simplify our notation by suppressing $\otimes\Q$ and
e.g.\ write $E^k_{pq}$ instead of ${E}^k_{pq}\otimes\Q$. We
consider first a more abstract situation.

\sn{\bf Definition 5.2} Consider a double complex
$(E^0_{pq},\delta_{pq},\partial_{pq})$ of finitely generated
$\Q$-modules with its spectral sequence $(E^k)_{k\in\Nat_0}$,
i.e.\ $E^1_{pq}:=H_\partial(E^0_{pq})$ and $\df^1=\delta_*$. We
write $\d_{pq} = \dim E^1_{pq}$. Assume that
\item{E1} There exist finite $M,N$ such that
$E^0_{pq}$ is non-trivial only for $0\leq q\leq M$ and $0\leq p \leq N$,
\item{E2} $\d_{p0}=0$ for $p \geq 1$.

\noindent
In particular the spectral sequence converges and we denote
$N_k:= \sum_{p+q=k} \dim E^\infty_{pq}$.

\sn For $M=1,2,3$, we now determine the ranks $\d_{pM}$ in terms
of $\d_{pq}$, $q<M$, and $N_k$ and the ranks of the higher
differentials $\df^k$.

\sn{\bf Lemma 5.3} {\it If $M=1$ then, with the notation of
Definition~5.2, $$ \d_{p1}=N_{p+1} +
\eos{\d_{00}-N_0}{p=0}{0}{p>0.} $$} \bew\ If {$M=1$} then (5.1)
implies that $\df^k=0$ for $k\geq 2$ and hence $E^\infty=E^2$.
With the notation $$a_{pq}:= \rk
\df^1_{p\,q\!+\!1}:E^1_{p\,q\!+\!1}\to E^1_{pq} \eqno (5.2)$$ we
get $$ \dim E^\infty_{pq} = \left\{\eqalign{ \d_{00}-a_{00} &
\quad \hbox{ for }\quad {p=q=0},\cr \d_{01}-a_{00} & \quad \hbox{
for }\quad {p=0,q=1},\cr \d_{p1} & \quad \hbox{ for }\quad {N\geq
p\geq 1,q=1},\cr 0 &\quad  \hbox{ otherwise, } }\right. $$ which
when $a_{00}$ is eliminated for $N_0$ yields the statement.\qed

\sn Let $$b_{pq}:= \rk
\tilde\delta_{p\,q\!+\!1\,*}:E^1_{p\,q\!+\!1}\to
H_\partial(Z_\delta(E^0_{pq})), \eqno(5.3)$$ where
$\tilde\delta_{p\,q\!+\!1\,*}$ is the map induced from
$\tilde\delta_{p\,q\!+\!1}$ which in turn is obtained from
$\delta_{p\,q\!+\!1}$ by restricting its target space, namely it
forms the exact sequence $0\to E^0_{p\,q+2}\bra{\delta_{p\,q+2}}
E^0_{p\,q+1}\bra{\tilde\delta_{p\,q+1}} Z_\delta(E^0_{pq})\to 0$.

\sn{\bf Lemma 5.4} {\it
If $M=2$ then,  with the notation of Definition~5.2,
$$
\d_{p2}=N_{p+2} + \d_{p1} - b_{p\!+\!1\,0}-b_{p0},\quad p>0
$$
$$
\d_{02}=N_{2} + \d_{01} - b_{10}-N_1+N_0-d_{00}.
$$}
\bew\ If $M=2$ then
(5.1) still implies that $\df^k=0$ for $k\geq 2$,
because of E2. But now, with the notation (5.2),
$$
\dim E^\infty_{pq} = \left\{\eqalign{
\d_{00}-a_{00} & \quad \hbox{ for }\quad {p=q=0},\cr
\d_{01}-a_{00}-a_{01} & \quad \hbox{ for }\quad {p=0,q=1},\cr
\d_{p1}-a_{p1} & \quad \hbox{ for }\quad {N\geq p\geq 1,q=1},\cr
\d_{p2}-a_{p1} & \quad \hbox{ for }\quad {N\geq p\geq 0,q=2},\cr
0 &\quad  \hbox{ otherwise. } }\right.
$$
Recall that $\df^1_{p\,q\!+\!1}:E^1_{p\,q\!+\!1}\to E^1_{pq}$
is the map induced on the homology groups from $\delta_{p\,q\!+\!1}$. From
the exact sequence
$0\to E^0_{p2}\bra{\delta_{p2}}
E^0_{p1}\bra{\tilde\delta_{p1}}
Z_\delta(E^0_{p0})\to 0$
we conclude
$$
a_{p1} = \dim\ker \tilde\delta_{p1\,*} = \d_{p1} - b_{p0}
$$
which now implies the statement.\qed

\sn{\bf Lemma 5.5} {\it
If $M=3$ then,  with the notation of Definition~5.2 and (5.3)
$$ \d_{p3} =
N_{p+3}+\d_{p2}+\d_{p\!+\!1\,1}-b_{p1}-b_{p\!+\!1\,1}-b_{p\!+\!1\,0}
-b_{p\!+\!2\,0},\quad p>0, $$
$$ \d_{03} = N_3-N_2+N_1-N_0+\d_{02}
+\d_{11}-\d_{01}+\d_{00}-b_{11}-b_{20}. $$} \bew\ Now (5.1)
implies that $\df^k=0$ only for $k>2$ so that $E^3=E^\infty$. With
the notation $$\mu_{pq}:= \rk
\df^2_{p\!-\!1\,q\!+\!2}:E^2_{p\!-\!1\,q\!+\!2}\to E^2_{pq}$$ and
(5.2) we get $$ \dim E^\infty_{pq} = \left\{\eqalign{
\d_{00}-a_{00} & \quad \hbox{ for }\quad {p=q=0},\cr
\d_{01}-a_{00}-a_{01} & \quad \hbox{ for }\quad {p=0,q=1},\cr
\d_{p1}-a_{p1}-\mu_{p1} & \quad \hbox{ for }\quad {N\geq p\geq
1,q=1},\cr \d_{p2}-a_{p1}-a_{p2} & \quad \hbox{ for }\quad {N\geq
p\geq 0,q=2},\cr \d_{p3}-a_{p2}-\mu_{p\!+\!1\,1} & \quad \hbox{
for }\quad {N> p\geq 0,q=3},\cr \d_{N3}-a_{N2} & \quad \hbox{ for
}\quad {p=N,q=3},\cr 0 &\quad  \hbox{ otherwise. } }\right. $$ Let
us recall the definition of $\df^2$. Look at the commuting diagram
$$ \matrix{ 0&\to &E^0_{p3}&\bra{\delta_{p3}}
&E^0_{p2}&\bra{\delta_{p2}}&E^0_{p1}& \bra{\delta_{p1}}&\cdots \cr
&&&&{\downarrow }{\tilde\delta_{p2}}&\buildrel{\imath_{p0}}\over{
\nearrow}& {\downarrow
}{\tilde\delta_{p1}}&\buildrel{\imath_{p0}}\over \nearrow& \cr
&&&& Z_\delta(E^0_{p1}) & & Z_\delta(E^0_{p0})&& \cr} . $$ Let
$\theta:H_\partial(Z_\delta(E^0_{p1})) \to
H_\partial(E^0_{p\!-\!1\,3})=E^1_{p\!-\!1\,3}$ be the connecting
map of the exact sequence $0\to E^0_{p3}\bra{\delta_{p3}}
E^0_{p2}\bra{\tilde\delta_{p2}} Z_\delta(E^0_{p1})\to 0$. Then
$\df^2_{p\!-\!1\,3}:E^2_{p\!-\!1\,3}\to E^2_{p1}$ is the map
induced by $\imath_{p1\,*}\circ\theta^{-1}:E^1_{p\!-\!1\,3}\to
E^1_{p1}$. Since $\theta^{-1}$ identifies
$E^2_{p\!-\!1\,3}=Z_{d^1}(E^1_{p\!-\!1\,3})
\buildrel{\theta^{-1}}\over{\cong} H_\partial(Z_\delta(E^0_{p1}))/
\im\tilde\delta_{p2\,*}$ we get $$\mu_{p1} = \rk\imath_{p1\,*} -
\rk\delta_{p2\,*} = \dim\ker \tilde\delta_{p1\,*} -
\rk\delta_{p2\,*}= (\d_{p1}-b_{p0}) - a_{p1}.$$ Furthermore, using
$a_{p2} = \dim\ker \tilde\delta_{p2\,*} = \d_{p2}-b_{p1}$ one
obtains the statement.\qed

\sn{\bf Definition 5.6}
We now specify to double complexes of the form
$E^0_{pq} = \Lambda_p\Gamma\otimes C^q$,
$\partial_{pq}=\partial_p\otimes 1$
and $\delta_{pq} = (-1)^p\otimes \delta_q$
as they arrise in our application
however tacitly changing the ring to be $\Q$.
Below, $\d_{pq}$ shall always denote the dimension of
$E^1_{pq}:=H_\delta (\Lambda_p\Gamma\otimes C^q)$.

\sn
Clearly, $E^0_{pq}$ is non-trivial only for
$0\leq p\leq \rk \Gamma$ and $0\leq q\leq \dpe$. Furthermore,
$E^1_{p0}=H_\delta (\Lambda_p\Gamma\otimes C^0)=H_p(\Gamma,C^0)$.
Since $C^0$ is always a free $\Gamma$-module one has
$H_p(\Gamma,C^0)=0$ for $p>0$ and the dimension of
$H_0(\Gamma,C^0)$ is simply $|I_0|$, the number of orbits of in $\P$:
$$
\d_{00} = L_0.
$$
In particular E1 and E2 are satisfied. Furthermore,
$N_k = {{\rk \Gamma} \choose k}$.

The case $\dpe = 1$ (2.1,2.2)
can now be immediately obtained using Lemma~5.3
and Proposition~5.1.

\sn {\bf Proof of Theorem~2.5} A proof can also be found in
\FHKphys\ but we repeat it here using the framework of spectral
sequences. So let $\dpe=2$. We seek to apply Lemma~5.4. For that
we need to determine $\d_{p1}$ and $b_{p0}$.

To determine $\d_{p1}$ we use Corollary~4.2.
$H(\Gamma^\alpha,C^1_\alpha)$ can be computed
using the double complex $(\Lambda_p\Gamma^\alpha\otimes C_\alpha^q,
\partial_{pq}^\alpha,\delta^\alpha_{pq})$.
So let us consider the spectral sequence which arrises when
starting with homology in $\partial^\alpha$. As for the case
$\dpe=1$ one concludes that $E^1_{p0}$ is non-trivial only for
$p=0$ and the dimension of $E^1_{p0}$ equal to $L^\alpha_0$, the
number of orbits of points in $\P\cap W_{\alpha}$. Since $\rk
\Gamma^\alpha=\nu$ we get $$\d_{p1} = L_1\left(\nu \atop p+1
\right),\quad p>0,$$ $$\d_{01} = L_1(\nu -1) + \sum_{\alpha\in
I_1}L^\alpha_0.$$ It remains to compute $ b_{p0} = \rk
\tilde\delta_{p1\,*}:H_p(\Gamma,C^1)\to
H_p(\Gamma,Z_\delta(C^0))$. For that look at the following
commuting diagram of exact sequences $$ \matrix{ 0 \to &
C_\alpha^1\otimes \Z[\Gamma/\Gamma^\alpha] &
\bra{\delta_1^\alpha\otimes 1} & C_\alpha^0\otimes
\Z[\Gamma/\Gamma^\alpha] & {\to} & \Z[\Gamma/\Gamma^\alpha] &\to 0
\cr & \downarrow \tilde\delta_{1}^\alpha\otimes 1 & & \downarrow &
& \downarrow \epsilon^\alpha & \cr 0 \to & Z_\delta(C^0) &
\hookrightarrow & C^0 & \bra{\epsilon} & \Z &\to 0 \cr} $$ where
the middle vertical arrow is the inclusion, the right vertical
arrow the sum of the coefficients, $\epsilon^\alpha[\gamma]=1$,
and the direct sum over all $\alpha\in I_1$ of the left vertical
arrows is $\tilde\delta_{p1}$. The above commutative diagram gives
rise to two long exact sequences of homology groups together with
vertical maps, all commuting, $(\tilde\delta_{p1}^\alpha\otimes
1)_*: H_p(\Gamma,C^1_\alpha\otimes\Z[\Gamma/\Gamma^\alpha])\to
H_p(\Gamma,Z_\delta(C^0))$ being one of them. Now we use that for
$p>0$, $H_p(\Gamma,C^0_\alpha\otimes\Z[\Gamma/\Gamma^\alpha])=0$
and $H_p(\Gamma,C^0)=0$ which implies
$H_p(\Gamma,C^1_\alpha\otimes\Z[\Gamma/\Gamma^\alpha])\cong
H_{p+1}(\Gamma,\Z[\Gamma/\Gamma^\alpha])\cong
H_{p+1}(\Gamma^\alpha,\Z)$ and $H_p(\Gamma,Z_\delta(C^0))\cong
H_{p+1}(\Gamma,\Z)$ and $$ (\tilde\delta_{p1}^\alpha\otimes
1)_*=\epsilon_{p+1}^\alpha.$$ By Lemma~4.1 the map
$\epsilon_{p}^\alpha$ can be identified with the inclusion
$\Lambda_p\Gamma^\alpha\hookrightarrow\Lambda_p\Gamma$. For $p>0$
therefore, $$ b_{p 0}=\rk \erz{\im
(\tilde\delta_{p1}^\alpha\otimes 1)_*: \alpha\in I_1}=\rk
\erz{\Lambda_{p+1}\Gamma^\alpha: \alpha\in I_1}. \eqno(5.4) $$
Direct application of Lemma~5.4 yields now the formulas for the
dimensions stated in Theorem~2.2. \qed

\sn{\bf Proof of Theorem~2.6} Let $\dpe=3$. To apply Lemma~5.5 we
need to determine $\d_{p1},\d_{p2},b_{p0},b_{p1}$. To calculate
$d_{p1}$ we use Corollary~4.2 and Lemma~5.3 to obtain, for
$\Theta\in I_1$, $$\dim H_p(\Gamma^\Theta,C^1_\Theta) = \left(
\nu\atop p+1\right) + \eos{L^\Theta_0-1}{p=0}{0}{p>0.} $$ Hence $$
\d_{p1} = L_1\left( \nu\atop p+1\right) +
\eos{-L_1+\sum_{\Theta\in I_1} L^\Theta_0}{p=0}{0}{p>0.} $$ To
calculate $d_{p2}$ we consider the complexes, $\alpha\in I_2$,
$$0\to C^{2}_\alpha\bra{\delta^\alpha} C^{1}_\alpha
\bra{\delta^\alpha} C^{0}_\alpha\to 0.$$ The appropriate double
complex to consider is therefore $(\Lambda_p\Gamma^\alpha\otimes
C_\alpha^q,\partial^\alpha_{pq}, \delta^\alpha_{pq})$. Repeating
the arguments of the proof of Theorem~2.2 (but now using
Corollary~4.3 in place of 4.2) we get $$ \dim H_p(\Gamma^\alpha,
C_\alpha^2) = \left( 2\nu\atop p+2\right) + L^\alpha_1\left(
\nu\atop p+1\right)-r_p^\alpha-r_{p+1}^\alpha \quad p>0, $$ $$
\dim H_0(\Gamma^\alpha, C_\alpha^2) = \left( 2\nu\atop 2\right)
-2\nu+1 + L^\alpha_1(\nu-1)-r^\alpha_1 -L^\alpha_0+\sum_{\Theta\in
I_1^\alpha}L^\Theta_0 $$ where $$r^\alpha_{p}
=\rk\erz{\Lambda_{p+1}\Gamma^{\Theta}:\Theta\in I_1^\alpha}.$$
Writing $\F=\sum_{\alpha\in I_2} L^\alpha_1$ we thus get $$
\d_{p2} = L_2\left( 2\nu\atop p+2\right) + \F\left( \nu\atop
p+1\right)-\sum_{\alpha\in I_2} (r_p^\alpha+r_{p+1}^\alpha),\quad
p>0, $$ $$ \d_{02} = L_2\left(\left( 2\nu\atop
p+2\right)-2\nu+1\right)+\F \left(\nu-1\right)+ \sum_{\alpha\in
I_2}\left(\sum_{\Theta\in I_1^\alpha}
L^\Theta_0-L^\alpha_0-r_1^\alpha\right). $$ We now compute
$b_{p0}$ and $b_{p1}$. For the first we can again repeat the
arguments which yield (5.4) to obtain
$$b_{p0}=\rk\erz{\Lambda_{p+1}\Gamma^{\Theta}:\Theta\in I_1}.$$
$b_{p1}$ is the rank of $\tilde\delta_{p2\,*}$, i.e.\ the rank of
the map induced from $\tilde\delta_2:C^2\to Z_\delta(C^1)$. To
determine it we consider the commutative diagram with exact rows
$$ \matrix{ 0 \to & C_\alpha^{2}\otimes \Z[\Gamma/\Gamma^\alpha] &
\bra{\delta_2^\alpha\otimes 1} & C_\alpha^{1}\otimes
\Z[\Gamma/\Gamma^\alpha] & \bra{\tilde\delta_1^\alpha\otimes 1} &
Z_\delta (C_\alpha^{0})\otimes \Z[\Gamma/\Gamma^\alpha] &\to 0 \cr
& \downarrow \tilde\delta_2^\alpha\otimes 1 & & \downarrow & &
\downarrow  & \cr 0 \to & Z_\delta(C^{1}) & \hookrightarrow &
C^{1} & \bra{\delta_1} & Z_\delta(C^{0}) &\to 0\cr} . $$ The
middle and right vertical maps are the obvious inclusions. By
Corollary~4.2 the direct sum of the left vertical arrows of the
diagrams which arrise if $\alpha\in I_2$ is $\tilde\delta_2$. The
diagram gives rise to a commutative diagram of which the following
is one degree $$\matrix{
\ho{p}{C_\alpha^2\otimes\Z[\Gamma/\Gamma^\alpha]} & \to &
\ho{p}{C_\alpha^1\otimes\Z[\Gamma/\Gamma^\alpha]} &
\bra{(\tilde\delta_{p1}^\alpha\otimes 1)_*} &
\ho{p}{Z_\delta(C_\alpha^0)\otimes\Z[\Gamma/\Gamma^\alpha]} \cr
\downarrow (\tilde\delta_{p2}^\alpha\otimes 1)_* & & \downarrow &
& \downarrow \cr \ho{p}{Z_\delta(C^1)} & \to & \ho{p}{C^1} &
\bra{\delta_{p1\,*}} & \ho{p}{Z_\delta(C^0)} \cr} \eqno(5.5) $$
also with exact sequences. In particular, $b_{p1} = \rk\erz{\im
(\tilde\delta_{p2}^\alpha\otimes 1)_* :\alpha\in I_2}$. Now we can
follow the same analysis as in the proof of Theorem~2.2 to see
that, for $p>0$, (5.5) is isomorphic to $$\matrix{ {\to} &
\ho{p}{C_\alpha^2\otimes\Z[\Gamma/\Gamma^\alpha]} & \to &
\bigoplus_{\Theta\in I_1^\alpha} \Lambda_{p+1}\Gamma^{\Theta} &
\bra{\beta^\alpha_p} & \Lambda_{p+1}\Gamma^{\alpha} & {\to}\cr
&\downarrow (\tilde\delta_{p2}^\alpha\otimes 1)_* & & \downarrow
j^\alpha_{p} & & \downarrow \imath^\alpha_p&\cr {\to} &
\ho{p}{Z_\delta(C^1)} & \to & \bigoplus_{\Theta\in I_1}
\Lambda_{p+1}\Gamma^{\Theta} & \bra{\gamma_p} &
\Lambda_{p+1}\Gamma & {\to}\cr} \eqno(5.6) $$ where
$\beta_p^\alpha$ is the direct sum over $\Theta\in I_1^\alpha$ of
the inclusions $\Lambda_{p+1}\Gamma^{\Theta} \hookrightarrow
\Lambda_{p+1}\Gamma^{\alpha}$, $\gamma_p$ is the direct sum over
$\Theta\in I_1$ of the inclusions $\Lambda_{p+1}\Gamma^{\Theta}
\hookrightarrow \Lambda_{p+1}\Gamma$, and $j_p^\alpha$ and
$\imath^\alpha_p$ are the obvious inclusions. As we are working
with finitely generated $\Q$-modules we obtain from (5.{6}) $$ \im
(\tilde\delta_{p2}^\alpha\otimes 1)_*  \cong  \im
\imath^\alpha_{p+1}/ (\im \imath^\alpha_{p+1}\cap \im
\gamma_{p+1}) \oplus j^\alpha_p(\ker \beta^\alpha_p). $$ Since
$I_1=\bigcup_{\alpha_\in I_2}I_1^\alpha$ the direct sum over
$\alpha\in I_1^\alpha$ of $j^\alpha_p$ is surjective. Hence $\im
\gamma_p\subset \erz{\im \imath^\alpha_p:\alpha\in I_2}$ and
therefore $$ \erz{\im \imath^\alpha_{p}/ (\im
\imath^\alpha_{p}\cap \im \gamma_{p}):\alpha\in I_2} = \erz{\im
\imath^\alpha_p:\alpha\in I_2}/\im \gamma_p. $$ Furthermore, using
also that $\imath^\alpha_p$ is injective, we get $j^\alpha_p(\ker
\beta^\alpha_p)=\im j^\alpha_p\cap\ker \gamma_p$ and $$ \erz{\im
j^\alpha_p\cap\ker \gamma_p:\alpha\in I_2}= \ker \gamma_p. $$ Thus
we have $$ \erz{\im(\tilde\delta_{p2}^\alpha\otimes 1)_*:\alpha\in
I_2} \cong \erz{\im \imath^\alpha_{p+1}:\alpha\in I_2}/ \im
\gamma_{p+1}\oplus\ker \gamma_p $$ which, since
$\rk\gamma_p=b_{p0}$,  implies $$b_{p1} =
\rk\erz{\Lambda_{p+2}\Gamma^{\alpha}:\alpha\in I_2}
-b_{p\!+\!1\,0}+L_1\left(\nu\atop p+1\right)- b_{p0}.$$ With these
results we obtain from Lemma~5.5 the ranks stated in Theorem~2.3.
\qed
\bigskip

\sn{\sect 6 Example: Ammann-Kramer tilings}

\sn In Chapter IV we showed that generically, canonical projection
tilings have infinitely generated cohomology. Nevertheless, all
tilings known to us which are used to describe quasicrystals have
finitely generated cohomology. In \GaKe\ a list of results for the
known $2$-dimensional tilings used by quasicrystallographers was
presented. Here we present the first  result for a $3$-dimensional
projection tiling with finitely generated cohomology. It is the
Ammann-Kramer tiling. The tiling was invented before the discovery
of quasicrystals \kn\ and rediscovered in \DuneauKatz\ \Elser . It
has been used to describe icosahedral quasicrystals in
\ElserHenley\ \HenleyElser. It is sometimes also called
$3$-dimensional Penrose tiling because it generalizes in a way
Penrose's $2$-dimensional tilings. A short description of this
tiling and further $3$-dimensional tilings related icosahedral
symmetry is given in \KramerPapadopolos. The tiling is obtained by
the canonical projection method from the data $(\Z^6,E)$, $\Z^6$
being the lattice generated by an orthonormal basis of $\R^6$, and
$E$ being obtained from symmetry considerations involving the
representation theory of the icosahedral group. Projecting $\Z^6$
orthogonally onto the orthocomplement $V$ of $E$ one obtaines the
(dense) lattice $\Gamma$ generated by the six vectors
$$\left(\matrix{\tau \cr 1 \cr 0} \right),
 \left(\matrix{ -\tau \cr 1 \cr 0 } \right) ,
 \left(\matrix{ 0 \cr \tau \cr -1 } \right) ,
 \left(\matrix{ 0 \cr \tau \cr 1 } \right) ,
 \left(\matrix{ 1 \cr 0 \cr -\tau } \right) ,
 \left(\matrix{ 1 \cr 0 \cr \tau } \right).$$
with respect to an orthonormal basis of $V$.
Here $\tau = {{\sqrt 5 - 1}\over 2}$
and we have the usual relation $\tau^2+\tau-1=0$.
The acceptance domain is the orthogonal
projection of the $6$-dimensional unit cube into $V$ and
forms the triacontrahedron whose vertices are the
linear combinations with coefficients $0,1$ of the above vectors.
It has $30$ faces which are all triangles
and the affine hyperplanes $\W$ are the planes spanned by these
triangles. Although we have an action of the icosahedral group
(which even permutes the faces)
the  geometry of their intersections is quite complex, but
can be summarized as follows.

\noindent {\bf Singular subspaces in $V$:}
The two triangles which are in opposite
position of the triacontrahedron
belong to the same orbit (under the action of $\Gamma$)
and representatives of the $15$
orbits of singular planes are given by the spans of all possible
pairs of vectors from our list
of $6$ vectors above. Thus $I_2$ has $L_2=15$ elements.
The $15$ planes can be gathered into five orthogonal triples.
Pairs of these planes intersect in singular lines and triples intersect in
singular points.
Singular lines occur in three classes: Intersections of $5$ planes
simultaneously at a line (class I - $6$ possible orientations):
Intersections of $3$ planes only (class II - $15$ possible orientations):
Intersections of $2$ planes only (class III - $10$ possible orientations).
Each of the class I and III orientations is a single   orbit, but each
of the class II orientations splits into two   orbits. Thus if we write
$I_1 = I_1^{\rm  I}\cup I_1^{\rm  II}\cup I_1^{\rm  III}$,
a disjoint union according to
the class, we have $|I_1^{\rm  I}| = 6$, $|I_1^{\rm  II}| =
2\times 15$, and
$|I_1^{\rm  III}| = 10$. In particular,  we have $L_1=46$
orbits of singular lines, arranged in $31$ possible orientations.
A similar calculation finds $L_0=32$ orbits of singular points formed at
the intersection of three or more singular planes: Briefly, every singular
point lies at the intersection of some orthogonal triple of planes. Given a
particular othogonal triple, those points which are at the intersection of
three planes in this orientation split into $8$ orbits (because
intersections of orthogonal planes give class II lines). However, two of
those orbits are found at intersections of planes parallel to every
other orientation (e.g.\ the orbit of the origin), and the remaining $6$ are
unique to the particular orthogonal triple chosen. This adds up to $6
\times 5 +2 = 32$ orbits therefore.

\noindent {\bf Singular subspaces in $W_i$, $i\in \JJ_2$:}
Since all $W_i$ lie in one orbit of the semidirect product of
$\pi^\perp(\Z^6)$ with the icosahedral group it doesn't
matter which $i$ to take.
Within a singular plane $W_i$ there are $2$ directions of each class of
lines, but each of the class II lines is in two orbits.
Thus if we denote the class of $i$ in $I_2$ by $\alpha$ and decompose
$I_1^\alpha = I_1^{\alpha\,{\rm  I}}\cup I_1^{\alpha\,{\rm  II}}\cup
I_1^{{\alpha\,\rm  III}}$ disjointly according to
the classes of its singular lines we have $|I_1^{\alpha\,{\rm  I}}| = 2$,
$|I_1^{{\alpha\,\rm  II}}| = 2\times 2$,
and $|I_1^{\alpha\,{\rm  III}}| = 2$, giving a total
of $L_1^\alpha=8$
  orbits of singular lines in $6$ possible directions.
In particular $\tilde L_1 = 15\times 8 - 46 = 74$.
A careful
calculation finds the singular points arranged into $L_0^\alpha=8$ orbits.

\noindent {\bf Singular subspaces in $W_{\hat\Theta}$,
$\hat\Theta\in \JJ_1$:} In each class I and class III line we have
exactly two orbits of singular points. In each class II line we
have exactly $4$ orbits. Hence $L_0^\Theta = 2$ if $\Theta$, the
class of $\hat\Theta$ in $I_1$, is a class I or class III line
whereas otherwise $L_0^\Theta = 4$. This gives $\sum_{\Theta\in
I_1} L_0^\Theta = 2\times 6 + 4\times 30 + 2\times 10 = 152$ and
$\sum_{\Theta\in I_1^\alpha} L_0^\Theta = 2\times 2 + 4\times 4 +
2\times 2 = 24$. Altogether $e=120$.

To find the ranks $R_1,R_2$, a more detailed analysis of the generators of the
stabilizers and their inner products must be made. Since this involves
finding the rank of matrices of up to $15 \times 15$ in size, it is best
checked by computer: $R_1=69$ and $R_2=9$.

Specializing Theorem 2.3 to $\dim V = 3$ and $\rk \Gamma = 6$ we obtain
as the only nonzero ranks
$$ D_3  =  1 $$
$$ D_2  =  6 + L_2 - R_2 $$
$$ D_1  = 15 + 4 L_2 + \tilde L_1 - R_1 - R_{2} $$
$$ D_{0} = 10 + 3 L_2 + \tilde L_1 + e  - R_1. $$

This yields the following numbers for the ranks of the homology groups:
$$ D_{0} =  180,\quad D_1  =  71,\quad  D_2   =  12, \quad D_3   =  1.$$ So to
conclude we have: $$K_0(C^*(\G\T)) \cong \Z^{192} \ \ \ \ \ \ \ K_1(C^*(\G\T))
\cong \Z^{72}$$ for the Ammann-Kramer tiling $\T$.
We are grateful to Franz G\"ahler for assistence in
completing the last step on the computer
and for verifying these results generally.

\newpage

\headline={\ifnum\pageno=1\hfil\else
{\ifodd\pageno\rightheadline\else\leftheadline\fi}\fi}
\def\rightheadline{\tenrm\hfil{}\hfil\folio}
\def\leftheadline{\tenrm\folio\hfil{\smallfont
FORREST HUNTON KELLENDONK}\hfil} \voffset=2\baselineskip

\sn{\sect References}{\smallfont
\bigskip

\item \ap J.E.Anderson, I.F.Putnam. Topological invariants for substitution
tilings and their associated C*-algebras.
Ergodic Theory and Dynamical Systems 18 (1998) 509--537.

\item \bks  M. Baake, R. Klitzing and M. Schlottmann,
Fractally shaped acceptance domains of quasiperiodic square-triangle tilings
with dodecagonal symmetry,
Physica A 191 (1992) 554-558.

\item \bjks M. Baake, D. Joseph, P. Kramer and M. Schlottmann,
Root lattices and quasicrystals,
J. Phys. A 23 (1990) L1037-L1041.

\item \bsj M.Baake, M.Schlottman and P.D.Jarvis. Quasiperiodic tilings
with tenfold symmetry and equivalence with respect to local derivability.
J.Phys.A 24, 4637-4654 (1991).

\item \beltwo J.Bellissard. K-theory for C*algebras in solid state
physics. Lect. Notes Phys. 257. Statistical Mechanics and Field Theory,
Mathematical Aspects, 1986, pp. 99-256.

\item \belone J.Bellissard. Gap labelling theorems for Schr\"odinger
operators. from Number Theory and Physics. Luck, Moussa, Waldschmit, eds.
Springer, 1993.

\item \BBG J. Bellissard, A. Bovier and J. M. Ghez. Gap labelling theorems for
one dimensional discrete Schr\"odinger operators, Rev. Math. Phys. 4 (1992),
1-38.

\item \bcl J.Bellissard, E.Contensou, A.Legrand. K-th\'eorie des
quasicristaux, image par la trace, le cas du r\'eseau octagonal.
C.R.Acad.Sci. Paris, t.326, S\'erie I, p.197-200, 1998.

\item \bhz J.~Bellissard, M.~Zarrouati and D.J.L.~Herrmann.
Hull of aperiodic solids and gap labelling theorems.
In: Directions in Mathematical Quasicrystals, eds.\ M. Baake and R.V. Moody,
CRM Monograph Series, Amer. Math. Soc. Providence, RI (2000).

\item{\br}
K.S. Brown.
Cohomology of Groups. Springer-Verlag, 1982.

\item \db N.G.de Bruijn. Algebraic Theory of Penrose's nonperiodic tilings
of the plane. Kon.Nederl.Akad.Wetensch. Proc. Ser.A 84 (Indagationes Math.
43) (1981), 38-66.

\item \dbtwo N.G.de Bruijn. Quasicrystals and their Fourier transforms.
Kon.Nederl.Akad.Wetensch.Proc. Ser.A 89 (Indagationes Math. 48) (1986),
123-152.

\item \con A.Connes. Non-commutative Geometry. Academic Press, 1994.

\item{\DuneauKatz}
M.~Duneau and A.~Katz. Phys. Rev. Lett., 54:2688, 1985.

\item \du F.Durand. PhD Thesis, 1997, Marseille.

\item{\Elser}
V.\ Elser. Acta Cryst., A42:36, 1986.

\item{\ElserHenley}
V\ Elser and Henley. Phys. Rev. Lett., 55:2883, 1985.

\item \fortwo A.H.Forrest. A Bratteli diagram for commuting
homeomorphisms of the Cantor set. Preprint 15.97 NTNU Trondheim.

\item{\fh} A.H.Forrest, J.Hunton. Cohomology and K-theory of Commuting
Homeomorphisms of the Cantor Set. Ergodic Theory and Dynamical Systems 19
(1999), 611-625.

\item{\fhk} A.H.Forrest, J.R.Hunton, J.Kellendonk. Projection
Quasicrystals I, II, and III. SFB-288-preprints, nos.\ 340, 396, and 459,
Berlin, August 1998, June  1999, and March 2000.

\item{\FHKphys} A.H.Forrest, J.Hunton, J.Kellendonk.
Cohomology of canonical projection tilings.
SFB-preprint No.\ 395, June 1999.

\item \furst H.Furstenberg. The structure of distal flows. Amer.J.Math. 85
(1963), 477-515.

\item{\GaKe}
F.\ G\"ahler and J.\ Kellendonk.
Cohomology groups for projection tilings of codimension $2$.
SFB-preprint No.\ 406, August 1999. To appear in the
Proceedings of the International conference on Quasicrystals 7.

\item \gljj  C. Godr\`eche, J.-M. Luck, A. Janner and T. Janssen,
Fractal atomic surfaces of self-similar quasiperiodic tilings of the plane.
J. Physique I  3 (1993) 1921.

\item{} J.-M. Luck, C. Godr\`eche, A. Janner and T. Janssen,
The nature of the atomic surfaces of quasiperiodic self-similar structures.
J. Phys. A26  (1993) 1951.

\item{\gps} T. Giordano, I. Putnam and C. Skau.
Topological orbit equivalence and C*-crossed products.
J. reine und angewandte Math., 469 (1995) 51-111.

\item \gs B.Gr\"unbaum, G.Shephard. Tilings and Patterns, San Francisco.
W.Freeman, 1987.

\item{\hm} G.A.Hedlund, M.Morse. Symbolic Dynamics II. Sturmian
trajectories. Am.J.Math.62 (1940), 1-42.

\item{\HenleyElser}
Henley and V.\ Elser. Phil. Mag., B53:L59, 1986.

\item{\her}
D.J.L.~Herrmann. Quasicrystals and Denjoy homeomorphisms
Preprint, 1999.

\item \hof A.Hof. Diffraction by Aperiodic Structures. in The Mathematics
of Long Range Aperiodic Order, 239-268. Kluwer 1997.

\item \kd A.Katz, M.Duneau. Quasiperiodic patterns and icosahedral
symmetry. Journal de Physique, Vol.47, 181-96. (1986).

\item \kelzero  J. Kellendonk.
 Non Commutative Geometry of Tilings and Gap Labelling.
 Rev. Math. Phys. 7 (1995) 1133-1180.

\item \kelone J. Kellendonk. The local structure of Tilings and their
Integer group of coinvariants.
Commun. Math. Phys. 187 (1997) 115-157.

\item \keltwo J. Kellendonk. Topological Equivalence of Tilings.
J.Math.Phys., Vol 38, No4, April 1997.

\item  \kelput
J. Kellendonk and I.F. Putnam. Tilings, C*-algebras and K-theory.
In: Directions in Mathematical Quasicrystals, eds.\ M. Baake and R.V. Moody,
CRM Monograph Series, Amer. Math. Soc. Providence, RI (2000).

\item \kn P. Kramer and R. Neri. On Periodic and Non-periodic Space Fillings
of $E^m$ Obtained by Projection. Acta Cryst. A 40 (1984) 580-7.

\item{\KramerPapadopolos}
P.~Kramer and Z.~Papadopolos. Models of icosahedral
quasicrystals from 6D lattices.
In: Proceedings of the International conference on Aperiodic Crystals,
Aperiodic '94, eds.\ G.Chapuis et al.,  pages 70--76. World Scientific,
Singapore 1995.

\item \lag J.C.Lagarias. Meyer's concept of quasicrystal and quasiregular sets.
Comm.Math.Phys. 179 (1996) 365-376.

\item \lagtwo J.C.Lagarias. Geometric Models for Quasicrystals I.
Delone Sets of Finite Type. To appear in Discrete \& Computational Geometry.

\item \le T.T.Q.Le. Local Rules for Quasiperiodic Tilings. in The
Mathematics of Long Range Aperiodic Order, 331-366. Kluwer 1997.

\item \mack A.Mackay. Crystallography and the Penrose pattern. Physics 114A
(1982), 609-613.

\item{\massey}
W.S. Massey.
Singular Homology Theory. Springer-Verlag, 1980.

\item \mrw P.S.Muhly, J.N.Renault, D.P.Williams. Equivalence and
isomorphism for groupoid C* algebras. J.Operator Theory 17 (1987) 3-22.

\item \okd C.Oguey, A.Katz, M.Duneau. A geometrical approach to
quasiperiodic tilings.
Commun. Math. Phys. 118 (1988) 99-118.

\item \pt W.Parry, S.Tuncel. Classification Problems in Ergodic Theory. LMS
Lecture Notes series. Vol 65. Cambridge University Press 1982.

\item \pat A.Paterson. Groupoids, inverse semigroups and their operator algebras.
Birkh\"auser 2000.

\item \pen R.Penrose. Pentaplexity. Mathematical Intelligencer. 2(1979),
32-37.

\item{\pr} N.M.Priebe. Detecting Hierarchy in Tiling Dynamical Systems via
Derived Voronoi Tessellations. PhD Thesis, University of North Carolina at
Chapel Hill, 1997.

\item{\pss} I.F.Putnam, K.Schmidt, C.Skau, C$^*$-algebras
associated with Denjoy homeomorphisms of the circle. J. of
Operator Theory 16 (1986), 99-126.

\item \rad C.Radin, M.Wolff. Space tilings and local isomorphisms.
Geometricae Dedicata 42, 355-360.

\item \rei M.A.Rieffel Applications of strong Morita equivalence to
transformation group C*-algebras. Proc.Symp. Pure Math 38 (1) (1982) 299-310.

\item \ren J.Renault. A Groupoid Approach to C*algebras. Lecture Notes in
Mathematics 793. Springer-Verlag 1980.

\item \robone E.A.Robinson Jr.. The dynamical theory of tilings and
quasicrystallography. in Ergodic Theory of Zd actions (M.Pollicott and
K.Schmidt, eds.) 451-473,  Cambridge University Press, 1996.

\item \robtwo E.A.Robinson Jr.. The dynamical properties of Penrose
tilings. Trans.Amer.Math.Soc. 348 (1994), 4447-4469.

\item \rudolf D.J.Rudolf. Markov tilings of $R^n$ and representations of
$R^n$ actions. Contemporary Mathematics 94, (1989) 271-290.

\item \schtwo M.Schlottmann. Generalized model sets and dynamical systems.
In: Directions in Mathematical Quasicrystals, eds.\ M. Baake and R.V. Moody,
CRM Monograph Series, Amer. Math. Soc. Providence, RI (2000).

\item \sen M.Senechal. Quasicrystals and Geometry.
Cambridge University Press 1995.

\item \st M.Senechal, J.Taylor. Quasicrystals: The view from Les Houches.
from The Mathematical Intelligencer, Vol 12, No 2, 1990, pp.54-64.

\item \schecht D. Schechtman, L. Blech, D. Gratias, J.V. Cahn. Metallic phase
with long range orientational order and no translational symmetry.
Phys. Rev. Lett., 53 (1984) 1951-1953.

\item \smith A.P. Smith,
The sphere packing problem in quasicrystals.
J. Non-Cryst. Solids 153\&154 (1993), 258.

\item \soc J.E.S.Socolar. Simple octagonal and dodecagonal quasicrystals.
Phys. Rev. B, 39(15), 10519-10551.

\item \solone B.Solomyak. Dynamics of selfsimilar tilings.
Ergodic Theory and Dynamical Systems 17 (1997), 695-738.

\item \sol B.Solomyak. Spectrum of Dynamical Systems Arising from Delone
Sets. in Quasicrystals and discrete geometry, ed J.Patera, Fields
Institute Monographs, vol.10, Amer.Math.Soc., Providence, RI,
1998, pp.247-264.

\item{\cbt} C. B. Thomas, Characteristic classes and the
cohomology of finite groups, Cambridge studies in advanced
mathematics, Cambridge University Press, 1986.

\item\wil R.F. Williams. The Penrose, Ammann, and DA tiling spaces are
Cantor set fiber bundles. Preprint 2000.
\item \zobetz E. Zobetz,
A pentagonal quasiperiodic tiling with fractal acceptance domain.
Acta Cryst. A48 (1992) 328.

\item \colone "The Physics of Quasicrystals", editors P.J. Steinhardt
and S. Ostlund,
World Scientific Publishing Co. 1987.

\item \coltwo  Proceedings of the NATO Advanced Study Institute on
"The Mathematics of Long-Range Aperiodic Order", editor R.V. Moody,
Kluwer Academic Publishers 1997
}

\end